\documentclass[12pt]{article}
\usepackage{amsfonts}
\usepackage{amsmath}
\usepackage{latexsym}
\usepackage{epsfig}
\renewcommand{\thepage}{\hfill \arabic{page} \hfill}

\newtheorem{theorem}{Theorem}[section]
\newtheorem{lemma}{Lemma}[section]

\newcommand{\newsection}[1]{\setcounter{equation}{0} \setcounter{theorem}{0}
\section{#1}}
\renewcommand{\theequation}{\thesection .\arabic{equation}}
\renewcommand{\thetheorem}{\thesection .\arabic{theorem}}
\renewcommand{\baselinestretch}{1.3}

\def \spec#1 {\mathop{#1}}
\def\til#1{\displaystyle{\spec #1 _{\sim}}}
\def\top #1{\displaystyle{\spec #1 ^{\sim}}}
\def \b #1 {\bf #1}
\def \qed {\hfill \vrule height6pt width 6pt depth 0pt}

\begin{document}
\begin{center}
{\bf \Large {The structure of finite clusters in high intensity Poisson
 Boolean stick process}}\\
$\;$
\\
{\sc Rahul Roy\footnote{RR is grateful to Chiba University for its warm
hospitality and acknowledges the financial support of JSPS.
 }
 and Hideki Tanemura}
\\
{\it Indian Statistical Institute} and {\it Chiba University}
\end{center}

\vspace{1.5cm}

\noindent {\footnotesize {\bf Abstract} 
Sticks at one of different orientation are placed in an
i.i.d. fashion at points of a Poisson point process of intensity $\lambda$.
Sticks of the same direction have the same length, while sticks in
different directions may have different lengths.
We study the geometry of finite cluster as $\lambda \to \infty$.
The asymptotic shape of the custer being determined by the 
probabilities of the sticks in various direction and their lengths
and orientations.
We also obtain the limiting geometric structure of this component.
}

\vspace{1.5cm}

\newsection{Introduction}
Consider one dimensional sticks placed at random locations and with random
orientations in the two dimensional plane. In the language of stochastic
geometry we have a planar  fibre process whose {\it grains} are two
dimensional linear segments and whose germs are the random locations. The most
commonly studied fibre process model which incorporates these features is when
the germs arise as realisations of a Poisson point process of intensity
$\lambda$ on $\Bbb R^2$ and each germ is the centre of a
stick of either fixed length or a random length and having a random
orientation, with the distribution of the length and orientation of a
stick being independent of the underlying Poisson process. This is the
Poisson Boolean stick process, a particular instance of the more general
planar Boolean fibre process.  Hall [1990] (Chapter 4),
Stoyan Kendall and Mecke [1995] (Chapter 9)
discuss the geometric and statistical aspects of this process.

While the stochastic geometry study of these processes was motivated by its
application in geology, viz., the subterranean earthquake faults are
modelled as a Poisson Boolean stick process (see, e.g., Weber [1977]); the 
interest
in the physics community of this model led to a probabilistic study of its
percolative properties. Suppose mirrors are placed
randomly on the plane and we are interested in the path of a ray of light
in this set-up. Clearly the geometry of the mirrors on the
plane determine the trajectory of the ray of light. This model is
a modern equivalent of the Ehrenfest wind-tree model  which was introduced by
Ehrenfest [1957] to study the Lorentz lattice gas model (see Grimmett
[1998] for an exposition of the mathematical study of this model).
This model has also been studied for its percolative properties (in particular,
the critical phenomenon it exhibits and the corresponding critical parameters) 
by Domany
and Kinzel [1984], Hall [1985], Menshikov [1986], Roy [1991] and Harris
[1997].

Here we study the geometric features of finite clusters
in the Poisson Boolean stick process when the intensity
of the underlying Poisson process is high. More particularly, consider a Poisson 
point
process of intensity $\lambda$ on $\Bbb R^2$ conditioned to have a point
at the origin. At each
point $x_i$ we centre a stick of length $r_i$ and orientation $\theta_i$
measured anticlockwise w.r.t. the $x$-axis. We suppose that\\
(i) $r_1, r_2, \ldots $ is an i.i.d. sequence of random variables,\\
(ii) $\theta_1, \theta_2, \ldots $ is an i.i.d. sequence of random
variables, and\\
(iii) the sequences $\{r_i\}$ and $\{\theta_i\}$ and the underlying Poisson
process are independent of each other.\\
Consider the cluster of the origin (which is
the connected component formed by sticks containing the stick at
the origin).
For the above model Hall [1985] has shown that if  the random variable $r_1$
is bounded, and the random variable $\theta_1$ is non-degenerate then there
is a critical intensity $\lambda_c$ such that, for $\lambda > \lambda_c$,
with positive probability the cluster defined above is unbounded. Moreover
this probability goes to $1$ as $\lambda \rightarrow \infty$. Given the rare 
event
that
this cluster contains exactly $m$ sticks, we investigate its structure as
the intensity $\lambda \rightarrow \infty$.

In the case of the Boolean model which consists of an underlying Poisson
point process of intensity $\lambda$ on $\Bbb R^d$ and each
point of the process is the centre of a $d$-dimensional ball of radius $r$,
Alexander [1991] showed that conditional on the cluster of the origin (i.e.
the connected component of balls containing a ball which
covers the origin) being finite and consisting of $m$ balls, the event
that these balls are centred in a small region of radius $O (\lambda^{-1})$
has a probability which tends to $1$
as $\lambda \rightarrow \infty$. This region where the balls are centred has 
volume
$O(\lambda ^{-d})$ whereas the ambient density is $\lambda$, thereby
giving rise to the phenomenon of compression wherein many more Poisson
points are accomodated in this region than the ambient density allows.
Sarkar [1998] showed that in case the balls forming the Boolean model are
allowed to be of varying sizes, then given that the cluster of the origin
contains $m$ balls, not all of the same size, the phenomenon of
rarefaction occurs, wherein the biggest sized balls remain compressed in a very
small region, but the other balls are sparsely placed in the region covered
by the biggest sized balls.

In our model the phenomenon of compression also occurs, however that is
of secondary interest. Instead we look at the geometry and the distribution
of the sticks of various orientation in the finite cluster.

In this paper we restrict ourselves to the study of the model when the sticks 
have exactly two or three possible orientations and sticks of the same
orientation have the same length.
In the case of two possible orientations 
the asymptotic distribution was shown to be independent of 
the angle and the length of the sticks 
-- a result which is not surprising in view of the affine invariance 
of the model. 
However, if three or more orientations are allowed 
then the affine invariance breaks down and the 
asymptotic distribution do depend on the angles. 
In this case we show that the asymptotic shape consists of sticks 
with only two orientations. 
The orientations which ``survive'' are chosen 
according to the lengths and angles of the possible orientations
and the probabilities of the sticks in various directions. 

The paper is organised as follows:-- in the next section we present a formal 
defintion of the process as well as the statements of our results and in 
Sections 3 and 4 we prove the results.

\newsection{Preliminaries and statement of results}
\subsection{Notation}
Let ${\cal R}= \Bbb R ^2 \times [0,\pi ) \times (0,\infty )$,
and
$$
{\cal M} ={\cal M}({\cal R})
:=\left\{ \xi = \{ \xi_i, i \in \Bbb N \} :
\xi_i = ( x_i, \theta _i, r_i )\in {\cal R} \right \}.
$$
For $(x,\theta,r) \in {\cal R}$,
$S(x, \theta , r) =  \{ x + u e_{\theta}, u \in [-r, r]\}$
is the stick with centre $x$, angle
$\theta $ and length $2r$,
where $e_{\theta} = (\cos \theta , \sin \theta )$.
We define the collection of sticks for
$\xi \in {\cal M}$ as
$S(\xi )=\{ S(x, \theta , r) : (x, \theta , r) \in \xi \}$.

We say two sticks $S$ and $S^\prime $ are connected and write
$S \stackrel {\xi }{\leftrightarrow} S^\prime $ if
there exists $S_1, S_2, \ldots S_k \in S (\xi )$ such that
$S \cap S_1 \neq \emptyset , ~ S^\prime \cap S_k \neq \emptyset$ and
$ S_i \cap S_{i+1} \neq \emptyset$ for every $i=1,2,\ldots , k-1$.
If $S(\xi )$ contains a stick $S_{\bf 0}$ centred at the origin ${\bf 0}$, 
we denote by
$C_{\bf 0} (\xi )$ the cluster of sticks containing $S_{\bf 0}$, i.e.
$$
C_{\bf 0}(\xi )=\{ y \in S : S \in S(\xi ),
S\stackrel {\xi }{\leftrightarrow }S_{\bf 0}  \}. 
$$
(We put $C_{\bf 0} (\xi ) = \emptyset $, if $S (\xi )$ does not
contain any stick with centre ${\bf 0}$).

Let $\rho $ be the Radon measure on ${\cal R}$ defined by
\begin{equation}
\label{ldirection}
\rho (dx d\theta dr )
=dx \sum_{j=1}^{d} p_j \delta_{\alpha_j}(d\theta) \delta_{R_j}(dr),
\end{equation}
where $\alpha_1 =0< \alpha_2 <\alpha_3 <\dots <\alpha_d <\pi$,
$p_j \ge 0, \sum_{j=1}^{d} p_j =1$, $R_j>0$, $j=1,2,\dots,d$
and
$\delta_*$ denotes the usual Dirac delta measure.
We denote by $\mu _{\rho} $
the Poisson point process on ${\cal M} ({\cal R})$ 
with intensity measure $\rho $.
Let
\begin{equation}
\label{gamma0}
\Gamma_0 := \{\xi \in {\cal {M}}: 
 ({\bf 0}, \alpha_j, R_j)\in \xi \mbox{ for some }j =1,2,\dots,d\}.
\end{equation}
For $w_i =(x_i, \theta_i, r_i)$, $i=1,2,\dots,m$,
let 
\begin{equation}
\label{p-def}
{ \bf w}_m := (w_1, w_2, \dots, w_m), \{ {\bf w}_m \} := \{ w_1, w_2, \dots, w_m \}, 
C_{\bf 0}({\bf w}_m):= C_{\bf 0}(\{{\bf w}_m\}).
\end{equation}

For ${\bf k}=(k_1, k_2,\dots, k_d)\in ({\Bbb N} \cup \{0\})^{d}$, 
we denote by $\Lambda ({\bf k})$ the set of clusters
containing exactly $|{\bf k}|=\sum_j^{d} k_j$ sticks
with $k_j$ sticks at an orientation $\alpha_j$,
$j=1,2,\dots,d$.

For $\alpha, \beta > 0$, $R_{\alpha},R_{\beta}>0$, $e_{\alpha}= (\cos \alpha, \sin \alpha)$,
 and
${{\bf x}}_{m} =(x_1, x_2, \cdots, x_{m})\in (\Bbb R^2)^m$, we define the following regions:-
\begin{eqnarray*}
& &B^{\alpha, \beta} _{R_{\alpha},R_{\beta}} 
:=\{ x^\alpha e_{\alpha} + x^\beta e_{\beta}: 
(x^\alpha,x^\beta) \in [-R_\alpha, R_\alpha] \times [-R_\beta , R_\beta ] \},
\\
& &B^{\alpha, \beta} _{R_{\alpha},R_{\beta}} (x)
:=B^{\alpha, \beta} _{R_{\alpha},R_{\beta}} +x,
\quad
x\in \Bbb R^2,
\\
& &B^{\alpha, \beta} _{R_{\alpha},R_{\beta}} ({\bf x}_m)
:=\bigcup_{j=1}^{m} B^{\alpha, \beta} _{R_{\alpha},R_{\beta}} (x_j).
\end{eqnarray*}

\subsection {Sticks of two types}

In this subsection we assume that \\
(i) {\it there are sticks with only two orientations, and }\\
(ii) {\it  sticks of the same orientation are of the same length
but sticks along different directions could be of different lengths.}\\
\noindent
Without loss of generality we assume that sticks are either horizontal or
at an angle $\alpha \in (0, \pi ]$. 
Sticks which are horizontal are of length $R_0$ and
sticks at an angle $\alpha$ are of length $R_\alpha$.

In this case $\Lambda(k,\ell)$ is the set of clusters containing 
$k$ horizontal sticks and $\ell$ sticks at an angle $\alpha$
with respect to the $x$-axis.
We show that
\begin{theorem}
\label{th2sticks}
Let $m = k+\ell,~ k,\ell\ge 1, ~ \alpha \in (0,\pi )$ and
$ 0<R_0, R_\alpha$.  As $\lambda \rightarrow \infty$, we have

\noindent {\rm (i)}  \qquad
$\mu _{\lambda \rho } (C_0 \in \Lambda(k,\ell) \mid \Gamma_0 )$
$$
\sim \left ( \frac{1}{\lambda |B^{0,\alpha}_{R_{0},R_{\alpha }}|}\right )^{m-3}
e^{-\lambda |B^{0,\alpha }_{R_{0}, R_{\alpha }}|} 
(pq)^{-2(m-1)} mp^{3k} k! q^{3l} \ell!, 
$$
where 
$a(\lambda)\sim b(\lambda)$ means that 
$\frac{a(\lambda)}{b(\lambda)}\to 1$ as $\lambda \to \infty$;

\noindent {\rm (ii)} \qquad
$p_{\lambda,m}(k,\ell) := \mu _{\lambda \rho } (\# C_0 = (k,\ell) \mid \# C_0 =
(k^\prime , \ell^\prime ), ~k^\prime + \ell^\prime = m)$
$$
\sim \frac {p^{3k} k! q^{3\ell} \ell!}{\mathop{\sum} _{k+\ell=m}
p^{3k}k!q^{3\ell} \ell!}.
$$
\end{theorem}

An interesting observation from (ii) above is that
asymptotically, as $\lambda \rightarrow \infty$,
the conditional probability $p_{\lambda,m}(k,\ell)$ of the sticks 
comprising the {\it finite}\/ cluster $C_0$,  
is independent of both the angle $\alpha$ as well as
$R_0$ and $R_\alpha$, the lengths of the sticks.
This is not surprising because the model is invariant under
affine transformations.
Now let $p_m(k,\ell) := \lim_{\lambda \rightarrow \infty}
p_{\lambda, m}(k,\ell)$.
We also observe from Theorem \ref{th2sticks} (ii) that, 
as $m \rightarrow \infty$,
\begin{eqnarray*}
p_m(m-1, 1) \rightarrow 1 \qquad \qquad & \mbox{ for } p > q,\\
p_m(1, m-1) \rightarrow 1 \qquad \qquad & \mbox{ for } p < q,\\
p_m(1, m-1) = p_m(m-1, 1)  \rightarrow \frac{1}{2} & \mbox{ for } p = q.
\end{eqnarray*}
Moreover, let $k$ and $m$ both approach infinity in such a way that
$(k/m) \rightarrow s$, for some $s \in [0,1]$, then we have
$$
\mathop{\lim}_{m \rightarrow \infty \atop (k/m) \rightarrow s}
\frac{1}{m} \log p_m(k,\ell) = H(s),
$$
where
$$
H(s) = s \log s +  (1-s )\log (1 - s ) +
\begin{cases}
3 (1-s) \log (q/p),
\quad &\hbox{if } p > q,
\\
3 s \log (p/q),
\quad &\hbox{if } p < q,
\\
0, \quad &\hbox{if } p=q,
\end{cases}
$$
from which we may deduce that as $ m \rightarrow \infty$,
for $0 \leq a \leq b \leq 1$,\\
\noindent $P(\mbox{the proportion $(k/m)$ of horizontal sticks in the cluster lies
between $a$ and $b$})$\\
$\sim \exp\{\sup_{s \in (a,b)}H(s)\}$.

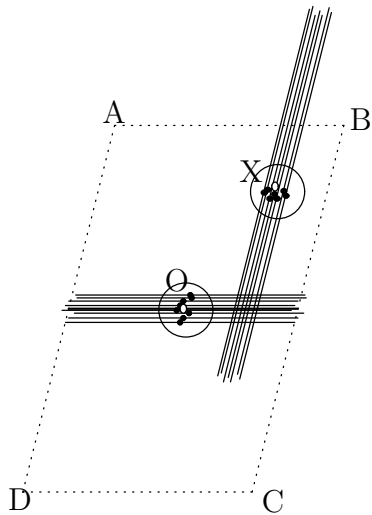
\begin{figure}[h]
\label{Fig1}
\begin{center}
\unitlength 0.1in
\begin{picture}( 17.8800, 25.4400)(  0.7000,-29.6400)
%
\special{pn 8}%
\special{pa 628 1044}%
\special{pa 1828 1044}%
\special{pa 1348 2964}%
\special{pa 148 2964}%
\special{pa 148 2964}%
\special{pa 628 1044}%
\special{dt 0.045}%
%
\special{pn 8}%
\special{pa 1750 428}%
\special{pa 1270 2348}%
\special{fp}%
%
\special{pn 13}%
\special{sh 1}%
\special{ar 1510 1388 10 10 0  6.28318530717959E+0000}%
\special{sh 1}%
\special{ar 1516 1388 10 10 0  6.28318530717959E+0000}%
%
\special{pn 8}%
\special{pa 424 1932}%
\special{pa 1624 1932}%
\special{fp}%
%
\special{pn 13}%
\special{sh 1}%
\special{ar 1024 1932 10 10 0  6.28318530717959E+0000}%
\special{sh 1}%
\special{ar 1024 1932 10 10 0  6.28318530717959E+0000}%
%
\special{pn 8}%
\special{pa 352 2012}%
\special{pa 1552 2012}%
\special{fp}%
%
\special{pn 13}%
\special{sh 1}%
\special{ar 952 2012 10 10 0  6.28318530717959E+0000}%
\special{sh 1}%
\special{ar 952 2012 10 10 0  6.28318530717959E+0000}%
%
\special{pn 8}%
\special{pa 370 1988}%
\special{pa 1570 1988}%
\special{fp}%
%
\special{pn 13}%
\special{sh 1}%
\special{ar 970 1988 10 10 0  6.28318530717959E+0000}%
\special{sh 1}%
\special{ar 970 1988 10 10 0  6.28318530717959E+0000}%
%
\special{pn 8}%
\special{pa 388 2052}%
\special{pa 1588 2052}%
\special{fp}%
%
\special{pn 13}%
\special{sh 1}%
\special{ar 988 2052 10 10 0  6.28318530717959E+0000}%
\special{sh 1}%
\special{ar 988 2052 10 10 0  6.28318530717959E+0000}%
%
\special{pn 8}%
\special{pa 388 1964}%
\special{pa 1588 1964}%
\special{fp}%
%
\special{pn 13}%
\special{sh 1}%
\special{ar 988 1964 10 10 0  6.28318530717959E+0000}%
\special{sh 1}%
\special{ar 988 1964 10 10 0  6.28318530717959E+0000}%
%
\special{pn 8}%
\special{pa 418 2028}%
\special{pa 1618 2028}%
\special{fp}%
%
\special{pn 13}%
\special{sh 1}%
\special{ar 1018 2028 10 10 0  6.28318530717959E+0000}%
\special{sh 1}%
\special{ar 1018 2028 10 10 0  6.28318530717959E+0000}%
%
\special{pn 8}%
\special{pa 370 2076}%
\special{pa 1570 2076}%
\special{fp}%
%
\special{pn 13}%
\special{sh 1}%
\special{ar 970 2076 10 10 0  6.28318530717959E+0000}%
\special{sh 1}%
\special{ar 970 2076 10 10 0  6.28318530717959E+0000}%
%
\special{pn 8}%
\special{pa 430 1948}%
\special{pa 1630 1948}%
\special{fp}%
%
\special{pn 13}%
\special{sh 1}%
\special{ar 1030 1948 10 10 0  6.28318530717959E+0000}%
\special{sh 1}%
\special{ar 1030 1948 10 10 0  6.28318530717959E+0000}%
%
\special{pn 8}%
\special{pa 1762 452}%
\special{pa 1282 2372}%
\special{fp}%
%
\special{pn 13}%
\special{sh 1}%
\special{ar 1522 1412 10 10 0  6.28318530717959E+0000}%
\special{sh 1}%
\special{ar 1528 1412 10 10 0  6.28318530717959E+0000}%
%
\special{pn 8}%
\special{pa 1714 468}%
\special{pa 1234 2388}%
\special{fp}%
%
\special{pn 13}%
\special{sh 1}%
\special{ar 1474 1428 10 10 0  6.28318530717959E+0000}%
\special{sh 1}%
\special{ar 1480 1428 10 10 0  6.28318530717959E+0000}%
%
\special{pn 8}%
\special{pa 1678 468}%
\special{pa 1198 2388}%
\special{fp}%
%
\special{pn 13}%
\special{sh 1}%
\special{ar 1438 1428 10 10 0  6.28318530717959E+0000}%
\special{sh 1}%
\special{ar 1444 1428 10 10 0  6.28318530717959E+0000}%
%
\special{pn 8}%
\special{pa 1702 444}%
\special{pa 1222 2364}%
\special{fp}%
%
\special{pn 13}%
\special{sh 1}%
\special{ar 1462 1404 10 10 0  6.28318530717959E+0000}%
\special{sh 1}%
\special{ar 1468 1404 10 10 0  6.28318530717959E+0000}%
%
\special{pn 8}%
\special{pa 1666 420}%
\special{pa 1186 2340}%
\special{fp}%
%
\special{pn 13}%
\special{sh 1}%
\special{ar 1426 1380 10 10 0  6.28318530717959E+0000}%
\special{sh 1}%
\special{ar 1432 1380 10 10 0  6.28318530717959E+0000}%
%
\special{pn 8}%
\special{pa 1648 436}%
\special{pa 1168 2356}%
\special{fp}%
%
\special{pn 13}%
\special{sh 1}%
\special{ar 1408 1396 10 10 0  6.28318530717959E+0000}%
\special{sh 1}%
\special{ar 1414 1396 10 10 0  6.28318530717959E+0000}%
%
\special{pn 20}%
\special{sh 1}%
\special{ar 988 2004 10 10 0  6.28318530717959E+0000}%
\special{sh 1}%
\special{ar 988 2004 10 10 0  6.28318530717959E+0000}%
%
\special{pn 20}%
\special{sh 1}%
\special{ar 1468 1364 10 10 0  6.28318530717959E+0000}%
\special{sh 1}%
\special{ar 1468 1364 10 10 0  6.28318530717959E+0000}%
%
\special{pn 13}%
\special{pa 388 2004}%
\special{pa 1588 2004}%
\special{fp}%
\put(5.6200,-10.2800){\makebox(0,0)[lb]{A}}%
\put(18.5800,-10.6800){\makebox(0,0)[lb]{B}}%
\put(13.9600,-30.6800){\makebox(0,0)[lb]{C}}%
\put(0.7000,-30.6000){\makebox(0,0)[lb]{D}}%
\put(8.9000,-19.1000){\makebox(0,0)[lb]{O}}%
\put(12.8000,-13.4000){\makebox(0,0)[lb]{X}}%
%
\special{pn 8}%
\special{sh 0}%
\special{ar 1468 1364 18 24  0.0000000 6.2831853}%
%
\special{pn 8}%
\special{ar 988 2004 18 24  0.0000000 6.2831853}%
%
\special{pn 8}%
\special{sh 0}%
\special{ar 988 2004 18 24  0.0000000 6.2831853}%
%
\special{pn 8}%
\special{ar 1480 1390 142 142  0.0000000 6.2831853}%
%
\special{pn 8}%
\special{ar 1000 2010 142 142  0.0000000 6.2831853}%
\end{picture}%
\end{center}
\caption{ \it{
The finite cluster for large $\lambda$. 
The region $X$ which contains the centres of the sticks 
at an angle $\alpha$ w.r.t. the 
$x$-axis is uniformly distributed in the parallelogram $ABCD$.}}
\end{figure}

>From the proof of the above theorem we also observe that
the centres of the horizontal sticks comprising the cluster $C_0$ lie in a
neighbourhood whose  area is of the order $o(\lambda^{-1+(\delta / 2)})$.
Similarly the centres of the sticks of orientation $\alpha$
comprising the cluster $C_0$ lie in another
neighbourhood whose  area is of the order $o(\lambda^{-1+(\delta /2)})$.
(See Figure 1.)
\subsection{Sticks of three types}
In this subsection we assume that \\
(i) {\it there are sticks with only three orientations --  
$0,\alpha$ and  $\beta$,}\\
(ii) {\it  sticks of the same orientation are of the same length.}\\
Here the results are significantly different from those obtained in the previous section. In particular the absence of any affine invariance leads to the dependence of the results on both the length and orientation of the sticks through the following quantities
\begin{equation}
\label{H}
H_{\alpha} = \frac{R_\alpha}{\sin \beta},\qquad
H_{\beta} = \frac{R_{\beta}}{\sin \alpha},\qquad
H_0 = \frac{R_0}{\sin (\beta-\alpha)}.
\end{equation}
By a suitable scaling we take 
\begin{equation}
\label{scalingH}
H_0=1 \mbox{ and let }H_\alpha = a, \;\;H_\beta = b \mbox{ after the scaling.}
\end{equation}
As the following theorem exhibits, 
the asymptotic (as $\lambda \to \infty$) composition of the finite cluster 
contains sticks of only two distinct orientation, 
while the third does not figure at all.
Here we use the shorthand ``$A(x,y)$ occurs" to mean that 
as $\lambda \rightarrow \infty$ 
the asymptotic shape of $C_0$ consists of sticks 
only in the directions $x$ and $y$.


\begin{theorem}
\label{pshape}
Given that $C_0$ consists of $m$ sticks,
\begin{itemize}
\item[{\rm (1)}] for $a,b \geq 2$;
\begin{itemize}
\item[\rm (i)] 
if $(ab -a +1/4)p_\beta +a < (ab-b+1/4)p_\alpha + b$, then 
$A(0,\alpha)$ occurs,
\item[{\rm (ii)}] 
if $(ab -a +1/4)p_\beta +a > (ab-b+1/4)p_\alpha + b$, then 
$A(0,\beta)$ occurs, and
\item[{\rm (iii)}]  
if $(ab -a +1/4)p_\beta +a = (ab-b+1/4)p_\alpha + b$, then 
both $A(0,\alpha)$ and $A(0,\beta)$ have positive 
probabilities of occurrence;
\end{itemize}

\item[{\rm (2)}]
for $1/2 < \min\{a,b\} < 2$ and $a \neq b,\;\;a,b \neq 1$ and for $x,y,z \in \{0,\alpha,\beta\}$
let 
$$
f(x,y,z) := p_x H_x \max\{H_y, H_z\} 
+ p_x \min\{H_y, H_z\}^2/4 + (1-p_x) H_y H_z,
$$
\begin{itemize}
\item[{\rm (i)}] $A(\alpha,\beta)$ occurs when 
$f(0,\alpha,\beta)< \min\{ f(\beta, 0,\alpha), f(\alpha,\beta, 0)\}$
\item[{\rm (ii)}] $A(0,\alpha)$ and $A(0,\beta)$ have positive 
probabilities of occurrence, when 
$f(\beta, 0,\alpha)= f(\alpha,\beta, 0) < f(0,\alpha,\beta) $, and
\item[{\rm (iii)}] $A(\alpha,\beta)$, $A(0,\alpha)$ and $A(0,\beta)$ 
all have positive probabilities of occurrence when 
$f(\beta, 0,\alpha)= f(\alpha,\beta, 0) = f(0,\alpha,\beta) $;
\end{itemize}

\item[{\rm (3)}] for  $ 0 < a=b < 1$, and, 
\begin{itemize}
\item[{\rm (i)}] for $p_0 \leq \min\{p_\alpha, p_\beta\}$, $A(\alpha, \beta)$ occurs, 
\item[{\rm (ii)}] for $p_0 > \min\{p_\alpha, p_\beta\}$, 
\\
if 
$a< {\bf l}_1(p_0,p_\alpha, p_\beta)
:= 1-\frac{p_0 -\min\{p_\alpha, p_\beta \}}{4-3p_0 -\min\{p_\alpha,p_\beta\}}$, then $A(\alpha, \beta)$ and fixation occurs,  while,  
\\
if 
$a \ge {\bf l}_1(p_0,p_\alpha, p_\beta)$, $A(0,\alpha)$ occurs for $p_\alpha > 
p_\beta$ and both $A(0,\alpha)$ and $A(0,\beta)$ have positive probability of 
occurrence for $p_\alpha = p_\beta$;
\end{itemize}

\item[{\rm (4)}] for  $ 1< a=b < 2$, and,
\begin{itemize}
\item[{\rm (i)}] for $p_0 < \min\{p_\alpha, p_\beta\}$, 
\\
if $a < {\bf l}_2(p_0,p_\alpha, p_\beta)
:= \frac{2\max\{p_\alpha,p_\beta\} + \sqrt{4\max\{p_\alpha,p_\beta\}^2 + 
4p_\alpha p_\beta + p_0 \min\{p_\alpha,p_\beta\}}}{4\max\{p_\alpha,p_\beta\} + p_0}$,  
then $A(\alpha, \beta)$ and fixation occurs,  while,\\
if $a \ge {\bf l}_2(p_0,p_\alpha, p_\beta)$, $A(0,\alpha)$ occurs for $p_\alpha > 
p_\beta$ and both $A(0,\alpha)$ and $A(0,\beta)$ have positive probability of 
occurrence for $p_\alpha = p_\beta$,
\item[{\rm (ii)}] for $\min\{p_\alpha, p_\beta\} \leq p_0$, 
$A(0,\alpha)$ occurs for $p_\alpha > p_\beta$ and both $A(0,\alpha)$ and 
$A(0,\beta)$ have positive probability of occurrence for $p_\alpha = p_\beta$;
\end{itemize}

\item[{\rm (5)}] for $a = b = 1$, fixation always occurs and
\begin{itemize}
\item[{\rm (i)}] $A(x,y)$ occurs when $p_z < \min\{p_x, p_y\}$,
\item[{\rm (ii)}] with equal probability $A(x,y)$ and $A(x,z)$ occur 
when $p_y = p_z < p_x$, and 
\item[{\rm (iii)}] with equal probability $A(x,y)$, $A(y,z)$ and 
$A(z,x)$ occur when $  p_x = p_y = p_z $;
\end{itemize}

\end{itemize}
\end{theorem}

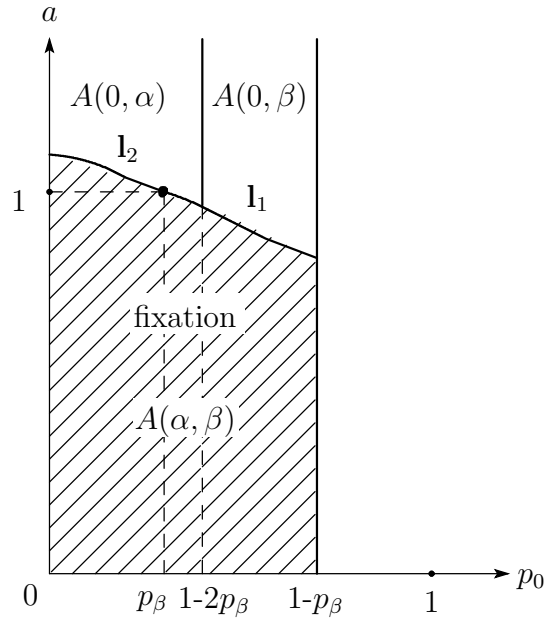
\begin{figure}[h]
\label{Fig2}
\begin{center}
\unitlength 0.1in
\begin{picture}( 26.7500, 30.8000)(  3.7500,-36.6500)
%
\special{pn 8}%
\special{pa 600 3600}%
\special{pa 3000 3600}%
\special{fp}%
\special{sh 1}%
\special{pa 3000 3600}%
\special{pa 2934 3580}%
\special{pa 2948 3600}%
\special{pa 2934 3620}%
\special{pa 3000 3600}%
\special{fp}%
%
\special{pn 8}%
\special{pa 600 3600}%
\special{pa 600 800}%
\special{fp}%
\special{sh 1}%
\special{pa 600 800}%
\special{pa 580 868}%
\special{pa 600 854}%
\special{pa 620 868}%
\special{pa 600 800}%
\special{fp}%
\special{pn 13}%
\special{pa 600 1406}%
\special{pa 606 1406}%
\special{pa 610 1408}%
\special{pa 616 1408}%
\special{pa 620 1410}%
\special{pa 626 1412}%
\special{pa 630 1412}%
\special{pa 636 1414}%
\special{pa 640 1414}%
\special{pa 646 1416}%
\special{pa 650 1418}%
\special{pa 656 1418}%
\special{pa 660 1420}%
\special{pa 666 1420}%
\special{pa 670 1422}%
\special{pa 676 1424}%
\special{pa 680 1424}%
\special{pa 686 1426}%
\special{pa 690 1428}%
\special{pa 696 1428}%
\special{pa 700 1430}%
\special{pa 706 1430}%
\special{pa 710 1432}%
\special{pa 716 1434}%
\special{pa 720 1434}%
\special{pa 726 1436}%
\special{pa 730 1438}%
\special{pa 736 1438}%
\special{pa 740 1440}%
\special{pa 746 1442}%
\special{pa 750 1442}%
\special{pa 756 1444}%
\special{pa 760 1446}%
\special{pa 766 1446}%
\special{pa 770 1448}%
\special{pa 776 1450}%
\special{pa 780 1450}%
\special{pa 786 1452}%
\special{pa 790 1454}%
\special{pa 796 1456}%
\special{pa 800 1456}%
\special{pa 806 1458}%
\special{pa 810 1460}%
\special{pa 816 1460}%
\special{pa 820 1462}%
\special{pa 826 1464}%
\special{pa 830 1466}%
\special{pa 836 1466}%
\special{pa 840 1468}%
\special{pa 846 1470}%
\special{pa 850 1472}%
\special{pa 856 1472}%
\special{pa 860 1474}%
\special{pa 866 1476}%
\special{pa 870 1478}%
\special{pa 876 1478}%
\special{pa 880 1480}%
\special{pa 886 1482}%
\special{pa 890 1484}%
\special{pa 896 1486}%
\special{pa 900 1486}%
\special{pa 906 1488}%
\special{pa 910 1490}%
\special{pa 916 1492}%
\special{pa 920 1494}%
\special{pa 926 1494}%
\special{pa 930 1496}%
\special{pa 936 1498}%
\special{pa 940 1500}%
\special{pa 946 1502}%
\special{pa 950 1502}%
\special{pa 956 1504}%
\special{pa 960 1506}%
\special{pa 966 1508}%
\special{pa 970 1510}%
\special{pa 976 1512}%
\special{pa 980 1514}%
\special{pa 986 1514}%
\special{pa 990 1516}%
\special{pa 996 1518}%
\special{pa 1000 1520}%
\special{pa 1006 1522}%
\special{pa 1010 1524}%
\special{pa 1016 1526}%
\special{pa 1020 1528}%
\special{pa 1026 1530}%
\special{pa 1030 1530}%
\special{pa 1036 1532}%
\special{pa 1040 1534}%
\special{pa 1046 1536}%
\special{pa 1050 1538}%
\special{pa 1056 1540}%
\special{pa 1060 1542}%
\special{pa 1066 1544}%
\special{pa 1070 1546}%
\special{pa 1076 1548}%
\special{pa 1080 1550}%
\special{pa 1086 1552}%
\special{pa 1090 1554}%
\special{pa 1096 1556}%
\special{pa 1100 1558}%
\special{pa 1106 1560}%
\special{pa 1110 1562}%
\special{pa 1116 1564}%
\special{pa 1120 1566}%
\special{pa 1126 1568}%
\special{pa 1130 1570}%
\special{pa 1136 1572}%
\special{pa 1140 1574}%
\special{pa 1146 1576}%
\special{pa 1150 1578}%
\special{pa 1156 1580}%
\special{pa 1160 1582}%
\special{pa 1166 1584}%
\special{pa 1170 1588}%
\special{pa 1176 1590}%
\special{pa 1180 1592}%
\special{pa 1186 1594}%
\special{pa 1190 1596}%
\special{pa 1196 1598}%
\special{pa 1200 1600}%
\special{sp}%
\special{pn 13}%
\special{pa 1200 1600}%
\special{pa 1206 1602}%
\special{pa 1210 1604}%
\special{pa 1216 1606}%
\special{pa 1220 1606}%
\special{pa 1226 1608}%
\special{pa 1230 1610}%
\special{pa 1236 1612}%
\special{pa 1240 1614}%
\special{pa 1246 1616}%
\special{pa 1250 1616}%
\special{pa 1256 1618}%
\special{pa 1260 1620}%
\special{pa 1266 1622}%
\special{pa 1270 1624}%
\special{pa 1276 1626}%
\special{pa 1280 1626}%
\special{pa 1286 1628}%
\special{pa 1290 1630}%
\special{pa 1296 1632}%
\special{pa 1300 1634}%
\special{pa 1306 1636}%
\special{pa 1310 1638}%
\special{pa 1316 1640}%
\special{pa 1320 1640}%
\special{pa 1326 1642}%
\special{pa 1330 1644}%
\special{pa 1336 1646}%
\special{pa 1340 1648}%
\special{pa 1346 1650}%
\special{pa 1350 1652}%
\special{pa 1356 1654}%
\special{pa 1360 1654}%
\special{pa 1366 1656}%
\special{pa 1370 1658}%
\special{pa 1376 1660}%
\special{pa 1380 1662}%
\special{pa 1386 1664}%
\special{pa 1390 1666}%
\special{pa 1396 1668}%
\special{pa 1400 1670}%
\special{pa 1406 1672}%
\special{pa 1410 1674}%
\special{pa 1416 1676}%
\special{pa 1420 1676}%
\special{pa 1426 1678}%
\special{pa 1430 1680}%
\special{pa 1436 1682}%
\special{pa 1440 1684}%
\special{pa 1446 1686}%
\special{pa 1450 1688}%
\special{pa 1456 1690}%
\special{pa 1460 1692}%
\special{pa 1466 1694}%
\special{pa 1470 1696}%
\special{pa 1476 1698}%
\special{pa 1480 1700}%
\special{pa 1486 1702}%
\special{pa 1490 1704}%
\special{pa 1496 1706}%
\special{pa 1500 1708}%
\special{pa 1506 1710}%
\special{pa 1510 1712}%
\special{pa 1516 1714}%
\special{pa 1520 1716}%
\special{pa 1526 1718}%
\special{pa 1530 1720}%
\special{pa 1536 1722}%
\special{pa 1540 1724}%
\special{pa 1546 1726}%
\special{pa 1550 1728}%
\special{pa 1556 1730}%
\special{pa 1560 1732}%
\special{pa 1566 1734}%
\special{pa 1570 1736}%
\special{pa 1576 1738}%
\special{pa 1580 1740}%
\special{pa 1586 1742}%
\special{pa 1590 1744}%
\special{pa 1596 1746}%
\special{pa 1600 1748}%
\special{pa 1606 1750}%
\special{pa 1610 1752}%
\special{pa 1616 1756}%
\special{pa 1620 1758}%
\special{pa 1626 1760}%
\special{pa 1630 1762}%
\special{pa 1636 1764}%
\special{pa 1640 1766}%
\special{pa 1646 1768}%
\special{pa 1650 1770}%
\special{pa 1656 1772}%
\special{pa 1660 1774}%
\special{pa 1666 1776}%
\special{pa 1670 1780}%
\special{pa 1676 1782}%
\special{pa 1680 1784}%
\special{pa 1686 1786}%
\special{pa 1690 1788}%
\special{pa 1696 1790}%
\special{pa 1700 1792}%
\special{pa 1706 1796}%
\special{pa 1710 1798}%
\special{pa 1716 1800}%
\special{pa 1720 1802}%
\special{pa 1726 1804}%
\special{pa 1730 1806}%
\special{pa 1736 1810}%
\special{pa 1740 1812}%
\special{pa 1746 1814}%
\special{pa 1750 1816}%
\special{pa 1756 1818}%
\special{pa 1760 1820}%
\special{pa 1766 1824}%
\special{pa 1770 1826}%
\special{pa 1776 1828}%
\special{pa 1780 1830}%
\special{pa 1786 1834}%
\special{pa 1790 1836}%
\special{pa 1796 1838}%
\special{pa 1800 1840}%
\special{pa 1806 1842}%
\special{pa 1810 1846}%
\special{pa 1816 1848}%
\special{pa 1820 1850}%
\special{pa 1826 1854}%
\special{pa 1830 1856}%
\special{pa 1836 1858}%
\special{pa 1840 1860}%
\special{pa 1846 1864}%
\special{pa 1850 1866}%
\special{pa 1856 1868}%
\special{pa 1860 1870}%
\special{pa 1866 1874}%
\special{pa 1870 1876}%
\special{pa 1876 1878}%
\special{pa 1880 1882}%
\special{pa 1886 1884}%
\special{pa 1890 1886}%
\special{pa 1896 1890}%
\special{pa 1900 1892}%
\special{pa 1906 1894}%
\special{pa 1910 1898}%
\special{pa 1916 1900}%
\special{pa 1920 1904}%
\special{pa 1926 1906}%
\special{pa 1930 1908}%
\special{pa 1936 1912}%
\special{pa 1940 1914}%
\special{pa 1946 1916}%
\special{pa 1950 1920}%
\special{pa 1956 1922}%
\special{pa 1960 1926}%
\special{pa 1966 1928}%
\special{pa 1970 1930}%
\special{pa 1976 1934}%
\special{pa 1980 1936}%
\special{pa 1986 1940}%
\special{pa 1990 1942}%
\special{pa 1996 1946}%
\special{pa 2000 1948}%
\special{sp}%
%
\special{pn 8}%
\special{pa 600 1600}%
\special{pa 1200 1600}%
\special{da 0.070}%
%
\special{pn 8}%
\special{pa 1200 1600}%
\special{pa 1200 3600}%
\special{da 0.070}%
%
\special{pn 13}%
\special{pa 2000 800}%
\special{pa 2000 3600}%
\special{fp}%
%
\special{pn 13}%
\special{pa 1400 800}%
\special{pa 1400 1670}%
\special{fp}%
%
\special{pn 8}%
\special{pa 1400 3600}%
\special{pa 1400 1670}%
\special{da 0.070}%
%
\special{pn 8}%
\special{pa 1990 2210}%
\special{pa 1400 2800}%
\special{fp}%
\special{pa 1990 2330}%
\special{pa 1400 2920}%
\special{fp}%
\special{pa 1990 2450}%
\special{pa 1400 3040}%
\special{fp}%
\special{pa 1990 2570}%
\special{pa 1400 3160}%
\special{fp}%
\special{pa 1990 2690}%
\special{pa 1400 3280}%
\special{fp}%
\special{pa 1990 2810}%
\special{pa 1400 3400}%
\special{fp}%
\special{pa 1990 2930}%
\special{pa 1400 3520}%
\special{fp}%
\special{pa 1990 3050}%
\special{pa 1440 3600}%
\special{fp}%
\special{pa 1990 3170}%
\special{pa 1560 3600}%
\special{fp}%
\special{pa 1990 3290}%
\special{pa 1680 3600}%
\special{fp}%
\special{pa 1990 3410}%
\special{pa 1800 3600}%
\special{fp}%
\special{pa 1990 3530}%
\special{pa 1920 3600}%
\special{fp}%
\special{pa 1990 2090}%
\special{pa 1400 2680}%
\special{fp}%
\special{pa 1990 1970}%
\special{pa 1400 2560}%
\special{fp}%
\special{pa 1920 1920}%
\special{pa 1400 2440}%
\special{fp}%
\special{pa 1840 1880}%
\special{pa 1400 2320}%
\special{fp}%
\special{pa 1760 1840}%
\special{pa 1400 2200}%
\special{fp}%
\special{pa 1680 1800}%
\special{pa 1400 2080}%
\special{fp}%
\special{pa 1600 1760}%
\special{pa 1400 1960}%
\special{fp}%
\special{pa 1520 1720}%
\special{pa 1400 1840}%
\special{fp}%
\special{pa 1430 1690}%
\special{pa 1400 1720}%
\special{fp}%
%
\special{pn 8}%
\special{pa 1200 2640}%
\special{pa 600 3240}%
\special{fp}%
\special{pa 1200 2760}%
\special{pa 600 3360}%
\special{fp}%
\special{pa 1200 2880}%
\special{pa 600 3480}%
\special{fp}%
\special{pa 1200 3000}%
\special{pa 610 3590}%
\special{fp}%
\special{pa 1200 3120}%
\special{pa 720 3600}%
\special{fp}%
\special{pa 1200 3240}%
\special{pa 840 3600}%
\special{fp}%
\special{pa 1200 3360}%
\special{pa 960 3600}%
\special{fp}%
\special{pa 1200 3480}%
\special{pa 1080 3600}%
\special{fp}%
\special{pa 1200 2520}%
\special{pa 600 3120}%
\special{fp}%
\special{pa 1200 2400}%
\special{pa 600 3000}%
\special{fp}%
\special{pa 1200 2280}%
\special{pa 600 2880}%
\special{fp}%
\special{pa 1200 2160}%
\special{pa 600 2760}%
\special{fp}%
\special{pa 1200 2040}%
\special{pa 600 2640}%
\special{fp}%
\special{pa 1200 1920}%
\special{pa 600 2520}%
\special{fp}%
\special{pa 1200 1800}%
\special{pa 600 2400}%
\special{fp}%
\special{pa 1200 1680}%
\special{pa 600 2280}%
\special{fp}%
\special{pa 1160 1600}%
\special{pa 600 2160}%
\special{fp}%
\special{pa 1040 1600}%
\special{pa 600 2040}%
\special{fp}%
\special{pa 920 1600}%
\special{pa 600 1920}%
\special{fp}%
\special{pa 800 1600}%
\special{pa 600 1800}%
\special{fp}%
\special{pa 680 1600}%
\special{pa 600 1680}%
\special{fp}%
%
\special{pn 8}%
\special{pa 1390 2810}%
\special{pa 1200 3000}%
\special{fp}%
\special{pa 1390 2930}%
\special{pa 1200 3120}%
\special{fp}%
\special{pa 1390 3050}%
\special{pa 1200 3240}%
\special{fp}%
\special{pa 1390 3170}%
\special{pa 1200 3360}%
\special{fp}%
\special{pa 1390 3290}%
\special{pa 1200 3480}%
\special{fp}%
\special{pa 1390 3410}%
\special{pa 1210 3590}%
\special{fp}%
\special{pa 1390 3530}%
\special{pa 1320 3600}%
\special{fp}%
\special{pa 1390 2690}%
\special{pa 1200 2880}%
\special{fp}%
\special{pa 1390 2570}%
\special{pa 1200 2760}%
\special{fp}%
\special{pa 1390 2450}%
\special{pa 1200 2640}%
\special{fp}%
\special{pa 1390 2330}%
\special{pa 1200 2520}%
\special{fp}%
\special{pa 1390 2210}%
\special{pa 1200 2400}%
\special{fp}%
\special{pa 1390 2090}%
\special{pa 1200 2280}%
\special{fp}%
\special{pa 1390 1970}%
\special{pa 1200 2160}%
\special{fp}%
\special{pa 1390 1850}%
\special{pa 1200 2040}%
\special{fp}%
\special{pa 1390 1730}%
\special{pa 1200 1920}%
\special{fp}%
\special{pa 1340 1660}%
\special{pa 1200 1800}%
\special{fp}%
\special{pa 1250 1630}%
\special{pa 1200 1680}%
\special{fp}%
%
\special{pn 8}%
\special{pa 810 1470}%
\special{pa 680 1600}%
\special{fp}%
\special{pa 900 1500}%
\special{pa 800 1600}%
\special{fp}%
\special{pa 990 1530}%
\special{pa 920 1600}%
\special{fp}%
\special{pa 1080 1560}%
\special{pa 1040 1600}%
\special{fp}%
\special{pa 720 1440}%
\special{pa 600 1560}%
\special{fp}%
%
\special{pn 13}%
\special{sh 1}%
\special{ar 2600 3600 10 10 0  6.28318530717959E+0000}%
\special{sh 1}%
\special{ar 2600 3600 10 10 0  6.28318530717959E+0000}%
%
\special{pn 13}%
\special{sh 1}%
\special{ar 600 1600 10 10 0  6.28318530717959E+0000}%
\special{sh 1}%
\special{ar 600 1600 10 10 0  6.28318530717959E+0000}%
\put(4.0000,-17.0000){\makebox(0,0)[lb]{1}}%
\put(25.6000,-38.0000){\makebox(0,0)[lb]{1}}%
\put(30.5000,-36.8000){\makebox(0,0)[lb]{$p_{0}$}}%
\put(6.0000,-6.7000){\makebox(0,0){$a$}}%
\put(14.6000,-37.5000){\makebox(0,0){1-2$p_{\beta}$}}%
\put(11.4000,-37.5000){\makebox(0,0){$p_{\beta}$}}%
\put(4.6000,-37.6000){\makebox(0,0)[lb]{0}}%
\put(9.6000,-11.1000){\makebox(0,0){$A(0,\alpha)$}}%
\put(17.0000,-11.1000){\makebox(0,0){$A(0,\beta)$}}%
\put(16.9000,-16.3000){\makebox(0,0){${\bf l}_{1}$}}%
\put(10.1000,-13.7000){\makebox(0,0){${\bf l}_{2}$}}%
%
\special{pn 8}%
\special{sh 0}%
\special{pa 1040 2680}%
\special{pa 1590 2680}%
\special{pa 1590 2880}%
\special{pa 1040 2880}%
\special{pa 1040 2680}%
\special{ip}%
\put(13.1000,-28.0000){\makebox(0,0){$A(\alpha,\beta)$}}%
%
\special{pn 8}%
\special{sh 0}%
\special{pa 1020 2130}%
\special{pa 1620 2130}%
\special{pa 1620 2340}%
\special{pa 1020 2340}%
\special{pa 1020 2130}%
\special{ip}%
\put(13.1000,-22.5000){\makebox(0,0){fixation}}%
\put(19.9000,-37.5000){\makebox(0,0){1-$p_{\beta}$}}%
%
\special{pn 20}%
\special{sh 1}%
\special{ar 1200 1600 10 10 0  6.28318530717959E+0000}%
\special{sh 1}%
\special{ar 1200 1600 10 10 0  6.28318530717959E+0000}%
%
\special{pn 20}%
\special{sh 1}%
\special{ar 1190 1590 10 10 0  6.28318530717959E+0000}%
\special{sh 1}%
\special{ar 1200 1590 10 10 0  6.28318530717959E+0000}%
\special{sh 1}%
\special{ar 1190 1600 10 10 0  6.28318530717959E+0000}%
\special{sh 1}%
\special{ar 1190 1610 10 10 0  6.28318530717959E+0000}%
\end{picture}%
\end{center}
\caption{ \it{The diagram in the case that $a=b$ and $p_\beta \in (0,1/3)$.
The curved line is the line 
${\bf l}_1 1_{\{0 \leq {\bf l}_1 \leq 1\}} + {\bf l}_2 1_{\{1 \leq {\bf l}_1 \leq 2\}}$. 
For $p_0 > 0$ and $a$ below this line $A(\alpha, \beta) \ occurs$, 
while for $a$ above the line $A(0, \beta)$ occurs 
when $p_\alpha < p_\beta$. 
At $p_0=0$, only $A(\alpha, \beta)$ occurs.}}
\end{figure}

Observe that for $\min{a,b} \leq 1/2$:

\noindent (A) \ If $b, 1 \ge 2a$, then by the scaling which transforms $a$ to $1$, $b$ to $b/a$ and $1$ to $1/a$, the resulting asymptotic cluster may be read from (1) of Theorem \ref{pshape}. Similarly if $a,1 \ge 2b$, we may scale suitably to obtain a situation as in (1) of Theorem \ref{pshape}.

\noindent (B) \ If either $ a/2 < \min \{ 1,b \} <2a, a\not= b, \ a,b\not=1$, or $b/2< \min \{1,a\} <2b, a\not= b, \ a,b\not=1$, then scaling shows that (2) of Theorem \ref{pshape} may be used to yield the asymptotic shape.

\noindent (C) \ If either $0 < b=1<a $ or $0< a=1 <b$, then scaling shows that (3) of Theorem \ref{pshape} may be used to yield the asymptotic shape.

\noindent (D) \ If either $a < b=1< 2a $ or $b < a=1 < 2b$, then scaling shows that (4) of Theorem \ref{pshape} may be used to yield the asymptotic shape.

\noindent Thus the above four observations demonstrate that Theorem \ref{pshape} yields the asymptotic shapes for all possible values of $a$ and $b$.


\begin{figure}[h]
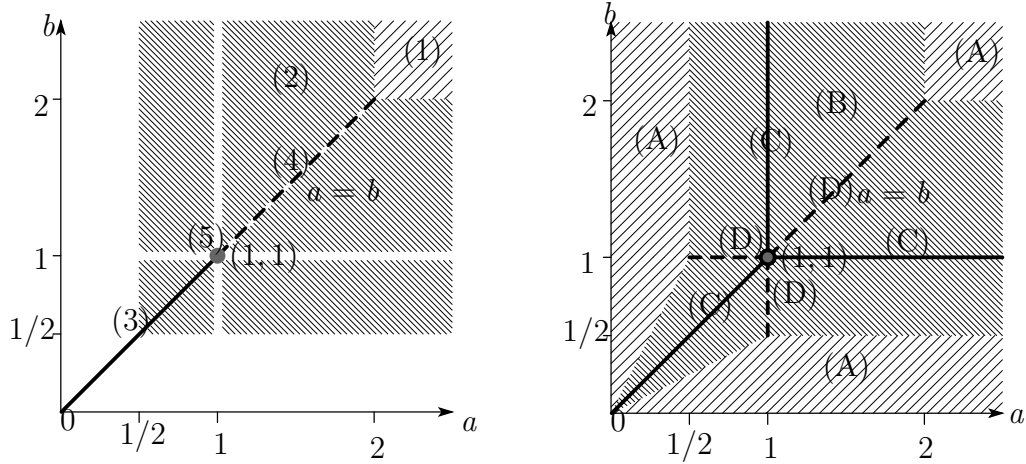

\label{Fig3}
\begin{center}
\unitlength 0.1in
%
\end{center}
\caption{\it{The various regions where Theorem the various parts of Theorem 2.2 hold.}}
\end{figure}

To prove  the above theorem we need to know the conditional probability of 
the composition of a cluster given that it is finite. 
This is obtained in the next two sections.

\newsection{Proof of Theorem 2.1}
\subsection{General set-up}
For ${\bf k} \in (\mathbb N \cup 0)^d$, $d \geq 2$, with $|{\bf k}| = m$, let $\Lambda({\bf k})$ and $
\Gamma_0$ be as in Section 2.1.  First we calculate $\mu_{\lambda\rho}(C_{\bf 0}\in\Lambda({\bf k})|
\Gamma_0)$.
Suppose that $w_m = ({\bf 0}, \alpha_{j_0}, R_{j_0})$ 
for some $j_0 \in \{ 1,2,\dots,d \}$. We have
\begin{eqnarray*}
&&\mu_{\lambda \rho} ( C_{\bf 0} \in \Lambda({\bf k}) \mid ~ w_m \in \xi )
\\
&&=\int \limits _{{\cal M}} 
\mu _{\lambda \rho} (d\xi)
\displaystyle{\sum _{ \{ {\bf w}_{m-1} \} \subset \xi} }
1_{\Lambda({\bf k})} (C_{\bf 0} ( {\bf w}_m ))
1_{\{S(\xi \setminus \{{\bf w}_m \})\cap S(\{ {\bf w}_m \}) =\emptyset \}},
\end{eqnarray*}
where ${\bf w_m}$, $\{{\bf w_m}\}$ and $C_{\bf 0} ( {\bf w}_m )$ are as defined in (\ref{p-def}). Thus,
\begin{eqnarray*}
&&\mu_{\lambda \rho} ( C_{\bf 0} \in \Lambda({\bf k}) \mid ~ w_m \in \xi )
\\
&&=\frac{\lambda^{m-1}}{(m-1)!}
\int \limits _{{\cal M}} \mu _{\lambda \rho} (d\eta)
\int \limits _{{\cal R}^{m-1} } \rho^{\otimes (m-1)}(d {\bf w}_{m-1})
1_{\Lambda({\bf k})} (C_{\bf 0} ( {\bf w}_m ))
1_{\{ S(\eta) \cap S(\{ {\bf w}_m \}) = \emptyset \}}
\\
&&=\frac{\lambda^{m-1}}{(m-1)!}
\int \limits _{{\cal R}^{m-1} } 
\rho^{\otimes (m-1)}(d {\bf w}_{m-1})
1_{\Lambda({\bf k})} (C_{\bf 0} ( {\bf w}_m ))
e^{-\lambda \rho ( w : S(w) \cap S(\{ {\bf w}_m \}) \not= \emptyset)}.
\end{eqnarray*}
Note that  $S(x,\theta,r)\cap S(\{ {\bf w}_m \})\not=\emptyset$
if and only if
$x \in \cup_{i=1}^{m}B^{\theta_i, \theta}_{r_i,r}(x_i)$ where
$w_i = (x_i, \theta_i, r_i)$, $i=1,2,\dots, m$.
Hence, 
$$
\rho ( w : S(w) \cap S(\{ {\bf w}_m \}) \not= \emptyset)
= \sum_{j=1}^d p_j | \bigcup_{i=1}^m
B^{\theta_i, \alpha_j}_{r_i, R_j} (x_i)|,
$$
and so
\begin{equation*}
\begin{split}
\mu_{\lambda \rho}( C_{\bf 0} \in \Lambda ({\bf k}) \mid ~ w_m \in \xi )
&=\frac{\lambda^{m-1}}{(m-1)!}
\int \limits _{{\cal R}^{m-1} } 
\rho^{\otimes (m-1)}(d {\bf w}_{m-1})
1_{\Lambda({\bf k})} (C_{\bf 0} ( {\bf w}_m ))
\\
&\quad \times
\exp \left[ -\lambda 
\sum_{j=1}^d p_j | \bigcup_{i=1}^m 
B^{\theta_i, \alpha_j}_{r_i, R_j} (x_i)|
\right].
\end{split}
\end{equation*}
Let
\begin{equation*}
\begin{split}
F_\lambda^{\alpha_{j_0}} ({\bf k})
&=\int \limits _{(\Bbb R ^2)^{k_1} } d{\bf x}_{1,k_1}
\int \limits _{(\Bbb R ^2)^{k_2} } d{\bf x}_{2,k_2}
\cdots
\int \limits _{(\Bbb R ^2)^{k_{j_0}-1} } d{\bf x}_{j_0,k_{j_0}-1}
\cdots
\int \limits _{(\Bbb R ^2)^{k_d} } d{\bf x}_{d,k_d}
\\
&\qquad \times 1_{\Lambda({\bf k})} 
(C_{\bf 0} ( {\bf x} ))
\exp \left[ -\lambda 
\sum_{j=1}^{d} p_j 
| \bigcup_{i=1, k_i\not= 0}^{d} 
B^{\alpha_i, \alpha_j} _{R_i,R_j} ({\bf x}_{i,k_i})|
\right],
\end{split}
\end{equation*}
where
$C_{\bf 0}({\bf x})
= C_{\bf 0}({\bf x}_{1,k_1}, {\bf x}_{2,k_2},\dots,{\bf x}_{d,k_d})
=C_{\bf 0}(\bigcup_{j=1}^{d}\{ (x_{j,i},\alpha_j,R_j): i=1,\dots,k_j \})$.
>From the translation invariance of Lebesgue measure
it is obvious that 
if $k_j, k_{j'}\ge 1$, then 
$ F_\lambda^{\alpha_{j}} ({\bf k})=F_\lambda^{\alpha_{j'}} ({\bf k})$.
Thus  writing $ F_\lambda ({\bf k})$ for $F_\lambda^{\alpha_{j}} ({\bf k})$,
since 
$\mu_{\lambda \rho}((0,\alpha_j,R_j) \in \xi \mid ~ \Gamma_{{\bf 0}})
= p_j$, we have
\begin{equation}
\label{lsticks}
\mu_{\lambda \rho}(C_{\bf 0}\in\Lambda({\bf k})\mid ~ \Gamma_{{\bf 0}})
= \frac{\lambda^{m-1}}{(m-1)!}\prod_{j=1}^d \frac{m!}{k_j!}p_j^{k_j}
F_\lambda ({\bf k})
=\lambda^{|{\bf k}|-1}|{\bf k}|
\prod_{j=1}^d \frac{p_j^{k_j}}{k_j!}
F_\lambda ({\bf k}).
\end{equation}

\subsection{Proof of Theorem \ref{th2sticks}}
To prove Theorem \ref{th2sticks}, observe first that in the case when we have sticks with only two orientations, the Radon measure $\rho $ is given by
\begin{equation}
\label{rho2}
\rho (dx ~d {\theta } ~dr ) =
dx \{ p\delta _0 (d {\theta }) \delta _{R_{0}} (dr )
+ q\delta _\alpha ( d {\theta }) \delta _{R_{\alpha }} (dr ) \}.
\end{equation}
>From (\ref{lsticks}) we have
\begin{eqnarray*}
\label{2stick-1}
\mu_{\lambda \rho}( C_{\bf 0}\in\Lambda(k,\ell)\mid ~ \Gamma_{\bf 0})
&=& \lambda^{k+\ell -1}(k+\ell)
\frac{p^k q^{\ell}}{k! \ell !}F^0_\lambda((k,\ell))\nonumber
\\
&=&\lambda^{k+\ell -1}(k+\ell)
\frac{p^k q^{\ell}}{k! \ell !}
e^{-\lambda |B_{R_0,R_{\alpha}}^{0,\alpha}|} f_\lambda(k,\ell),
\end{eqnarray*}
where
$$
f_\lambda (k,\ell) := \int \limits_{(\Bbb R^{2})^{k-1}} d{{\bf x}} _{k-1}
\int \limits_{(\Bbb R^{2})^{l}} d{{\bf y}} _\ell ~
1_{\Lambda (k,\ell)}(C_{\bf 0}({\bf x}_k,{\bf y}_\ell))
\chi^{0,\alpha}_{p\lambda} ({\bf y}_{\ell}) 
\chi^{0,\alpha}_{q\lambda} ({\bf x}_{k}),
$$
\begin{eqnarray}
\chi^{\theta_1,\theta_2}_c ({\bf x}) 
= \exp \left[ -c
\{|B^{\theta_1,\theta_2}_{R_{\theta_1},R_{\theta_2}}({\bf x})|
-|B^{\theta_1,\theta_2}_{R_{\theta_1},R_{\theta_2}}|\} \right]
\end{eqnarray}
(note here that $x_k = {\bf 0}$).
Now consider the event $A({\bf x}_k, {\bf y}_\ell, k,\ell):= \{C_0$ contains exactly $m$ sticks $({\bf 0}, 0, 1/2), 
(x_1,0,1/2), \ldots, (x_{k-1},0,1/2), (y_1,\frac{\pi}{2},1/2), \ldots, (y_\ell,\frac{\pi}{2},1/2)\}$.
By the affine invariance of the Lebesgue measure
\begin{eqnarray}
\label{fund2}
f_\lambda (k,\ell) &=& |B_{R_0,R_{\alpha}}^{0,\alpha}|^{m-1}
\int \limits _{(\Bbb R^{2})^{k-1}} d{{\bf x}} _{k-1}
\int \limits _{(\Bbb R^{2})^{\ell}} d{{\bf y}} _\ell ~
1_{A({\bf x}_k, {\bf y}_\ell, k,\ell)}
\nonumber
\\
&&\times \exp [ -\lambda p |B_{R_0,R_{\alpha}}^{0,\alpha}|
\{|B_{\frac{1}{2}} ({\bf y}_{\ell}) | - |B_{\frac{1}{2}}| \} ]
\nonumber
\\
&&\times\exp [ -\lambda q |B_{R_0,R_{\alpha}}^{0,\alpha}|
\{|B_{\frac{1}{2}}({\bf x}_k) |- |B_{\frac{1}{2}}| \}],
\end{eqnarray}
where $B_R = [-R, R]^2$, $B_R(x)= B_R + x$
and $B_R ({\bf x}_k) = \cup_{i=1}^k B_R (x_i)$.

For the proof of Theorem \ref{th2sticks} we will 
obtain lower and upper bounds of
$f_\lambda(k,l)$ which we later show to agree 
as $\lambda \rightarrow \infty$.
To this end we need the following lemma whose proof is given in the appendix. 
For each $x \in {\Bbb R}^2$ we take
$x^\alpha, x^\beta \in {\Bbb R}$
such that 
$x = x^\alpha e_\alpha  + x^\beta e_\beta$.
Note that $(x^\alpha, x^\beta)$ is just the representation
of $x \in \Bbb R^2$ in the base given by the axes parallel to the orientation of the sticks.
Let $h_\alpha (x) = \frac{x^\alpha}{\sin \beta}$,
$h_\beta (x)= \frac{x^\beta}{\sin \alpha}$ and 
$$
h_\theta ({\bf x}_k) 
= (h_\theta (x_1), h_\theta (x_2),\dots,h_\theta (x_k)),
\quad {\bf x}_k= (x_1,x_2,\dots,x_k)\in (\Bbb R^2)^k.
$$
We put 
$$
M({\bf u_k})=\max_{1\le i,j \le k}|u_i-u_j|, 
\quad {\bf u}_k= (u_1,u_2,\dots,u_k)\in (\Bbb R)^k.
$$
and $C_{\alpha,\beta}=\sin \alpha \sin \beta \sin (\alpha-\beta)$.

 \begin{lemma}
\label{lemma3.1}
Let ${{\bf x}}_k = (x_1, x_2, \cdots , x_k) \in (\Bbb R^2)^k$ 
with $x_k = {\bf 0}$.  Then
\begin{eqnarray}
\label{lem33_i}
|B^{\alpha,\beta}_{R_\alpha, R_\beta} ({\bf x}_k)
\backslash B^{\alpha,\beta}_{R_\alpha, R_\beta}|
&\le& 
2C_{\alpha,\beta}
\{  H_\alpha M (h_\beta ({\bf x}_k))+ H_\beta M(h_\alpha({\bf x}_k))\}
\nonumber
\\
&+& C_{\alpha,\beta}
M(h_\beta({\bf x}_k))M(h_\alpha({\bf x}_k)),
\end{eqnarray}
and, if $B^{\alpha,\beta}_{R_\alpha, R_\beta}({{\bf x}}_k)$ 
is connected, then we have 
\begin{eqnarray}
\label{lem33_iia}
|B^{\alpha,\beta}_{R_\alpha, R_\beta} ({{\bf x}}_k)
\backslash B^{\alpha,\beta}_{R_\alpha, R_\beta}|
&\ge&
C_{\alpha,\beta}
\{H_\alpha M (h_\beta ({\bf x}_k))+ H_\beta M(h_\alpha({\bf x}_k))\},
\end{eqnarray}
\begin{eqnarray}
\label{lem33_ii}
|B^{\alpha,\beta}_{R_\alpha, R_\beta} ({{\bf x}}_k)
\backslash B^{\alpha,\beta}_{R_\alpha, R_\beta}|
&\ge&
2C_{\alpha,\beta}
\{H_\alpha M (h_\beta ({\bf x}_k))+ H_\beta M(h_\alpha({\bf x}_k))\}\nonumber
\\
&-& C_{\alpha,\beta}
M(h_\beta({\bf x}_k))M(h_\alpha({\bf x}_k)).
\end{eqnarray}
\end{lemma}

Now we evaluate the bounds of $f_\lambda(k,\ell)$.\\
\noindent {\sc Lower bound} :
By (\ref{lem33_i}) of Lemma \ref{lemma3.1}, taking $x_k = {\bf 0}$  we have
\begin{eqnarray}
\label{1.5}
f_\lambda (k,\ell) & \ge & |B_{R_0,R_{\alpha}}^{0,\alpha}|^{m-1}
\int \limits _{(\Bbb R^{2})^{k-1}} d{{\bf x}} _{k-1} 
\int \limits _{(\Bbb R^{2})^{\ell}} d{{\bf y}}_\ell
\;\; 1_{A({\bf x}_k, {\bf y}_\ell, k,\ell)}
\nonumber 
\\
& & \times \exp [
-\lambda q |B_{R_0,R_{\alpha}}^{0,\alpha}|
 (M({\bf x}^1_k)+ M({\bf x}^2_k))]
\nonumber 
\\
& & \times \exp [
-\lambda p |B_{R_0,R_{\alpha}}^{0,\alpha}|
(M({\bf y}^1_\ell)+ M({\bf y}^2_\ell))]
\nonumber 
\\
& & \times \exp [-\lambda |B_{R_0,R_{\alpha}}^{0,\alpha}| \{ q M({\bf x}^1_k)M({\bf x}^2_k) + p M({\bf y}^1_\ell)M({\bf y}^2_\ell) \}].
\end{eqnarray}
Let $L(\lambda)$ be such that, as $\lambda \rightarrow \infty$,
$\lambda L(\lambda) \rightarrow \infty$ and
$\lambda (L(\lambda))^2 \rightarrow 0$.
If $\{x_i\}_{i=1}^{k-1}\subset B_{L(\lambda)}$
and $\{y_i\}_{i=1}^{\ell-1}\subset B_{L(\lambda)}(y_\ell)$,
then, for $x_k = {\bf 0}$,
$y_\ell \in B_{R-{L(\lambda)}}$
and for $\lambda$ sufficiently large, we have
$A({\bf x}_k, {\bf y}_\ell, k,\ell)$ occurs, and so
the expression on the right of the inequality (\ref{1.5}) is bounded from below by
\begin{eqnarray}
\label{star1}
|B_{R_0,R_{\alpha}}^{0,\alpha}|^{m-1}
\lefteqn{\int \limits _{(B_{L(\lambda)})^{k-1}} d{{\bf x}} _{k-1} 
\int \limits _{B_{1/2-L(\lambda)}} dy_\ell
\int \limits _{(B_{L(\lambda)}(y_{\ell}))^{\ell-1}} d{{\bf y}} _{\ell-1}} 
\nonumber\\
& & \quad \times \exp [
-\lambda q |B_{R_0,R_{\alpha}}^{0,\alpha}|
 (M({\bf x}^1_k)+ M({\bf x}^2_k))]
\nonumber\\
& & \quad \times \exp [
-\lambda p |B_{R_0,R_{\alpha}}^{0,\alpha}|
(M({\bf y}^1_k)+ M({\bf y}^2_k))]
\nonumber\\
& & \quad\times \exp [-\lambda |B_{R_0,R_{\alpha}}^{0,\alpha}|
\{ q M({\bf x}^1_k)M({\bf x}^2_k) + p M({\bf y}^1_k)M({\bf y}^2_k) \}]
\nonumber\\
& \ge & |B_{R_0,R_{\alpha}}^{0,\alpha}|^{m-1}
e^{-4(p+q)(L(\lambda))^2 } | B_{1/2-L(\lambda)}|
\int \limits_{(B_{L(\lambda)})^{k-1}} d{{\bf x}} _{k-1}
\int \limits _{(B_{L(\lambda)})^{\ell-1}}d{{\bf y}} _{\ell-1}
\nonumber\\
& & \quad \times \exp [
-\lambda q |B_{R_0,R_{\alpha}}^{0,\alpha}|
 (M({\bf x}^1_k)+ M({\bf x}^2_k))]
\nonumber\\
& & \quad \times \exp [
-\lambda p |B_{R_0,R_{\alpha}}^{0,\alpha}|
(M({\bf y}^1_k)+ M({\bf y}^2_k))]
\nonumber\\
& =& e^{-4\lambda (L(\lambda))^2} | B_{1/2-L(\lambda)}|
(q\lambda)^{-2(k-1)} (p\lambda)^{-2(\ell-1)}
|B_{R_0,R_{\alpha}}^{0,\alpha}|^{-(m-3)}
\nonumber\\
& & \quad \times 
\int \limits _{
(B_{q \lambda_{\alpha} L(\lambda)})^{k-1}} 
d{{\bf u}} _{k-1}
\exp [ - M({\bf u}_k^1) -M({\bf u}_k^2)]
\nonumber\\
& & \quad \times 
\int \limits _{(B_{p\lambda_{\alpha} L(\lambda)})^{\ell-1}} 
d{{\bf v}} _{\ell-1}
\exp [- M({\bf v}_k^1) -M({\bf v}_k^2)]
\end{eqnarray}
where ${{\bf u}}_k = (u_1, \ldots, u_k)$ and
${{\bf v}}_\ell = (v_1, \ldots, v_\ell)$ with $v_\ell = u_k = {\bf 0}$,
and $\lambda_{\alpha} 
=|B_{R_0,R_{\alpha}}^{0,\alpha}|\lambda$.
Then we have
\begin{eqnarray}
\label{1.6}
f_\lambda (k,\ell) 
& \ge & e^{-4\lambda (L(\lambda))^2} |B_{R-L(\lambda)}| \lambda ^{-2(m-2)}
|B_{R_0,R_{\alpha}}^{0,\alpha}|^{-(m-3)} q^{-2(k-1)} p^{-2(\ell-1)} 
\nonumber\\
& \times & \left [ 
\int \limits ^{q\lambda_{\alpha} L(\lambda)}_{-q\lambda_{\alpha} L(\lambda)} 
da_1\cdots
\int \limits ^{q\lambda_{\alpha} L(\lambda)}_{-q\lambda_{\alpha} L(\lambda)}
da_{k-1} 
~\exp \{ - \mathop{\max}_{1\le i, j\le k} |a_i-a_j|\} \right ]^2 \nonumber 
\\
& \times & \left [ 
\int \limits ^{p\lambda_{\alpha} L(\lambda)}_{-p\lambda_{\alpha} L(\lambda)}
db_1\cdots
\int \limits ^{p\lambda_{\alpha} L(\lambda)}_{-p\lambda_{\alpha} L(\lambda)}
db_{\ell-1} ~
\exp \{ - \mathop{\max}_{1\le i, j\le
\ell} |b_i-b_j|\} \right ]^2.
\end{eqnarray}

Since $e^{-4\lambda (L(\lambda))^2} = 1 - O(\lambda (L(\lambda))^2)$
as $\lambda \rightarrow 0$,
by (\ref{1.6}) and the above lemma we obtain that, as $\lambda \rightarrow 0$,
\begin{equation}
\label{1.7}
f_\lambda (k,\ell) \ge \left [\left ( \frac {1}{\lambda }\right )^{2(m-2)}
\left( \frac{1}{|B^{0,\alpha}_{R_0,R_{\alpha}}|} \right )^{m-3} 
q^{-2(k-1)}p^{-2(\ell-1)} (k!)^2 (\ell!)^2  \right ]
(1-O(\lambda (L(\lambda))^2)).
\end{equation}

Now we will obtain the upper bound of $f_\lambda(k,\ell)$.

\noindent {\sc Upper bound:}
For $L(\lambda)$ as earlier, consider the event\\
$E:=\{x_1, \ldots, x_{k-1} \in B_{L(\lambda)}, y_1, \ldots, y_{\ell-1} \in B_{L(\lambda)}(y_\ell)\}$.\\
If $x_k = {\bf 0}$, for $E \cap A({\bf x}_k, {\bf y}_\ell, k,\ell)$ to occur, we must have $y_\ell \in B_{(1/2) + L(\lambda)}$. Thus from (\ref{fund2}) we have
\begin{eqnarray}
\label{fund3}
f_\lambda (k,\ell) &\leq& |B_{R_0,R_{\alpha}}^{0,\alpha}|^{m-1}
\int \limits _{(\Bbb R^{2})^{k-1}} d{{\bf x}} _{k-1}
\int \limits_{\Bbb R^2} d y_\ell
\int \limits _{(\Bbb R^{2})^{\ell}} d{{\bf y}} _{\ell-1} ~ \nonumber\\
&\times&(1_{E \cap \{y_\ell \in B_{(1/2) + L(\lambda)}\}} + 1_{E^c} 1_{A({\bf x}_k, {\bf y}_\ell, k,\ell)})
\nonumber
\\
&& \quad \times \exp [ -\lambda p |B_{R_0,R_{\alpha}}^{0,\alpha}|
\{|B_{\frac{1}{2}} ({\bf y}_{\ell}) | - |B_{\frac{1}{2}}| \} ]
\nonumber
\\
&&\quad \times\exp [ -\lambda q |B_{R_0,R_{\alpha}}^{0,\alpha}|
\{|B_{\frac{1}{2}}({\bf x}_k) |- |B_{\frac{1}{2}}| \}].
\end{eqnarray}
On opening the parenthesis $(1_{E \cap \{y_\ell \in B_{(1/2) + L(\lambda)}\}} + 1_{E^c} 1_{A({\bf x}_k, {\bf y}_\ell, k,\ell)})$ in the expression on the right of the inequality (\ref{fund3}) above the term involving 
$1_{E \cap \{y_\ell \in B_{(1/2) + L(\lambda)}\}}$, 
for large $\lambda$, may be bounded from above by 
\begin{eqnarray}
\label{star2}
\lefteqn{e^{4\lambda (L(\lambda))^2} | B_{1/2+L(\lambda)}|
(q\lambda)^{-2(k-1)} (p\lambda)^{-2(\ell-1)}
|B_{R_0,R_{\alpha}}^{0,\alpha}|^{-(m-3)}}
\nonumber\\
& &  \times 
\int \limits _{
(B_{q \lambda_{\alpha} L(\lambda)})^{k-1}} 
d{{\bf u}} _{k-1}
\exp [ - M({\bf u}_k^1) -M({\bf u}_k^2)]
\nonumber\\
& &  \times 
\int \limits _{(B_{p\lambda_{\alpha} L(\lambda)})^{\ell-1}} 
d{{\bf v}} _{\ell-1}
\exp [- M({\bf v}_k^1) -M({\bf v}_k^2)].
\end{eqnarray}
(Here we have used the inequality (\ref{lem33_ii}) of Lemma \ref{lemma3.1} and calculations similar to those leading to (\ref{star1}).) 

Using the inequality (\ref{lem33_iia}) of Lemma \ref{lemma3.1} we bound the expression involving $1_{E^c} 1_{A({\bf x}_k, {\bf y}_\ell, k,\ell)}$ in the right of the inequality (\ref{fund3}) by  $|B_{R_0,R_{\alpha}}^{0,\alpha}|^{m-1}
\{ I_1 + I_2 \}$,
where
\begin{eqnarray*}
I_1 & = & \int\limits_{(\Bbb R^{2})^{k-1} \backslash (B_{L(\lambda)})^{k-1}}
d{{\bf x}} _{k-1}
\int\limits_{B_{m}} dy_\ell \int\limits_{(\Bbb R^{2})^{\ell-1}} d{{\bf y}} _{\ell-1} 
\\
& & \times \exp \{ - (q/2)\lambda (M({\bf x}^1_k) +M({\bf x}^2_k))\} 
\exp \{ - (p/2)\lambda (M({\bf y}^1_\ell) +M({\bf y}^2_\ell))\} 
\end{eqnarray*}
and
\begin{eqnarray*}
I_2 & = & \int\limits_{(\Bbb R^2)^{k-1}} d {{\bf x}} _{k-1}
\int\limits_{B_{m}} dy_\ell \int\limits_{(\Bbb R^{2})^{\ell-1} 
\backslash (B_{L(\lambda)})^{\ell-1}}d{{\bf y}} _{\ell-1} 
\\
& & \times \exp \{ - (q/2)\lambda (M({\bf x}^1_k) +M({\bf x}^2_k))\} 
\exp \{ - (p/2)\lambda (M({\bf y}^1_\ell) +M({\bf y}^2_\ell))\}.
\end{eqnarray*}

Let $a_k = 0$.  Then, it is easy to see that
$$
\int \limits _{\Bbb R^{k-1}} da_1 \cdots da_{k-1}
\exp\{ -\mathop{\max}_{1\le i,j\le k} |a_i-a_j|\} = k!.
$$
Using this equation and calculations as in 
(\ref{1.6}) and (\ref{1.7}), for $\lambda \to \infty$, the expression in (\ref{star2}) may be bounded above by
$$
\left [\left ( \frac {1}{\lambda }\right )^{2(m-2)}
\left( \frac{1}{|B^{0,\alpha}_{R_0,R_{\alpha}}|} \right )^{m-3} 
q^{-2(k-1)}p^{-2(\ell-1)} (k!)^2 (\ell!)^2  \right ]
(1+O(\lambda (L(\lambda))^2)).
$$
Thus to show that, asymptotically in $\lambda$ the lower bound (\ref{1.7}) of
$f(k,\ell)$ agrees with its upper bound it suffices to show that 
\begin{equation}
\label{1.10}
I_1 + I_2 = O (\lambda ^{-2m-3}) ~\mbox {as}~ \lambda \rightarrow \infty .
\end{equation}

To estimate the integrals $I_1$ and $I_2$, we use the symmetry of the
integrand in $I_1$ to obtain
\begin{eqnarray*}
I_1 & \le & 4 (k-1) \int \limits _{(\Bbb R^{2})^{k-2}}  d{{\bf x}} _{k-2}
\int \limits _{\Bbb R} dx^1_{k-1}
\int \limits ^\infty _{L(\lambda)} dx^2_{k-1} ~|B_{m}|
\int \limits _{(\Bbb R^{2})^{\ell-1}}
d{{\bf y}} _{\ell-1} 
\\
& & \times \exp \{ - (q/2)\lambda (M({\bf x}^1_k) +M({\bf x}^2_k))\} 
\exp \{ - (p/2)\lambda (M({\bf y}^1_\ell) +M({\bf y}^2_\ell))\} 
\\
& = & 4(k-1) ~|B_{m}| ~\left(\frac{q\lambda}{2}\right)^{-2(k-1)} \left(\frac{p\lambda}{2}\right)^{-2(\ell-1)} k! (\ell!)^2 
\\
& & \times \int \limits _{\Bbb R^{k-2}} da_1  \cdots da_{k-2}
\int \limits ^{\infty }_{q\lambda L(\lambda)} da_{k-1} 
\exp \{ - \mathop{\max}_{1\le i,j\le k}
|a_i-a_j|\}.
\end{eqnarray*}
Since $a_k = 0$, we have  the inequality 
$\mathop{\max}_{1\le i,j\le k} |a_i-a_j| \ge 
\frac {1}{2}
\mathop{\max}_{{1\le i,j\le k}\atop {i,j\neq k-1}} |a_i - a_j| 
+\frac {1}{2}|a_{k-1}|$,
which we use to obtain
\begin{eqnarray*}
\lefteqn {\int \limits _{\Bbb R^{k-2}} da_1  \cdots da_{k-2}
\int \limits ^\infty _{q\lambda L(\lambda)} d a_{k-1}
\exp \{ -\mathop{\max}_{1\le i,j\le k} |a_i-a_j|\} }
\\
& \le & 2^{k-1} \int \limits _{\Bbb R^{k-1}} da_1 da_2 \cdots da_{k-2}
\exp \{ -\mathop{\max}_{1\le i,j\le k} |a_i - a_j|\}
\int \limits ^\infty _{\frac{1}{2}q \lambda L(\lambda)} da_{k-1} e^{-a_{k-1}} 
\\
& = & 2^{k-1} (k-1)! e^{- \frac{1}{2}q \lambda L(\lambda)}.
\end{eqnarray*}
Hence
\begin{eqnarray*}
I_1 & \le & 2^{k+1} ~|B_{m}|~ \lambda^{-2(m-2)} \left(\frac{p}{2}\right)^{-2(\ell-1)} \left(\frac{q}{2}\right)^{-2(k-1)}
(k!)^2(\ell!)^2 e^{- \frac{1}{2}q\lambda L(\lambda)} 
\\
& = & o (e^{-\frac{1}{2}q \lambda L(\lambda)}) 
\qquad \mbox{ as } \lambda \rightarrow \infty.
\end{eqnarray*}
Similarly we obtain
$$
I_2 = o (e^{-\frac{1}{2}p \lambda L(\lambda)}) 
\qquad \mbox{ as } \lambda \rightarrow \infty. 
$$
Now fix $0< \delta < 1/2$
and take $L(\lambda) = \lambda ^{-1+ (\delta / 2)}$.
The bounds obtained above for $I_1$ and $I_2$
show that (\ref{1.10}) holds.

This proves Theorem \ref{th2sticks}(i).
The second part of Theorem \ref{th2sticks} is derived easily from the first part.

\newsection{Proof of Theorem \ref{pshape}}

We now prove Theorem \ref{pshape}. Towards this end we need some estimates on the areas of the unions of various parallelograms. These are presented in the next subsection. The proof of these results are given in the appendix.

\subsection{Area estimates}

Throughout this section we assume $0 < \alpha < \beta < \pi$.
\vskip 1em


\begin{lemma}
\label{lemma4.1}
{\rm (i)}~ If $H_{\alpha}, H_{\beta} >  2H_0$, then
\begin{equation*}
|B_{R_0, R_{\alpha}}^{0,\alpha} \cup B_{R_0, R_{\beta}}^{0,\beta}|
=4 C_{\alpha,\beta}H_0 ( H_\alpha + H_\beta - H_0 ).
\end{equation*}
\\
{\rm (ii)}~ If $\min\{ H_{\alpha}, H_\beta\} \le 2H_0$,
then
\begin{equation*}
|B_{R_0, R_{\alpha}}^{0,\alpha} \cup B_{R_0, R_{\beta}}^{0,\beta}|
= C_{\alpha,\beta}\{ 4H_0 \max \{ H_\alpha, H_\beta \}
+  \min \{ H_\alpha^2, H_\beta^2 \}  \}.
\end{equation*}
\end{lemma}



Next we will estimate

\begin{equation}
\triangle (x)
=
\frac{1}{C_{\alpha,\beta}}
\{
|B_{R_0, R_{\alpha}}^{0,\alpha} \cup B_{R_0, R_{\beta}}^{0,\beta}(x)|
-|B_{R_0, R_{\alpha}}^{0,\alpha} \cup B_{R_0, R_{\beta}}^{0,\beta}|
\},
\quad
x\in {\Bbb R^2}.
\end{equation}
Taking
$$
D^{\theta,\theta'}_{R, R'} 
:= \left ( 
\begin{array} {cc} 
R \cos \theta & R' \cos \theta'
\\
R \sin \theta & R' \sin \theta'
\end{array} 
\right )
$$ 
and
$$
A^{\theta,\theta'}_{R, R'} 
:= \left ( 
\begin{array} {cc} 
R' \sin \theta' & -R' \cos \theta'
\\ 
- R \sin \theta & R \cos \theta
\end{array} 
\right ),
$$
for $\theta, \theta' \in [0, \pi)$,
$R, R' >0$,
we have
$B^{\alpha,\beta}_{R_\alpha, R_\beta}
= D^{\alpha,\beta}_{R_\alpha, R_\beta} [-1,1]^2$,
and
$$
{ D^{\alpha,\beta}_{R_\alpha, R_\beta}}^{-1}
=\frac{1}{\sin (\beta-\alpha)R_\alpha R_\beta }
A^{\alpha,\beta}_{R_\alpha, R_\beta}.
$$
In this notation we have
\begin{equation}
\left(
\begin{array}{c}
h_\alpha (x) \\ h_\beta (x)
\end{array}
\right )
=
{ D^{\alpha,\beta}_{\sin \beta, \sin \alpha} }^{-1} x
= \frac{1}{C_{\alpha,\beta}}
\left ( 
\begin{array} {cc} 
\sin \alpha \langle x, e_{\beta-\frac{\pi}{2}} \rangle
\\ 
\sin \beta \langle x, e_{\alpha + \frac{\pi}{2}} \rangle
\end{array} 
\right )
\end{equation}
where $h_\alpha$ and $h_\beta$ are as defined prior to Lemma \ref{lemma3.1}.
Note that 
$$
(h_\alpha (x), h_\beta (x))
\in [-H_\alpha, H_\alpha]\times [-H_\beta,H_\beta],
\hbox{ if and only if }
x \in B^{\alpha,\beta}_{R_\alpha, R_\beta},
$$
and
$$
\overline{h}_0 (x) := 
\frac{\langle x, e_{\frac{\pi}{2}}\rangle}{\sin \alpha \sin \beta}
= h_\alpha (x) + h_\beta (x),
\qquad x\in {\Bbb R}^2.
$$
See Figure 4.


\begin{figure}[h]
\label{Fig4}
\begin{center}
\unitlength 0.1in
\begin{picture}( 35.2000, 17.3100)(  1.1000,-19.9100)
%
\special{pn 13}%
\special{pa 1462 1030}%
\special{pa 3630 1030}%
\special{da 0.070}%
%
\special{pn 13}%
\special{pa 1462 1990}%
\special{pa 3630 1990}%
\special{da 0.070}%
%
\special{pn 13}%
\special{pa 3382 1030}%
\special{pa 3382 390}%
\special{fp}%
\special{sh 1}%
\special{pa 3382 390}%
\special{pa 3362 458}%
\special{pa 3382 444}%
\special{pa 3402 458}%
\special{pa 3382 390}%
\special{fp}%
%
\special{pn 13}%
\special{pa 3382 390}%
\special{pa 3382 1030}%
\special{fp}%
\special{sh 1}%
\special{pa 3382 1030}%
\special{pa 3402 964}%
\special{pa 3382 978}%
\special{pa 3362 964}%
\special{pa 3382 1030}%
\special{fp}%
%
\special{pn 13}%
\special{pa 3382 1990}%
\special{pa 3382 1030}%
\special{fp}%
\special{sh 1}%
\special{pa 3382 1030}%
\special{pa 3362 1098}%
\special{pa 3382 1084}%
\special{pa 3402 1098}%
\special{pa 3382 1030}%
\special{fp}%
%
\special{pn 13}%
\special{pa 3382 1030}%
\special{pa 3382 1990}%
\special{fp}%
\special{sh 1}%
\special{pa 3382 1990}%
\special{pa 3402 1924}%
\special{pa 3382 1938}%
\special{pa 3362 1924}%
\special{pa 3382 1990}%
\special{fp}%
%
\special{pn 13}%
\special{pa 1302 1990}%
\special{pa 1302 390}%
\special{fp}%
\special{sh 1}%
\special{pa 1302 390}%
\special{pa 1282 458}%
\special{pa 1302 444}%
\special{pa 1322 458}%
\special{pa 1302 390}%
\special{fp}%
%
\special{pn 13}%
\special{pa 1302 390}%
\special{pa 1302 1990}%
\special{fp}%
\special{sh 1}%
\special{pa 1302 1990}%
\special{pa 1322 1924}%
\special{pa 1302 1938}%
\special{pa 1282 1924}%
\special{pa 1302 1990}%
\special{fp}%
%
\special{pn 13}%
\special{pa 1622 1990}%
\special{pa 3062 1030}%
\special{fp}%
\special{sh 1}%
\special{pa 3062 1030}%
\special{pa 2996 1050}%
\special{pa 3018 1060}%
\special{pa 3018 1084}%
\special{pa 3062 1030}%
\special{fp}%
%
\special{pn 13}%
\special{pa 3054 1022}%
\special{pa 2414 382}%
\special{fp}%
\special{sh 1}%
\special{pa 2414 382}%
\special{pa 2448 444}%
\special{pa 2452 420}%
\special{pa 2476 416}%
\special{pa 2414 382}%
\special{fp}%
%
\special{pn 13}%
\special{ar 1622 1990 320 320  5.6744960 6.2831853}%
\put(16.2200,-21.5000){\makebox(0,0)[lb]{0}}%
\put(19.8200,-19.2600){\makebox(0,0)[lb]{$\alpha$}}%
\put(31.1000,-8.7800){\makebox(0,0)[lb]{$\beta$}}%
\put(22.7800,-4.3000){\makebox(0,0)[lb]{$x$}}%
\put(34.7800,-7.9800){\makebox(0,0)[lb]{$h_\beta(x)\sin\alpha\sin\beta$}}%
\put(34.8600,-15.8200){\makebox(0,0)[lb]{$h_\alpha(x)\sin\alpha\sin\beta$}}%
\put(1.1000,-12.5400){\makebox(0,0)[lb]{$\overline{h_0}(x)\sin\alpha\sin\beta$}}%
%
\special{pn 13}%
\special{ar 3014 1102 244 244  4.3121483 5.9822651}%
%
\special{pn 8}%
\special{ar 3040 2864 1682 1682  3.6870677 4.0771133}%
%
\special{pn 8}%
\special{ar 3220 2714 1682 1682  4.2135113 4.6042713}%
\put(22.5000,-13.7400){\makebox(0,0){$x^\alpha$}}%
%
\special{pn 8}%
\special{ar 3120 324 714 714  2.5332676 3.0433187}%
%
\special{pn 8}%
\special{ar 3100 304 714 714  1.6221462 2.1323665}%
\put(26.1000,-8.2400){\makebox(0,0){$x^\beta$}}%
\end{picture}%
\end{center}
\caption{\it{The quantities $h_\alpha$, $h_\beta$ and $\overline{h}_0$.}}
\end{figure}
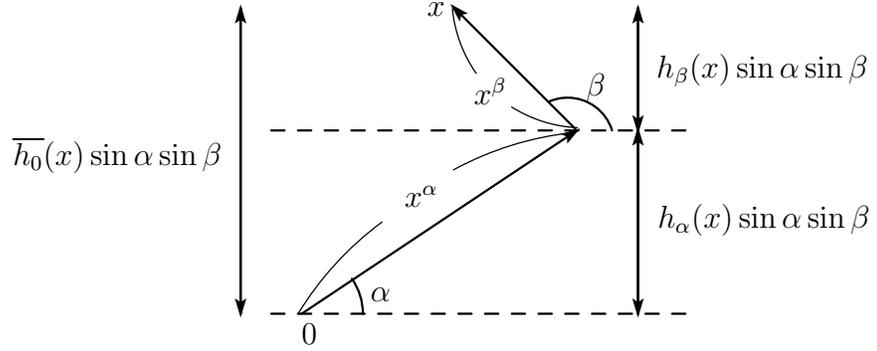
\begin{lemma}
\label{lemma4.2}
Assume that $x \in {\Bbb R}^2$ with
$h_\alpha (x) \in [-H_\alpha, H_\alpha]$,
$h_\beta (x)\in [-H_\beta, H_\beta ]$.

\noindent
{\rm (i)}~ Suppose that
$2H_0 < H_\alpha, H_{\beta}$. Then
\begin{eqnarray*}
\triangle (x)
& = & 
\frac{1}{2} 
\max \{ -h_\alpha(x)+2H_0-H_\alpha, h_\beta(x) + 2H_0 -H_\beta, 0\}^2
\\
&+& 
\frac{1}{2}
\max \{ h_\alpha(x)+2H_0-H_\alpha, -h_\beta(x)+ 2H_0 -H_\beta, 0\}^2.
\end{eqnarray*}

\noindent
{\rm (ii)}~ Suppose that
$2H_0 \ge \min\{H_\alpha, H_{\beta}\}$ and
$H_{\alpha} \ge H_{\beta}$.
\\
{\rm (a)}~When 
$|\overline{h}_0 (x)| \le H_\alpha -H_\beta$,
\begin{equation*}
\triangle(x) =
\begin{cases}
h_\beta (x)^2, & \text{if $|h_\beta (x)| \le 2H_0 - H_\beta$},
\\
h_\beta (x)^2 - \frac{1}{2} \{ |h_\beta (x)| - (2H_0 - H_\beta)\}^2,
& \text{if $|h_\beta (x)| >  2H_0 - H_\beta$}.
\end{cases}
\end{equation*}
\\
{\rm (b)}~When 
$|\overline{h}_0 (x)| > H_\alpha -H_\beta$
and $|h_\beta (x)| \le 2H_0 - H_\beta$,
\begin{eqnarray*}
& &\triangle(x) =h_\beta(x)^2 
+\frac{1}{2}\{ |\overline{h}_0(x)|-(H_\alpha - H_\beta)\}^2
\\
& & \: + 
\{ 2H_0 - H_\beta -\mathrm{sgn}(\overline{h}_0(x))h_\beta(x)\}
\{ |\overline{h}_0(x)|- (H_\alpha - H_\beta)\}.
\end{eqnarray*}
\\
{\rm (c)}~When 
$|\overline{h}_0(x)| > H_{\alpha}-H_{\beta}$,
$|h_\beta (x)| > 2H_0 - H_\beta$ and
$\overline{h}_0(x) h_\beta (x)>0$,
\begin{equation*}
\begin{split}
\triangle(x) 
&= h_\beta(x)^2 -\frac{1}{2}\{ |h_\beta (x)|-(2H_0 - H_\beta)\}^2
\\
&+\frac{1}{2}[2H_0 - H_{\alpha}+ \mathrm{sgn}(h_\beta(x))h_\alpha(x)]_+^2,
\end{split}
\end{equation*}
where $[a]_+ = \max \{a, 0\}$, $[a]_- = \max \{-a, 0\}$.
\\
{\rm (d)}~When 
$|\overline{h}_0(x)| > H_{\alpha}-H_{\beta}$, 
$|h_\beta (x)| > 2H_0 - H_\beta$ and $\overline{h}_0(x)h_\beta (x)<0$,
\begin{eqnarray*}
\triangle(x) 
&=&
h_\beta (x)^2 - \frac{1}{2} \{ |h_\beta (x)| - (2H_0 - H_\beta) \}^2
\\
&+&
\{ |\overline{h}_0(x)| - (H_\alpha - H_\beta)\}
\\
&\times&
[ 2H_0 - H_\beta + |h_\beta (x)| + \frac{1}{2}\{|\overline{h}_0(x)|-(H_\alpha-H_\beta)\}].
\end{eqnarray*}
\end{lemma}

\noindent
{\bf Remark 4.1.} \:
The area $\{ x\in {\Bbb R^2} : \triangle(x) =0 \}$ depends on angles 
$\alpha, \beta$ and stick lengths $R_0, R_\alpha,R_\beta$.
>From the above lemma we see that
\begin{equation}
\{ x \in {\Bbb R^2} : \triangle(x) = 0 \}
= B^{\alpha,\beta}_{R_\alpha -2R_0^\alpha, R_\beta -2R_0^\beta},
\quad\text{when $2H_0 < H_\alpha, H_\beta$},
\end{equation}
and 
\begin{equation}
\{ x \in {\Bbb R^2} : \triangle (x) = 0 \}
=B^{\alpha,\beta}_{[R_\alpha -R_\beta^\alpha]_+, [R_\beta -R_\alpha^\beta]_+},
\quad\text{when $2H_0 \ge \min\{H_\alpha, H_\beta\}$},
\end{equation}
where for $\theta = 0,\alpha,\beta$,
$R_\theta^0 = H_\theta \sin (\beta-\alpha)$,
$R_\theta^\alpha = H_\theta \sin \beta$,
$R_\theta^\beta = H_\theta \sin \alpha$.
In particular $R_\theta^\theta = R_\theta$.

Since
$$
A^{\alpha,\beta}_{R_\alpha, R_\beta} x
= \left ( 
\begin{array} {cc} 
R_\beta \langle x, e_{\beta-\frac{\pi}{2}} \rangle
\\ 
R_\alpha  \langle x, e_{\alpha + \frac{\pi}{2}} \rangle
\end{array} 
\right ),
$$
we have
\begin{eqnarray*}
&&M (A^{\alpha,\beta}_{R_\alpha, R_\beta} {\bf x}_k (0))
= R_\beta M ({\bf x}_k(\beta-\frac{\pi}{2} ))
=C_{\alpha,\beta} H_\beta M(h_\alpha ({\bf x}_k)),
\\
&&M (A^{\alpha,\beta}_{R_\alpha, R_\beta} {{\bf x}_k}(\frac{\pi}{2}))
= R_\alpha M ({\bf x}_k(\alpha + \frac{\pi}{2}))
=C_{\alpha,\beta} H_\alpha M(h_\beta ({\bf x}_k)).
\end{eqnarray*}
For ${\bf x}_k \in {\Bbb R^2}^k$, ${\bf y}_\ell \in {\Bbb R^2}^\ell$
and $u \in {\Bbb R^2}$ we write
$$
{\bf x}_k \cdot {\bf y}_\ell 
= (x_1, x_2, \dots x_k, y_1, y_2, \dots,y_\ell)
\in ({\Bbb R^2})^{k+\ell},
$$
and
${\bf x}_k + u = (x_1+u, x_2+u, \dots, x_k +u) \in ({\Bbb R^2})^k$.
We put 
$$
\triangle ({\bf x}_k, {\bf y}_\ell |u)
=
\frac{1}{C_{\alpha,\beta}}\{
|B^{0, \alpha}_{R_0, R_\alpha} ({\bf x}_k)
\cup B^{0,\beta}_{R_0, R_\beta} ({\bf y}_\ell+u)|
-|B^{0, \alpha}_{R_0, R_\alpha} \cup B^{0,\beta}_{R_0, R_\beta}(u)|\},
$$
and write $\triangle ({\bf x}_k, {\bf y}_\ell)$
for $\triangle ({\bf x}_k, {\bf y}_\ell |{\bf 0})$.
The following two lemmas are important to show the main theorem. Their proofs are given in the appendix.

\begin{lemma}
\label{lemma4.3}
Let ${\bf x}_k \in (\Bbb R^2)^k$ with $x_k = {\bf 0}$ and 
${\bf y}_\ell \in (\Bbb R^2)^\ell$ with $y_\ell = {\bf 0}$.  
\\
{\rm (i)}~ Suppose that $2H_0 < H_\alpha, H_\beta$. If
\begin{eqnarray}
\label{llcond1}
M(h_\alpha({\bf x}_k)) + M(h_\alpha({\bf y}_\ell))< H_\alpha -2H_0 \quad  \mbox{ and}\nonumber\\
 M(h_\beta({\bf x}_k)) + M(h_\beta({\bf y}_\ell))< H_\beta -2H_0
\end{eqnarray}
hold, then we have
\begin{equation*}
\triangle ({\bf x}_k, {\bf y}_\ell)
\le \frac{1}{C_{\alpha, \beta}}
\{ |B^{0, \alpha}_{R_0, R_\alpha - R_0^\alpha} ({\bf x}_k)
\setminus B^{0, \alpha}_{R_0, R_\alpha - R_0^\alpha}|
+|B^{0,\beta}_{R_0, R_\beta - R_0^\beta} ({\bf y}_\ell)
\setminus B^{0,\beta}_{R_0, R_\beta - R_0^\beta}|\},
\end{equation*}
\begin{eqnarray*}
\triangle ({\bf x}_k, {\bf y}_\ell)
&\ge& \frac{1}{C_{\alpha, \beta}}
\{|B^{0, \alpha}_{R_0, R_\alpha - R_0^\alpha} ({\bf x}_k)
\setminus B^{0, \alpha}_{R_0, R_\alpha - R_0^\alpha}|
+|B^{0,\beta}_{R_0, R_\beta - R_0^\beta} ({\bf y}_\ell)
\setminus B^{0,\beta}_{R_0, R_\beta - R_0^\beta}|\}
\\
&-& M(h_\alpha({\bf y}_\ell))M(h_\beta({\bf x}_k)).
\end{eqnarray*}
{\rm (ii)}~ Suppose that $2H_0 \ge \min\{H_\alpha, H_\beta\}$ 
and $H_\alpha > H_\beta$. If
$M(h_\alpha({\bf x}_k)) + M(h_\alpha({\bf y}_\ell))< H_\alpha -H_\beta$ and
$M(h_\beta({\bf x}_k)) + M(h_\beta({\bf y}_\ell))< H_\beta$ hold,
then we have
\begin{eqnarray}
\label{r-tricond}
\triangle ({\bf x}_k, {\bf y}_\ell)
&\le&\frac{1}{C_{\alpha, \beta}}
\{|B^{0, \alpha}_{R_0, R_\alpha - \frac{1}{2}R_\beta^\alpha} ({\bf x}_k)
\setminus
B^{0, \alpha}_{R_0, R_\alpha - \frac{1}{2}R_\beta^\alpha}|
+|B^{0,\beta}_{\frac{1}{2}R_\beta^0, \frac{1}{2}R_\beta} ({\bf y}_\ell)
\setminus
B^{0,\beta}_{\frac{1}{2}R_\beta^0, \frac{1}{2}R_\beta}|\}\nonumber
\\
&+&\frac{1}{2}M(h_\beta({\bf x}_k))^2+\frac{1}{2}M(h_\alpha({\bf y}_\ell))^2,
\end{eqnarray}
\begin{eqnarray*}
\triangle ({\bf x}_k, {\bf y}_\ell)
&\ge&
\frac{1}{C_{\alpha, \beta}}
\{|B^{0, \alpha}_{R_0, R_\alpha - \frac{1}{2}R_\beta^\alpha} ({\bf x}_k)
\setminus
B^{0, \alpha}_{R_0, R_\alpha - \frac{1}{2}R_\beta^\alpha}|
+|B^{0,\beta}_{\frac{1}{2}R_\beta^0, \frac{1}{2}R_\beta} ({\bf y}_\ell)
\setminus
B^{0,\beta}_{\frac{1}{2}R_\beta^0, \frac{1}{2}R_\beta}|\}
\\
&-& M(h_\beta(x_k))M(h_\beta(y_\ell)) - 
M(h_\beta(x_k))M(h_\alpha(y_\ell)) \\
&-&
(M(h_\beta(x_k)))^2 - (M(h_\alpha(y_\ell)))^2.
\end{eqnarray*}
{\rm (iii)}~ Suppose that $2H_0 \ge H_\alpha = H_\beta$. If
$M(h_\alpha({\bf x}_k)) + M(h_\alpha({\bf y}_\ell))< H_\alpha$ and
$M(h_\beta({\bf x}_k)) + M(h_\beta({\bf y}_\ell))< H_\beta$ hold,
 then we have
\begin{eqnarray*}
\triangle ({\bf x}_k, {\bf y}_\ell)
&\le&
\frac{1}{C_{\alpha, \beta}}
\{|B^{0, \alpha}_{\frac{1}{2}R_\alpha^0, \frac{1}{2}R_\alpha}({\bf x}_k)
\setminus
B^{0, \alpha}_{\frac{1}{2}R_\alpha^0, \frac{1}{2}R_\alpha}|
+|B^{0,\beta}_{\frac{1}{2}R_\beta^0, \frac{1}{2}R_\beta}({\bf y}_\ell)
\setminus B^{0,\beta}_{\frac{1}{2}R_\beta^0, \frac{1}{2}R_\beta}|\}
\\
&+& (2H_0 -H_\beta)M(\overline{h}_0({\bf x}_k \cdot {\bf y}_\ell))
+\frac{1}{2}M({h}_\beta({\bf x}_k))^2+\frac{1}{2}M({h}_\alpha({\bf y}_\ell))^2,
\end{eqnarray*}
and
\begin{eqnarray}
\label{r-trapcond}
\triangle ({\bf x}_k, {\bf y}_\ell)
&\ge&\frac{1}{C_{\alpha, \beta}}
\{ |B^{0, \alpha}_{\frac{1}{2}R_\alpha^0, \frac{1}{2}R_\alpha}({\bf x}_k)
\setminus
B^{0, \alpha}_{\frac{1}{2}R_\alpha^0, \frac{1}{2}R_\alpha}|
+|B^{0,\beta}_{\frac{1}{2}R_\beta^0, \frac{1}{2}R_\beta}({\bf y}_\ell)
\setminus B^{0,\beta}_{\frac{1}{2}R_\beta^0, \frac{1}{2}R_\beta}|\}
\nonumber \\
&+&(2H_0 -H_\beta )M(\overline{h}_0({\bf x}_k \cdot {\bf y}_\ell))
-\frac{1}{2}M(\overline{h}_0({\bf x}_k \cdot {\bf y}_\ell))^2
\nonumber \\
&-& \min \{ M(\overline{h}_0({\bf x}_k)),M(\overline{h}_0({\bf y}_\ell))\}
\{M({h}_\beta({\bf x}_k))+M({h}_\alpha ({\bf y}_\ell))\}.
\end{eqnarray}
\end{lemma}
\begin{lemma}
\label{lemma4.4}
Let ${\bf x}_k \in (\Bbb R^2)^k$ with $x_k = {\bf 0}$,
${\bf y}_\ell \in (\Bbb R^2)^\ell$ with $y_\ell = {\bf 0}$
and $u \in \Bbb R^2$.
\\
{\rm (i)}~ Suppose that
$2H_0 < H_\alpha, H_\beta$. If
$M(h_\alpha({\bf x}_k)) + M(h_\alpha({\bf y}_\ell))+ |h_\alpha (u)| 
< H_\alpha - 2H_0$ and 
$M(h_\beta({\bf x}_k)) + M(h_\beta({\bf y}_\ell))+ |h_\beta (u)| 
< H_\beta - 2H_0$ hold,
then we have
\begin{equation*}
\triangle ({\bf x}_k, {\bf y}_\ell |u)
=\triangle ({\bf x}_k, {\bf y}_\ell).
\end{equation*}
{\rm (ii)}~ Suppose that $2H_0 \ge \min\{H_\alpha, H_\beta\}$ 
and $H_\alpha > H_\beta$. If
$M(h_\alpha({\bf x}_k)) + M(h_\alpha({\bf y}_\ell))+ |h_\alpha (u)| 
< H_\alpha - H_\beta$ and
$M(h_\beta({\bf x}_k)) + M(h_\beta({\bf y}_\ell))+ |h_\beta (u)| 
< H_\beta$ hold,
then we have
\begin{equation*}
\left| \triangle ({\bf x}_k, {\bf y}_\ell | u)
- \triangle ({\bf x}_k, {\bf y}_\ell) \right|
\le h_\beta(u)^2.
\end{equation*}
{\rm (iii)}~ Suppose that $2H_0 \ge H_\alpha = H_\beta$. If
$M(h_\alpha({\bf x}_k)) + M(h_\alpha({\bf y}_\ell))+ |h_\alpha (u)| 
< H_\alpha$ and
$M(h_\beta({\bf x}_k)) + M(h_\beta({\bf y}_\ell))+ |h_\beta (u)| 
< H_\beta$ hold,
 then we have
\begin{eqnarray*}
&&\left| \;
\triangle ({\bf x}_k, {\bf y}_\ell | u)
-\triangle ({\bf x}_k, {\bf y}_\ell )
\right.
\\
&&\qquad \left. - ( 2H_0 -H_\beta )
\{ M(\overline{h}_0({\bf x}_k \cdot ({\bf y}_\ell +u))) - |\overline{h}_0 (u)|
-M(\overline{h}_0({\bf x}_k \cdot {\bf y}_\ell)) \}
\; \right|
\\
&&\le 
h_\alpha(u)^2 + h_\beta(u)^2
+| M(\overline{h}_0({\bf x}_k \cdot ({\bf y}_\ell +u)))- |\overline{h}_0 (u)|
-M(\overline{h}_0({\bf x}_k \cdot {\bf y}_\ell))|
\\
&&\times
\{M(h_\alpha({\bf x}_k)) + M(h_\alpha({\bf y}_\ell))+ |h_\alpha (u)|
+ M(h_\beta({\bf x}_k)) + M(h_\beta({\bf y}_\ell))+ |h_\beta (u)| \},
\end{eqnarray*}
if
$M(h_\alpha({\bf x}_k)) + M(h_\alpha({\bf y}_\ell))+ |h_\alpha (u)| 
< H_\alpha$,
$M(h_\beta({\bf x}_k)) + M(h_\beta({\bf y}_\ell))+ |h_\beta (u)| 
< H_\beta$.
\end{lemma}

\subsection {The asymptotic shape}

First, we examine the behaviour of the function 
$\mu_{\lambda\rho}(C_{\bf 0}\in\Lambda({\bf k})| \Gamma_0)$
as $\lambda\to\infty$ when ${\bf k}=(0,k_\alpha,k_\beta)$.
When ${\bf k}=(k_0,k_\alpha,0)$ or ${\bf k}=(k_0,0,k_\beta)$,
we can estimate similarly.
>From (\ref{lsticks}) we have
\begin{equation}
\mu_{\lambda \rho}(C_{\bf 0}\in\Lambda (0, k_\alpha, k_\beta)\mid ~ 
\Gamma_0 )
= \lambda^{|{\bf k}| -1}|{\bf k}|
\frac{p_\alpha^{k_\alpha} p_\beta^{k_\beta}}{(k_\alpha)!k_\beta !}
F_\lambda (0, k_\alpha, k_\beta),
\end{equation}
where
\begin{eqnarray*}
F_\lambda (0,k_\alpha, k_\beta )
&=& \int\limits_{(\Bbb R^2)^{k_\alpha-1} } d{\bf y}_{k_\alpha-1}
\int\limits_{(\Bbb R^2)^{k_\beta}} d{\bf z}_{k_\beta}
1_{\Lambda(0,k_\alpha, k_\beta)} 
(C_{\bf 0} ( {\bf y}_{k_\alpha},{\bf z}_{k_\beta}))
\\
&\times& e^{ -\lambda \{
p_0 | B^{\alpha,0}_{R_\alpha,R_0}({\bf y}_{k_\alpha})
\cup B^{\beta, 0}_{R_\beta,R_0} ({\bf z}_{k_\beta})|
+ p_\alpha |B^{\beta,\alpha}_{R_\beta,R_\alpha}({\bf z}_{k_\beta})|
+ p_\beta |B^{\alpha, \beta}_{R_\alpha, R_\beta} ({\bf y}_{k_\alpha})|\}}.
\end{eqnarray*}
We put
\begin{eqnarray}
&&\Phi ({\bf p})
= p_0 | B^{\alpha, 0}_{R_\alpha,R_0} 
\cup B^{\beta, 0}_{R_\beta ,R_0}|
+ p_\alpha |B^{\beta, \alpha}_{R_\beta ,R_\alpha}|
+ p_\beta |B^{\alpha, \beta}_{R_\alpha ,R_\beta}|,
\\
&&f_\lambda (0, k_\alpha, k_\beta) 
= F_\lambda (0, k_\alpha, k_\beta)
e^{ \lambda \Phi ({\bf p})}.
\nonumber
\end{eqnarray}
To examine the function $f_\lambda ({\bf k})$,
we introduce the following function 
\begin{equation}
\chi^{\theta_1,\theta_2,\theta_3}_c ({\bf x}, {\bf y} |z) 
=
e^{-c\{| B^{\theta_1, \theta_2}_{R_{\theta_1},R_{\theta_2}} ({\bf x})
\cup B^{\theta_1, \theta_3}_{R_{\theta_1} ,R_{\theta_3}} ({\bf y}+z)|
-
| B^{\theta_1, \theta_2}_{R_{\theta_1},R_{\theta_2}}
\cup B^{\theta_1, \theta_3}_{R_{\theta_1} ,R_{\theta_3}} (z)|\}},
\end{equation}
for $\theta_1, \theta_2, \theta_3 \in [0,\pi)$,
$c>0$, ${\bf x} \in ({\Bbb R^2})^k$
${\bf y} \in ({\Bbb R^2})^{k'}, k, k'\in {\Bbb N}$ and $z \in {\Bbb R^2}$.
We write $\chi^{\theta_1, \theta_2, \theta_3}_c ({\bf x},{\bf y})$ for
$\chi^{\theta_1,\theta_2,\theta_3}_c ({\bf x}, {\bf y} |{\bf 0})$.
By using these functions we obtain 
\begin{eqnarray*}
f_\lambda (0, k_\alpha, k_\beta)
&=&
\int \limits _{(\Bbb R ^2)^{k_\alpha -1}} d{\bf y}_{k_\alpha-1}
\int \limits _{(\Bbb R ^2)^{k_\beta }} d{\bf z}_{k_\beta}
1_{\Lambda(0, k_\alpha, k_\beta)}
(C_{\bf 0} ({\bf y}_{k_\alpha},{\bf z}_{k_\beta }))
\\
&&\quad\times
\chi^{0,\alpha,\beta}_{\lambda p_0 }({\bf y}_{k_\alpha},{\bf z}_{k_\beta})
\chi^{\alpha,\beta}_{\lambda p_\alpha} ({\bf z}_{k_\beta})
\chi^{\alpha,\beta}_{\lambda p_\beta}({\bf y}_{k_\alpha}).
\end{eqnarray*}
Putting
${\bf u}_{k_\alpha} = {\bf y}_{k_\alpha} - y_{k_\alpha}$,
${\bf v}_{k_\beta} = {\bf z}_{k_\beta} - z_{k_\beta}$
and $z_\beta = z$,
we have
\begin{equation*}
f_\lambda (0, k_\alpha, k_\beta)
= \int \limits _{\Bbb R ^2} dz
g_\lambda (0, k_\alpha, k_\beta, z)
\chi_{\lambda p_0}^{0,\alpha,\beta} ({\bf 0}, z),
\end{equation*}
where
\begin{eqnarray}
g_\lambda (0, k_\alpha, k_\beta, z)
&=&\int \limits_{(\Bbb R ^2)^{k_\alpha-1} } d{\bf u}_{k_\alpha-1}
\int \limits_{(\Bbb R ^2)^{k_\beta-1}} d{\bf v}_{k_\beta -1}
1_{\Lambda(0, k_\alpha, k_\beta)} 
(C_{\bf 0} ({\bf u}_{k_\alpha}, {\bf v}_{k_\beta}+z))
\nonumber
\\
&\times&
\chi^{0,\alpha,\beta}_{\lambda p_0}({\bf u}_{k_\alpha}, 
{\bf v}_{k_\beta} |z )
\chi^{\alpha,\beta}_{\lambda p_\alpha} ({\bf v}_{k_\beta})
\chi^{\alpha,\beta}_{\lambda p_\beta} ({\bf u}_{k_\alpha}).
\end{eqnarray}
Writing $g_\lambda ({\bf k})$ for
$g_\lambda ({\bf k},{\bf 0})$,
we have
\begin{eqnarray}
\label{asym1}
&&\mu_{\lambda \rho}
( C_{\bf 0}\in\Lambda(0, k_\alpha, k_\beta)\mid ~ \Gamma_{\bf 0})
\nonumber
\\
&=& e^{-\lambda \Phi ({\bf p})}\lambda^{|{\bf k}| -1}|{\bf k}|
\frac{p_\alpha^{k_\alpha} p_\beta^{k_\beta}}{k_\alpha!k_\beta!}
\int \limits _{\Bbb R ^2} dz
 g_\lambda (0, k_\alpha, k_\beta, z)
\chi_{\lambda p_0}^{0,\alpha,\beta} ({\bf 0}, z).
\end{eqnarray}
{\bf Remark 4.2.}
The function $\chi_{\lambda p_0}^{0,\alpha,\beta}$
determines the structure of finite clusters.
>From Remark 4.1 we see that
$\chi_{\lambda p_0}^{0,\alpha,\beta}(0,z)
=\exp [-\lambda p_0 C_{\alpha,\beta}\triangle (z)]=1$ 
if and only if 
\begin{eqnarray*}
&&z \in B^{\alpha,\beta}_{R_\alpha -2R_0^\alpha, R_\beta-2R_0^\beta},
\qquad
\text{when $H_{\alpha}, H_{\beta} > 2H_0$},
\\
&&z \in B^{\alpha, \beta}_{[R_\alpha - R_\beta^\alpha]_+, 
[R_\beta - R_\alpha^\beta]_+},
\quad
\text{ when $\min\{ H_{\alpha}, H_{\beta}\} \le 2H_0$}.
\end{eqnarray*}
We divide into four cases and obtain  estimates.

{\bf Case (1)} $2H_0 < H_\alpha, H_\beta$.
In this case we will show that
\begin{eqnarray}
\label{case1}
&&\mu_{\lambda\rho}(C_{\bf 0}\in\Lambda(0,k_\alpha,k_\beta)| \Gamma_0)
\nonumber
\\
&&\qquad \sim
\exp[ -4C_{\alpha,\beta}\lambda
\{p_0 H_0(H_\alpha+H_\beta-H_0)+(1-p_0)H_\alpha H_\beta\}]
\nonumber
\\
&&\qquad \times \left( 
\frac{1}{4C_{\alpha,\beta}\lambda}\right)^{|{\bf k}|-3}
|{\bf k}|H_\alpha H_\beta(H_\alpha-2H_0)(H_\beta-2H_0)
\nonumber
\\
&&\qquad \times
p_{\alpha}^{k_\alpha} k_\alpha !
G^{k_\alpha}(p_0 H_0 + p_\beta H_\beta, p_\beta H_\alpha, p_0(H_\alpha -H_0))
\nonumber
\\
&& \qquad \times
p_{\beta}^{k_\beta} k_\beta !
G^{k_\beta}(p_\alpha H_\beta, p_0 H_0 + p_\alpha H_\alpha, p_0(H_\beta -H_0)),
\end{eqnarray}
where for $c_1, c_2, c_3 >0$
\begin{eqnarray}
G^k(c_1,c_2,c_3)
&=& (\frac{1}{k!})^2 \int \limits_{(\Bbb R ^2)^{k-1}} d{\bf u}_{k-1}
\gamma^k(c_1,c_2,c_3)({\bf u}_k),
\\
\gamma^k(c_1, c_2, c_3)({\bf u}_k)
&=&\exp [-\{ c_1 M({\bf u}^1_k)+c_2 M({\bf u}^2_k)
+c_3M({\bf u}^1_k+{\bf u}^2_k)\}].
\end{eqnarray}
>From Remark 4.2 we see that the asymptotic shape of the cluster is given by
\begin{equation*}
\{ x\in {\Bbb R^2} : 
|h_\alpha (x)| \le H_\alpha - 2H_0, \; 
|h_\beta (x) | \le H_\beta - 2H_0 \}.
\end{equation*}
By Lemma \ref{lemma4.2} {\rm (i)} and Lemma \ref{lemma4.4} {\rm (i)} we have
\begin{equation}
f_\lambda (0,k_\alpha,k_\beta) 
\sim 
|B^{\alpha, \beta}_{R_\alpha - 2 R_0^\alpha, R_\beta - 2R_0^\beta}|
g_\lambda(0,k_\alpha,k_\beta),
\quad\text{as $\lambda\to\infty$}.
\end{equation}
By Lemma \ref{lemma4.3} {\rm (i)} we have
\begin{eqnarray*}
&&g_\lambda (0, k_\alpha, k_\beta)
\\
&\sim& \int \limits_{(\Bbb R ^2)^{k_\alpha-1} } d{\bf u}_{k_\alpha-1}
e^{-\lambda\{ p_0 |B^{0,\alpha}_{R_0,R_\alpha - R^\alpha_0}({\bf u}_{k_\alpha})
\setminus B^{0,\alpha}_{R_0,R_\alpha - R^\alpha_0}|
+p_\beta |B^{\alpha,\beta}_{R_\alpha,R_\beta}({\bf u}_{k_\alpha})
\setminus B^{\alpha,\beta}_{R_\alpha,R_\beta}|\}}
\\
&\times& \int \limits_{(\Bbb R ^2)^{k_\beta-1}} d{\bf v}_{k_\beta -1}
e^{-\lambda \{ p_0 |B^{0,\beta}_{R_0,R_\beta-R_0^\beta}({\bf v}_{k_\beta})
\setminus B^{0,\beta}_{R_0,R_\beta - R_0^\beta}|
+p_\alpha |B^{\alpha,\beta}_{R_\alpha,R_\beta}({\bf v}_{k_\beta})
\setminus B^{\alpha,\beta}_{R_\alpha,R_\beta}|\}}
\end{eqnarray*}
Using Lemma \ref{lemma3.1} and putting 
${\bf \hat{u}}=A_{2\lambda\sin\beta,2\lambda\sin\alpha}^{\alpha,\beta}{\bf u}$,
by a simple calculation we have
\begin{eqnarray*}
&&\int \limits_{(\Bbb R ^2)^{k_\alpha-1} } d{\bf u}_{k_\alpha-1}
e^{-\lambda\{ p_0 |B^{0,\alpha}_{R_0,R_\alpha - R^\alpha_0}({\bf u}_{k_\alpha})
\setminus B^{0,\alpha}_{R_0,R_\alpha - R^\alpha_0}|
+p_\beta |B^{\alpha,\beta}_{R_\alpha,R_\beta}({\bf u}_{k_\alpha})
\setminus B^{\alpha,\beta}_{R_\alpha,R_\beta}|\}}
\\
&\sim&\int\limits_{(\Bbb R ^2)^{k_\alpha-1} } d{\bf u}_{k_\alpha-1}
e^{-2C_{\alpha,\beta}\lambda
[(p_0 H_0 + p_\beta H_\beta) M(h_\alpha({\bf u}_{k_\alpha}))
+ p_0 (H_\alpha - H_0)M(\overline{h}_0({\bf u}_{k_\alpha}))
+ p_\beta H_\beta M(h_\alpha({\bf u}_{k_\alpha}))]}
\\
&=&\left(\frac{1}{4C_{\alpha,\beta}\lambda^2}\right)^{k_\alpha-1}
G^{k_\alpha}(p_0 H_0 + p_\beta H_\beta, p_\beta H_\alpha, p_0 (H_\alpha -H_0)).
\end{eqnarray*}
Similarly, we have
\begin{eqnarray*}
&&\int \limits_{(\Bbb R ^2)^{k_\beta-1}} d{\bf v}_{k_\beta -1}
e^{-\lambda \{ p_0 |B^{0,\beta}_{R_0,R_\beta-R_0^\beta}({\bf v}_{k_\beta})
\setminus B^{0,\beta}_{R_0,R_\beta - R_0^\beta}|
+p_\alpha |B^{\alpha,\beta}_{R_\alpha,R_\beta}({\bf v}_{k_\beta})
\setminus B^{\alpha,\beta}_{R_\alpha,R_\beta}|\}}
\\
&\sim&
\left(\frac{1}{4C_{\alpha,\beta}\lambda^2}\right)^{k_\beta-1}
G^{k_\beta}(p_\alpha H_\beta,p_0 H_0 + p_\alpha H_\alpha,p_0(H_\beta-H_0))
\end{eqnarray*}
Since by Lemma \ref{lemma4.1} (i)
\begin{equation}
\Phi ({\bf p}) = 4C_{\alpha,\beta}
\{ p_0 H_0 (H_\alpha + H_\beta -H_0) + (1-p_0) H_\alpha H_\beta \},
\end{equation}
we have (\ref{case1}) from (\ref{asym1}) and the above estimates.

{\bf Case (2)} $2H_0 \ge H_\beta$, $H_\alpha > H_\beta $.
In this case we will show that
\begin{eqnarray}
\label{case2}
&&\mu_{\lambda\rho}(C_{\bf 0}\in\Lambda(0,k_\alpha,k_\beta)| \Gamma_0)
\nonumber
\\
&&\qquad \sim\exp[-4 C_{\alpha,\beta}\lambda
\{p_0 H_0 H_\alpha + \frac{p_0}{4}H_\beta^2 +(1-p_0)H_\alpha H_\beta\}]
\nonumber
\\
&&\qquad \times 
\left(\frac{1}{4C_{\alpha, \beta}\lambda}\right)^{|{\bf k}|-\frac{5}{2}}
|{\bf k}||H_\alpha - H_\beta|(\frac{\pi}{p_0})^{\frac{1}{2}}
\nonumber
\\
&&\qquad \times p_\alpha^{k_\alpha } k_\alpha !
G^{k_\alpha}
(p_0 H_0 + p_\beta H_\beta, p_\beta H_\alpha,
p_0 ( H_\alpha - \frac{1}{2} H_\beta))
\nonumber
\\
&&\qquad \times p_\beta^{k_\beta } k_\beta !
G^{k_\beta}
(p_\alpha H_\beta,  \frac{1}{2}p_0 H_\beta + p_\alpha H_\alpha,
\frac{1}{2}p_0 H_\beta).
\end{eqnarray}

>From Remark 4.2 we see that the asymptotic shape of the cluster is given by
\begin{equation*}
\{ x\in {\Bbb R^2} : 
|h_\alpha (x)| \le H_\alpha - H_\beta, \; 
|h_\beta (x) | =0  \}.
\end{equation*}
By Lemma \ref{lemma4.4} {\rm (ii)} and a simple calculation we have
\begin{equation*}
g_\lambda(0,k_\alpha,k_\beta,z)
\sim g_\lambda(0,k_\alpha,k_\beta)
\quad\text{as $\lambda\to\infty$, }
\end{equation*}
when 
$|h_\alpha (z)|< H_\alpha - H_\beta$
$|h_\beta(z)|= o(1)$.
>From Lemma \ref{lemma4.2} {\rm (ii)} we have
\begin{equation*}
\chi_{\lambda p_0}^{0,\alpha.\beta} ({\bf 0}, z)
= e^{-p_0 C_{\alpha,\beta}\lambda h_\beta(z)^2},
\end{equation*}
if $|\overline{h}_0(z)|\le H_\alpha-H_\beta$, $|h_\beta(z)|\le 2H_0-H_\beta$.
Then we have
\begin{eqnarray}
f_\lambda (0,k_\alpha,k_\beta) 
&\sim&
g_\lambda(0,k_\alpha,k_\beta)\int \limits _{\Bbb R ^2} dz
\chi_{\lambda p_0}^{0,\alpha,\beta} ({\bf 0}, z)
\nonumber
\\
&\sim&
2|H_\alpha-H_\beta|(\frac{C_{\alpha,\beta}\pi}{p_0 \lambda})^{1/2}
g_\lambda(0,k_\alpha,k_\beta)\
\quad\text{as $\lambda\to\infty$}.
\end{eqnarray}
By Lemma \ref{lemma3.1} and Lemma \ref{lemma4.3} {\rm (ii)} and 
similar calculations as above, we have
\begin{eqnarray*}
g_\lambda (0, k_\alpha, k_\beta)
&\sim&
\left(\frac{1}{4C_{\alpha, \beta}\lambda^2}\right)^{k_\alpha -1}
G^{k_\alpha}
(p_0 H_0 + p_\beta H_\beta, p_\beta H_\alpha, 
p_0 (H_\alpha - \frac{1}{2}H_\beta))
\\
&\times&
\left(\frac{1}{4C_{\alpha, \beta}\lambda^2}\right)^{k_\beta -1}
G^{k_\beta}
(p_\alpha H_\beta, \frac{1}{2}p_0 H_\beta + p_\alpha H_\alpha,
\frac{1}{2}p_0 H_\beta).
\end{eqnarray*}
Since by Lemma \ref{lemma3.1} (ii)
\begin{equation}
\Phi ({\bf p}) = 4C_{\alpha,\beta}
\{ p_0 H_0 H_\alpha + \frac{p_0}{4}H_\beta^2 + (1-p_0) H_\alpha H_\beta \},
\end{equation}
we have (\ref{case2}) from (\ref{asym1}) and the above estimates

{\bf Case (3)} $2H_0=H_\alpha =H_\beta$.
In this case we will show that
\begin{eqnarray}
\label{case3}
&&\mu_{\lambda\rho}(C_{\bf 0}\in\Lambda(0,k_\alpha,k_\beta)| \Gamma_0)
\nonumber
\\
&&\qquad
\sim\exp[-4 C_{\alpha,\beta}\lambda (4-p_0 )H_0^2]
\nonumber
\\
&&\qquad 
\times
\left(\frac{1}{4C_{\alpha, \beta}\lambda}\right)^{|{\bf k}|-2}
|{\bf k}|\frac{3\pi +4}{2p_0}
\nonumber
\\
&&\qquad
\times p_\alpha^{k_\alpha } k_\alpha !
G^{k_\alpha}
((p_0 + 2p_\beta) H_0, 2 p_\beta H_0, p_0 H_0)
\nonumber
\\
&&\qquad
\times p_\beta^{k_\beta } k_\beta !
G^{k_\beta}
(2p_\alpha H_0, (p_0 + 2 p_\alpha) H_0, p_0 H_0).
\end{eqnarray}

>From Remark 4.2 we see that the asymptotic shape of the cluster is given by
\begin{equation*}
\{ x\in {\Bbb R^2} : 
|h_\alpha (x)| \le 0, \; 
|h_\beta (x) | \le 0 \} = \{ {\bf 0} \}.
\end{equation*}
By Lemma \ref{lemma4.4} {\rm (iii)} and a simple calculation we have
\begin{equation*}
g_\lambda(0,k_\alpha,k_\beta,z)
\sim g_\lambda(0,k_\alpha,k_\beta)
\quad\text{as $\lambda\to\infty$, }
\end{equation*}
when 
$|h_\alpha (z)|= o(1)$, $|h_\beta(z)|= o(1)$.
>From Lemma \ref{lemma4.2} {\rm (ii)} we have
\begin{equation}
\chi_{\lambda p_0}^{0,\alpha,\beta} ({\bf 0}, z)=
\begin{cases}
\exp [ - \frac{1}{2}C_{\alpha, \beta}p_0 \lambda
(h_\alpha(z)^2 + h_\beta(z)^2 ) ],
&\quad \overline{h}_0(z)h_\beta (z)>0,
\\
\exp [ - \frac{1}{2}C_{\alpha, \beta}p_0 \lambda h_\alpha(z)^2 ],
&\quad \overline{h}_0(z)h_\beta (z)>0,
\end{cases}
\end{equation}
if $|h_\alpha(z)|\le H_\alpha$, $|h_\beta(z)|\le H_\beta$.
Then we have
\begin{eqnarray}
f_\lambda (0,k_\alpha,k_\beta) 
&\sim&
g_\lambda(0,k_\alpha,k_\beta)\int \limits_{\Bbb R ^2} dz
\chi_{\lambda p_\beta}^{0,\alpha,\beta} ({\bf 0}, z)
\nonumber
\\
&\sim&
(\frac{3\pi + 4}{2p_0 \lambda})
g_\lambda(0,k_\alpha,k_\beta)
\quad\text{as $\lambda\to\infty$}.
\end{eqnarray}
By Lemma \ref{lemma3.1} and Lemma \ref{lemma4.3} {\rm (iii)} and 
similar calculations as above, we have
\begin{eqnarray*}
g_\lambda (0, k_\alpha, k_\beta)
&\sim&
\left(\frac{1}{4C_{\alpha, \beta}\lambda^2}\right)^{k_\alpha -1}
G^{k_\alpha}
(\frac{1}{2} p_0 H_\alpha + p_\beta H_\beta, p_\beta H_\alpha, 
\frac{1}{2}p_0 H_\alpha)
\\
&\times&
\left(\frac{1}{4C_{\alpha, \beta}\lambda^2}\right)^{k_\beta -1}
G^{k_\beta}
(p_\alpha H_\beta, \frac{1}{2}p_0 H_\beta + p_\alpha H_\alpha,
\frac{1}{2}p_0 H_\beta).
\end{eqnarray*}
Since by Lemma \ref{lemma3.1} (ii),
$\Phi ({\bf p}) = 4C_{\alpha,\beta}(4 - p_0)H_0^2)$,
we have (\ref{case3}) from (\ref{asym1}) and the above estimates

{\bf Case (4)} $2H_0 > H_\alpha =H_\beta $.
In this case we will show that 
\begin{eqnarray}
\label{case4}
&&\mu_{\lambda\rho}
(C_{\bf 0}\in\Lambda(0,k_\alpha,k_\beta)| \Gamma_0)
\nonumber
\\
&&\qquad \sim
\exp[-4 C_{\alpha,\beta}\lambda
\{p_0 H_0 H_\alpha + (1-\frac{3}{4}p_0)H_\alpha^2\}]
\nonumber
\\
&&\qquad \times
\left(\frac{1}{4C_{\alpha, \beta}\lambda}\right)^{|{\bf k}|-\frac{3}{2}}
|{\bf k}|
(\frac{2 \pi}{p_0})^{\frac{1}{2}}
p_\alpha^{k_\alpha } k_\alpha !
p_\beta^{k_\beta } k_\beta !
\\
&&\qquad \times
G_{\frac{1}{2}(2H_0 - H_\alpha)}^{k_\alpha, k_\beta}
((\frac{p_0}{2} + p_\beta) H_\alpha, p_\beta H_\alpha, 
\frac{p_0}{2} H_\alpha,
p_\alpha H_\alpha, (\frac{p_0}{2} + p_\alpha) H_\alpha,
\frac{p_0}{2} H_\alpha).
\nonumber
\end{eqnarray}
where
\begin{eqnarray}
G^{k,\ell}_{z}(c_1,c_2,c_3,c_4,c_5,c_6)
&=&
(\frac{1}{k!})^2 (\frac{1}{\ell!})^2
\int \limits _{(\Bbb R ^2)^{k_\alpha -1}} d{\bf u}_{k_\alpha-1}
\int \limits _{(\Bbb R ^2)^{k_\beta -1 }} d{\bf v}_{k_\beta-1}
\nonumber
\\
&&\times
J_{z}({\bf u}_{k_\alpha},{\bf v}_{k_\beta} )
\gamma(c_1,c_2,c_3)({\bf u}_{k_\alpha})
\gamma(c_4,c_5,c_6)({\bf v}_{k_\beta}),
\nonumber
\end{eqnarray}


>From Remark 4.2 we see that the asymptotic shape of the cluster is given by
\begin{equation*}
\{ x\in {\Bbb R^2} : 
|h_\alpha (x)| \le 0, \; |h_\beta (x) | \le 0 \} = \{ {\bf 0} \}.
\end{equation*}
By Lemma \ref{lemma4.3} {\rm (iii)}, Lemma \ref{lemma4.4} {\rm (iii)} 
and a simple calculation we have
\begin{eqnarray*}
\triangle ({\bf x}_k, {\bf y}_\ell | z)
&\sim &
\frac{1}{C_{\alpha, \beta}}
\{|B^{0, \alpha}_{\frac{1}{2}R_\alpha^0, \frac{1}{2}R_\alpha}({\bf x}_k)
\setminus
B^{0, \alpha}_{\frac{1}{2}R_\alpha^0, \frac{1}{2}R_\alpha}|
+|B^{0,\beta}_{\frac{1}{2}R_\beta^0, \frac{1}{2}R_\beta}({\bf y}_\ell)
\setminus
B^{0,\beta}_{\frac{1}{2}R_\beta^0, \frac{1}{2}R_\beta}|\}
\\
&+& (2H_0 - H_\beta)
\{ M(\overline{h}_0 ({\bf x}_k\times ({\bf y}_\ell + z))) - \overline{h}_0 (z)\}
\end{eqnarray*}
when 
$|h_\alpha (z)|= o(1)$, $|h_\beta(z)|= o(1)$.
>From Lemma \ref{lemma4.2} (ii)
\begin{equation*}
\triangle (z) 
= \frac{1}{2}(h_\alpha(z)^2 + h_\beta(z)^2)
+ (2H_0 - H_\beta)|\overline{h}_0 (z)|,
\end{equation*}
if $|h_\alpha(z)|\le H_\alpha$, $|h_\beta(z)|\le H_\beta$.
Then
\begin{eqnarray*}
f_\lambda (0, k_\alpha, k_\beta)
&\sim&
\int \limits _{(\Bbb R ^2)^{k_\alpha -1}} d{\bf u}_{k_\alpha-1}
\int \limits _{(\Bbb R ^2)^{k_\beta-1 }} d{\bf v}_{k_\beta-1}
K_\lambda ({\bf u}_{k_\alpha},{\bf v}_{k_\beta})
\\
&&\quad\times
\chi^{0,\alpha,\beta}_{\lambda p_0 }({\bf u}_{k_\alpha},{\bf v}_{k_\beta})
\chi^{\alpha,\beta}_{\lambda p_\alpha} ({\bf v}_{k_\beta})
\chi^{\alpha,\beta}_{\lambda p_\beta}({\bf u}_{k_\alpha}),
\end{eqnarray*}
where
\begin{eqnarray*}
K_\lambda ({\bf u}_{k_\alpha},{\bf v}_{k_\beta})
&=& \int_{{\Bbb R}^2}dz 
\exp [-\frac{1}{2} C_{\alpha,\beta}p_0
\lambda (h_\alpha(z)^2 + h_\beta(z)^2)]
\\
&&\times 
\exp [-\lambda C_{\alpha,\beta} p_0 (2H_0 - H_\beta)
M(\overline{h}_0 ({\bf u}_{k_\alpha}\cdot ({\bf v}_{k_\beta}+ z)))].
\end{eqnarray*}
By Lemma \ref{lemma3.1} and Lemma \ref{lemma4.3} {\rm (iii)} and 
similar calculations as above, we have
\begin{eqnarray*}
f_\lambda (0, k_\alpha, k_\beta)
&\sim&
\left(\frac{1}{4C_{\alpha, \beta}\lambda^2}\right)^{|k|-1}
\left(\frac{8\pi C_{\alpha, \beta}\lambda}{p_0} \right)^{\frac{1}{2}}
\\
&\times&
\int \limits _{(\Bbb R ^2)^{k_\alpha -1}} d{\bf u}_{k_\alpha-1}
\int \limits _{(\Bbb R ^2)^{k_\beta -1 }} d{\bf v}_{k_\beta-1}
J_{\frac{p_0}{2}(2H_0 - H_\beta)}({\bf u}_{k_\alpha},{\bf v}_{k_\beta} )
\nonumber
\\
&\times&
\gamma^{k_\alpha}
((\frac{1}{2} p_0 + p_\beta) H_\alpha, p_\beta H_\alpha, 
\frac{1}{2}p_0 H_\alpha)
\gamma^{k_\beta}
(p_\alpha H_\alpha, (\frac{1}{2}p_0  + p_\alpha) H_\alpha,
\frac{1}{2}p_0 H_\alpha).
\end{eqnarray*}
Since by Lemma \ref{lemma4.1} (ii),
$\Phi ({\bf p}) = 4C_{\alpha,\beta}
\{p_0 H_0 H_\alpha  (1-\frac{3}{4}p_0)H_\alpha^2 \}$,
we have (\ref{case4}) from (\ref{asym1}) and the above estimates.

\vspace{.5cm}

\noindent {\bf Proof of Theorem 2.2 }
First we examine the behaviour of the function 
$\mu_{\lambda\rho}(C_{\bf 0}\in\Lambda({\bf k})| \Gamma_0)$
as $\lambda\to\infty$ when ${\bf k}=(k_0,k_\alpha,k_\beta)$,
with $k_0, k_\alpha, k_\beta \in {\Bbb N}$.
>From (1.3) and an argument similar to that needed to obtain (4.1) we have
\begin{equation}
\mu_{\lambda \rho}(C_{\bf 0}\in\Lambda({\bf k})\mid ~  \Gamma_{\bf 0})
= \lambda^{|{\bf k}| -1}|{\bf k}|
\frac{p_0^{k_0} p_\alpha^{k_\alpha} p_\beta^{k_\beta}}{k_0!k_\alpha!k_\beta!}
F_\lambda ({\bf k}),
\end{equation}
where
\begin{eqnarray*}
F_\lambda ({\bf k})
&=&e^{-\lambda \{ 
p_0 |B^{\alpha, 0}_{R_\alpha,R_0}\cup B^{\beta, 0}_{R_\beta ,R_0}|
+p_\alpha |B^{0,\alpha}_{R_0,R_\alpha}\cup B^{\beta,\alpha}_{R_\beta,R_\alpha}|
+ p_\beta |B^{0,\beta}_{R_0,R_\beta}\cup B^{\alpha,\beta}_{R_\alpha,R_\beta}|
\}}
\\
&&\qquad\qquad\times
\int \limits _{(\Bbb R ^2)^{k_0-1} } d{\bf x}_{k_0-1}
\int \limits _{(\Bbb R ^2)^{k_\alpha}} d{\bf y}_{k_\alpha}
\int \limits _{(\Bbb R ^2)^{k_\beta }} d{\bf z}_{k_\beta}
1_{\Lambda({\bf k})} 
(C_{\bf 0} ( {\bf x}_{k_0},{\bf y}_{k_\alpha},{\bf z}_{k_\beta } ))
\\
&&\qquad\qquad\times
\chi^{0,\alpha,\beta}_{\lambda p_0 }({\bf y}_{k_\alpha},{\bf z}_{k_\beta})
\chi^{\alpha,\beta,0}_{\lambda p_\alpha} ({\bf z}_{k_\beta},{\bf x}_{k_0})
\chi^{\beta,0,\alpha}_{\lambda p_\beta}({\bf x}_{k_0}, {\bf y}_{k_\alpha}).
\end{eqnarray*}
>From the above we see that the probability that the cluster contains sticks of three distinct orientations is much smaller than that of only two distinct orientations.

For  {\bf case (1)}, when $a,b \ge 2$, from (\ref{case1}), (\ref{case3}) and (\ref{case2}) we have
\begin{eqnarray*}
&&\lim_{\lambda\to\infty} \frac{-1}{4C_{\alpha,\beta}\lambda}
\log \mu_{\lambda\rho}(C_0 \in \Lambda(0,k,\ell) | \Gamma_0)
= p_0(a+b-1)+(1-p_0)ab,
\\
&&\lim_{\lambda\to\infty} \frac{-1}{4C_{\alpha,\beta}\lambda}
\log \mu_{\lambda\rho}(C_0 \in \Lambda(k,0,\ell) | \Gamma_0)
=p_\alpha ab +\frac{p_\alpha}{4} + (1-p_\alpha)b,
\\
&&\lim_{\lambda\to\infty} \frac{-1}{4C_{\alpha,\beta}\lambda}
\log \mu_{\lambda\rho}(C_0 \in \Lambda(k,\ell,0) | \Gamma_0)
=p_\beta ab +\frac{p_\beta}{4} + (1-p_\beta)a.
\end{eqnarray*}
Since 
$$
p_0(a+b-1)+(1-p_0)ab > \min\{p_\alpha ab +\frac{p_\alpha}{4} 
+ (1-p_\alpha)b,p_\beta ab +\frac{p_\beta}{4} + (1-p_\beta)a\},
$$
we obtain Theorem \ref{pshape} (1) (i) and (ii).
>From (\ref{case2}) we see that 
$$
\mu_{\lambda\rho}(C_0 \in \Lambda(k,0,\ell) | \Gamma_0)
\exp\{\lambda\Phi({\bf p})\} \sim c \lambda^{k+\ell-5/2},
$$
and 
$$
\mu_{\lambda\rho}(C_0 \in \Lambda(k,\ell,0) | \Gamma_0)
\exp\{\lambda\Phi({\bf p})\} \sim c' \lambda^{k+\ell-5/2},
$$
with positive constants $c$ and $c'$ independent of $\lambda$.
Thus we have (iii).

For {\bf case (2)}, when $1/2 < \min\{a,b\} <2$, $a\not= b$, $a,b \not= 1$,
from (\ref{case2}) we have
\begin{eqnarray*}
&&\lim_{\lambda\to\infty} \frac{-1}{4C_{\alpha,\beta}\lambda}
\log \mu_{\lambda\rho}(C_0 \in \Lambda(0,k,\ell) | \Gamma_0)
= f(0,\alpha,\beta),
\\
&&\lim_{\lambda\to\infty} \frac{-1}{4C_{\alpha,\beta}\lambda}
\log \mu_{\lambda\rho}(C_0 \in \Lambda(k,0,\ell) | \Gamma_0)
=f(\beta,0,\alpha)
\\
&&\lim_{\lambda\to\infty} \frac{-1}{4C_{\alpha,\beta}\lambda}
\log \mu_{\lambda\rho}(C_0 \in \Lambda(k,\ell,0) | \Gamma_0)
=f(\alpha,\beta,0).
\end{eqnarray*}
Thus we obtain Theorem \ref{pshape} (2).

For {\bf case (3)}, when $0 <a=b<1$, from (\ref{case2}) and (\ref{case3})
we have
\begin{eqnarray*}
&&\lim_{\lambda\to\infty} \frac{-1}{4C_{\alpha,\beta}\lambda}
\log \mu_{\lambda\rho}(C_0 \in \Lambda(0,k,\ell) | \Gamma_0)
= p_0a + (1-\frac{3}{4}p_0)a^2,
\\
&&\lim_{\lambda\to\infty} \frac{-1}{4C_{\alpha,\beta}\lambda}
\log \mu_{\lambda\rho}(C_0 \in \Lambda(k,0,\ell) | \Gamma_0)
=\frac{1}{4}p_\alpha a^2 + a,
\\
&&\lim_{\lambda\to\infty} \frac{-1}{4C_{\alpha,\beta}\lambda}
\log \mu_{\lambda\rho}(C_0 \in \Lambda(k,\ell,0) | \Gamma_0)
=\frac{1}{4}p_\beta a^2 + a.
\end{eqnarray*}
If $p_\alpha \ge p_\beta$, $A(\alpha,\beta)$ occurs
whenever 
$$
p_0a + (1-\frac{3}{4}p_0)a^2 < \frac{1}{4}p_\beta a^2 + a,
$$
i.e., $a < {\bf l}_1(p_0,p_\alpha,p_\beta)$.
Since ${\bf l}_1(p_0,p_\alpha,p_\beta)\ge 1$ for $p_0 \le p_\beta$,
we obtain Theorem \ref{pshape} (3).

Finally for {\bf case (4)}, when $1 <a=b<2$, from (\ref{case2}) and (\ref{case3})
we have
\begin{eqnarray*}
&&\lim_{\lambda\to\infty} \frac{-1}{4C_{\alpha,\beta}\lambda}
\log \mu_{\lambda\rho}(C_0 \in \Lambda(0,k,\ell) | \Gamma_0)
= p_0a + (1-\frac{3}{4}p_0)a^2,
\\
&&\lim_{\lambda\to\infty} \frac{-1}{4C_{\alpha,\beta}\lambda}
\log \mu_{\lambda\rho}(C_0 \in \Lambda(k,0,\ell) | \Gamma_0)
=p_\alpha a^2 + \frac{1}{4}p_\alpha  + (1-p_\alpha)a,
\\
&&\lim_{\lambda\to\infty} \frac{-1}{4C_{\alpha,\beta}\lambda}
\log \mu_{\lambda\rho}(C_0 \in \Lambda(k,\ell,0) | \Gamma_0)
=p_\beta a^2 + \frac{1}{4}p_\beta  + (1-p_\beta)a.
\end{eqnarray*}
If $p_\alpha \ge p_\beta$, we see that $A(\alpha,\beta)$ occurs
whenever 
$$
p_0a+(1-\frac{3}{4}p_0)a^2 < p_\beta a^2 + \frac{1}{4}p_\beta +(1-p_\beta)a,
$$
i.e.,  $a < {\bf l}_2(p_0,p_\alpha,p_\beta)$.
Since ${\bf l}_2(p_0,p_\alpha,p_\beta)\le 1$ for $p_0 \ge p_\beta$,
we obtain Theorem \ref{pshape} (4).

Also for {\bf case (4)} $a=b=1$, from (\ref{case2}) and (\ref{case3})
we have Theorem \ref{pshape} (5), easily.

\section{Appendix}

\noindent{\bf Proof of Lemma \ref{lemma3.1}}: We bound the volume of $B^{\alpha,\beta}_{R_\alpha, R_\beta} ({\bf x}_k)$
by the volume of the smallest parallelogram containing it.
\begin{eqnarray*}
|B^{\alpha,\beta}_{R_\alpha, R_\beta} ({\bf x}_k)|
& \leq & (2R_\alpha + M({\bf x}_k^\alpha))(2R_\beta + M({\bf x}_k^\beta))\sin(\beta - \alpha)\\
&=& 2R_\alpha 2R_\beta \sin(\beta - \alpha) + (2R_\alpha M({\bf x}_k^\beta)  + 2R_\beta M({\bf x}_k^\alpha))\sin(\beta - \alpha)\\
& & \qquad + M({\bf x}_k^\alpha)M({\bf x}_k^\beta)\sin(\beta - \alpha)\\
&= & |B^{\alpha,\beta}_{R_\alpha, R_\beta}| +
2C_{\alpha,\beta}
\{H_\alpha M (h_\beta ({\bf x}_k))+ H_\beta M(h_\alpha({\bf x}_k))\}
\\
& & \qquad +  C_{\alpha,\beta}
M(h_\beta({\bf x}_k))M(h_\alpha({\bf x}_k))
\end{eqnarray*}
which yields (\ref{lem33_i}).

The inequality (\ref{lem33_iia}) follows on observing that \\
(i) $|B^{\alpha,\beta}_{R_\alpha, R_\beta} ({\bf x}_k)|$ must include an area $2R_\alpha \max\{x_1^\beta, \ldots , x_k^\beta\} \sin(\beta - \alpha)$ along the `length' of the connected cluster,\\
(ii) $|B^{\alpha,\beta}_{R_\alpha, R_\beta} ({\bf x}_k)|$ must must include an area $2R_\beta \max\{x_1^\alpha, \ldots , x_k^\alpha\} \sin(\beta - \alpha)$ along the `breadth' of the connected cluster.\\
Thus removing the double counting obtained when we consider the parallelograms along the breadth of the cluster we obtain (\ref{lem33_iia}). 

To show the last inequality we must estimate the double counting more precisely. Observe that
the two halves of the parallelograms on the extremes (in either of the two directions $\alpha$ or $\beta$) of the region $B^{\alpha,\beta}_{R_\alpha, R_\beta} ({\bf x}_k)$ constitute an area $|B^{\alpha,\beta}_{R_\alpha, R_\beta}|$. Also if $B^{\alpha,\beta}_{R_\alpha, R_\beta} ({\bf x}_k)$ is connected, then the area of this region between the lines $\{x \in \mathbb R^2: x^\alpha = \min\{x_1^\alpha, \ldots , x_k^\alpha\}\}$ and $\{x \in \mathbb R^2: x^\alpha = \max\{x_1^\alpha, \ldots , x_k^\alpha\}\}$  has an area at least $(2R_\alpha \max\{x_1^\beta, \ldots , x_k^\beta\} + 2R_\beta\max\{x_1^\alpha, \ldots , x_k^\alpha\}) \sin(\beta - \alpha) -  \max\{x_1^\alpha, \ldots , x_k^\alpha\}\max\{x_1^\beta, \ldots , x_k^\beta\} \sin(\beta - \alpha)$. Since $x_k = {\bf 0}$, (\ref{lem33_ii}) follows.
$\qquad$ $\qed$

\noindent {\bf Proof of Lemma \ref{lemma4.1}} 
If $2H_0 \geq H_\beta$ and 
$H_\alpha \geq H_\beta$. Then  
\begin{eqnarray*}
|B_{R_0, R_{\alpha}}^{0,\alpha} \cup B_{R_0, R_{\beta}}^{0,\beta}|
& = &
|B_{R_0, R_{\alpha}}^{0,\alpha} \setminus B_{R_0, R_{\beta}}^{0,\beta}| +
|B_{R_0, R_{\beta}}^{0,\beta} \setminus B_{R_0, R_{\alpha}}^{0,\alpha}| + 
|B_{R_0, R_{\alpha}}^{0,\alpha} \cap B_{R_0, R_{\beta}}^{0,\beta}|\\
& = &  2 R_0. 2R_\alpha \sin \alpha + R_\beta \sin (\pi - \beta) R_\beta \sin 
(\beta - \alpha)(\sin\alpha)^{-1}\\
& = & C_{\alpha, \beta} (4 H_0 H_\alpha + H_\beta^2).
\end{eqnarray*}
If $2H_0 \geq H_\alpha$ and $H_\beta
 \geq H_\alpha$. Then, similarly, we have 
\begin{eqnarray*}
|B_{R_0, R_{\alpha}}^{0,\alpha} \cup B_{R_0, R_{\beta}}^{0,\beta}| & = &  2 
R_0. 2R_\beta \sin \beta + R_\alpha\sin (\pi - \alpha) R_\alpha \sin (\beta - 
\alpha)(\sin\beta)^{-1}\\
& = & C_{\alpha, \beta} (4 H_0 H_\beta + H_\alpha^2).
\end{eqnarray*} 

Finally if $H_\alpha, H_\beta > 2H_0$, then
\begin{eqnarray*}
|B_{R_0, R_{\alpha}}^{0,\alpha} \cup B_{R_0, R_{\beta}}^{0,\beta}|
& = &|B_{R_0, R_{\alpha}}^{0,\alpha}| +
|B_{R_0, R_{\beta}}^{0,\beta}| - |B_{R_0, R_{\alpha}}^{0,\alpha} \cap B_{R_0, 
R_{\beta}}^{0,\beta}|\\
& = & 4R_0 R_\alpha \sin\alpha + 4 R_0 R_\beta \sin \beta - 4R_0^2\sin  \alpha 
\sin\beta (\sin(\beta-\alpha))^{-1}\\
&=& 4C_{\alpha, \beta} H_0(H_\alpha + H_\beta - H_0). 
\end{eqnarray*} 
This proves the lemma. $\qquad$ $\qed$

\noindent {\bf Proof of Lemma \ref{lemma4.2}} 
Suppose that $2H_0 \geq H_\beta$ and $H_\alpha \geq H_\beta$. 
Also assume that $|\overline{h}_0 (x)| \le 
H_\alpha -H_\beta$ and $|h_\beta (x)| \le 2H_0 - H_\beta$. 
In this case we have
$B_{R_0, R_{\alpha}}^{0,\alpha} \cup B_{R_0, R_{\beta}}^{0,\beta}$ 
represented as the union of the two parallelograms $ABCD$ and $EFGH$ 
in Figure 5, while 
$B_{R_0, R_{\alpha}}^{0,\alpha} \cup B_{R_0,R_{\beta}}^{0,\beta}(x)$ 
is the union of $ABCD$ and $IJKL$. 
The difference between these two regions is 
thus the difference of the ``dashed" triangles and 
the ``solid" triangles outside the parallelogram $ABCD$. 
It is easily seen that the sum of the area of the 
``dashed" triangles is 
$\frac{\sin \alpha \sin \beta}{\sin(\beta - \alpha)} 
[ \frac{R_\beta^2 \sin^2(\beta - \alpha)}{\sin^2\alpha} 
+ (x_1 - \frac{x_2}{\tan \alpha})]$, 
while the sum of the 
areas of the solid triangles is $\frac{R_\beta^2 \sin \beta \sin (\beta - 
\alpha)}{\sin\alpha}$. 
This proves the first case Lemma \ref{lemma4.2} (i).
By considering similar figures, the other parts of the lemma follow. $\qquad$ $\qed$


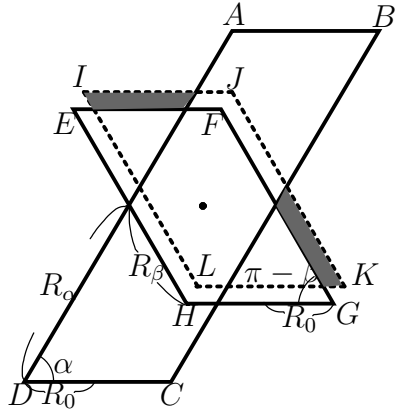
\begin{figure}[h]
\label{Fig5}
\begin{center}
\unitlength 0.1in
\begin{picture}( 22.1700, 20.9600)(  0.0400,-21.8800)
%
\special{pn 20}%
\special{sh 1}%
\special{ar 1300 1210 10 10 0  6.28318530717959E+0000}%
\special{sh 1}%
\special{ar 1300 1210 10 10 0  6.28318530717959E+0000}%
%
\special{pn 20}%
\special{pa 366 2132}%
\special{pa 1134 2132}%
\special{pa 2222 294}%
\special{pa 1454 294}%
\special{pa 1454 294}%
\special{pa 366 2132}%
\special{fp}%
%
\special{pn 8}%
\special{ar 372 2124 148 148  5.2539566 6.2831853}%
\put(5.4400,-16.3200){\makebox(0,0){$R_{\alpha}$}}%
\put(5.7000,-20.4800){\makebox(0,0){$\alpha$}}%
\put(17.1800,-15.6900){\makebox(0,0){$\pi-\beta$}}%
%
\special{pn 8}%
\special{pa 710 1382}%
\special{pa 728 1356}%
\special{pa 746 1330}%
\special{pa 766 1304}%
\special{pa 786 1280}%
\special{pa 808 1256}%
\special{pa 832 1236}%
\special{pa 856 1216}%
\special{pa 886 1204}%
\special{pa 908 1208}%
\special{sp}%
%
\special{pn 8}%
\special{pa 416 1876}%
\special{pa 402 1904}%
\special{pa 388 1934}%
\special{pa 376 1962}%
\special{pa 366 1992}%
\special{pa 358 2024}%
\special{pa 350 2054}%
\special{pa 348 2086}%
\special{pa 352 2118}%
\special{pa 366 2134}%
\special{sp}%
%
\special{pn 20}%
\special{pa 622 704}%
\special{pa 1396 704}%
\special{pa 1984 1722}%
\special{pa 1210 1722}%
\special{pa 1216 1722}%
\special{pa 622 704}%
\special{fp}%
\put(10.0500,-15.4200){\makebox(0,0){$R_{\beta}$}}%
%
\special{pn 8}%
\special{pa 1620 1716}%
\special{pa 1638 1742}%
\special{pa 1664 1760}%
\special{pa 1694 1770}%
\special{pa 1728 1772}%
\special{pa 1728 1772}%
\special{sp}%
%
\special{pn 8}%
\special{pa 366 2132}%
\special{pa 384 2158}%
\special{pa 410 2176}%
\special{pa 440 2186}%
\special{pa 474 2188}%
\special{pa 474 2188}%
\special{sp}%
%
\special{pn 8}%
\special{pa 1984 1728}%
\special{pa 1966 1754}%
\special{pa 1940 1774}%
\special{pa 1910 1782}%
\special{pa 1876 1786}%
\special{pa 1876 1786}%
\special{sp}%
%
\special{pn 8}%
\special{pa 736 2132}%
\special{pa 718 2158}%
\special{pa 692 2176}%
\special{pa 662 2186}%
\special{pa 628 2188}%
\special{pa 628 2188}%
\special{sp}%
\put(18.0500,-17.9800){\makebox(0,0){$R_0$}}%
\put(5.4400,-22.0100){\makebox(0,0){$R_0$}}%
%
\special{pn 20}%
\special{pa 678 614}%
\special{pa 1454 614}%
\special{pa 2042 1632}%
\special{pa 1268 1632}%
\special{pa 1274 1632}%
\special{pa 678 614}%
\special{dt 0.054}%
%
\special{pn 8}%
\special{pa 1204 1728}%
\special{pa 1176 1712}%
\special{pa 1150 1694}%
\special{pa 1124 1674}%
\special{pa 1100 1652}%
\special{pa 1078 1630}%
\special{pa 1062 1612}%
\special{sp}%
%
\special{pn 8}%
\special{pa 916 1210}%
\special{pa 916 1242}%
\special{pa 916 1274}%
\special{pa 918 1306}%
\special{pa 924 1338}%
\special{pa 932 1368}%
\special{pa 940 1398}%
\special{pa 952 1428}%
\special{pa 960 1452}%
\special{sp}%
\put(14.0200,-2.6200){\makebox(0,0)[lb]{$A$}}%
\put(21.8200,-2.7500){\makebox(0,0)[lb]{$B$}}%
\put(10.7500,-22.5200){\makebox(0,0)[lb]{$C$}}%
\put(2.7500,-22.5200){\makebox(0,0)[lb]{$D$}}%
\put(5.1200,-8.3200){\makebox(0,0)[lb]{$E$}}%
\put(12.8600,-8.3800){\makebox(0,0)[lb]{$F$}}%
\put(19.9700,-18.3000){\makebox(0,0)[lb]{$G$}}%
\put(11.3300,-18.3600){\makebox(0,0)[lb]{$H$}}%
\put(6.1400,-5.8800){\makebox(0,0)[lb]{$I$}}%
\put(14.1400,-5.9500){\makebox(0,0)[lb]{$J$}}%
\put(20.7400,-16.2500){\makebox(0,0)[lb]{$K$}}%
\put(12.6200,-15.7700){\makebox(0,0)[lb]{$L$}}%
%
\special{pn 4}%
\special{sh 0.600}%
\special{pa 678 618}%
\special{pa 1262 618}%
\special{pa 1198 698}%
\special{pa 726 706}%
\special{pa 678 618}%
\special{ip}%
%
\special{pn 4}%
\special{sh 0.600}%
\special{pa 1734 1114}%
\special{pa 1686 1210}%
\special{pa 1934 1634}%
\special{pa 2030 1626}%
\special{pa 1734 1114}%
\special{ip}%
%
\special{pn 8}%
\special{ar 1990 1738 190 190  3.2659476 4.1123388}%
\end{picture}%
\end{center}
\caption{\it{Figure accompanying proof of Lemma \ref{lemma4.2}.}}
\end{figure}

\noindent {\bf Proof of Lemma \ref{lemma4.3}} 
First we consider the situation when $y_1 = {\bf 0}$, $\ell =1$ and $k=2$ with $x_2={\bf 0}$ and $x_1$ such that 
\begin{equation}
\label{x1cond}
|x_1^\alpha | \leq R_\alpha - 2R_0^\alpha,
\quad |x_1^\beta | \leq R_\beta - 2R_0^\beta.
\end{equation}

We note here that this choice of $x_1$ ensures the existence of the hatched region in Figure 6 which is isomorphic to a parallelogram with sides making angles $\alpha$ and $\beta$ with the $x$-axis.

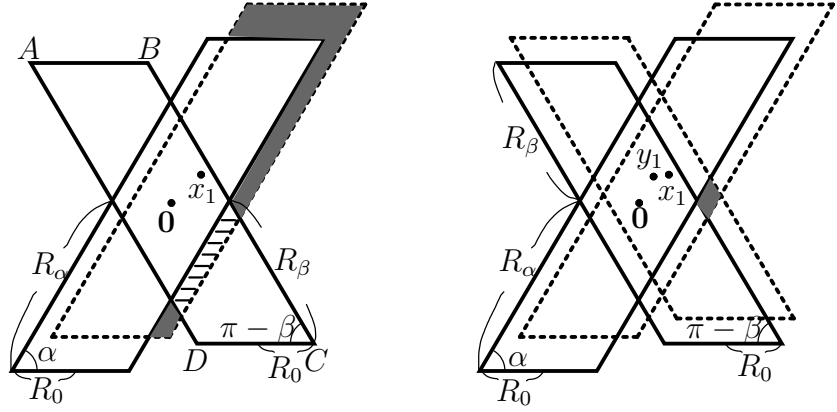
\begin{figure}[h]
\label{Fig6}
\begin{center}\leavevmode
\unitlength 0.1in
\begin{picture}( 46.3700, 19.7100)( -0.0800,-21.6100)
%
\special{pn 20}%
\special{sh 1}%
\special{ar 1168 1230 10 10 0  6.28318530717959E+0000}%
\special{sh 1}%
\special{ar 1168 1230 10 10 0  6.28318530717959E+0000}%
%
\special{pn 20}%
\special{pa 1914 1968}%
\special{pa 1300 1968}%
\special{pa 430 498}%
\special{pa 1044 498}%
\special{pa 1044 498}%
\special{pa 1914 1968}%
\special{fp}%
%
\special{pn 8}%
\special{ar 348 2106 118 118  5.2513655 6.2831853}%
%
\special{pn 8}%
\special{ar 1910 1962 120 120  3.1415927 4.1662224}%
\put(5.3200,-15.5800){\makebox(0,0){$R_{\alpha}$}}%
\put(5.1600,-20.2200){\makebox(0,0){$\alpha$}}%
\put(16.2000,-18.9000){\makebox(0,0){$\pi-\beta$}}%
%
\special{pn 8}%
\special{pa 1874 1762}%
\special{pa 1888 1792}%
\special{pa 1900 1820}%
\special{pa 1912 1850}%
\special{pa 1922 1882}%
\special{pa 1926 1912}%
\special{pa 1928 1944}%
\special{pa 1914 1968}%
\special{sp}%
\put(18.0000,-15.5000){\makebox(0,0){$R_{\beta}$}}%
%
\special{pn 8}%
\special{pa 666 2116}%
\special{pa 646 2142}%
\special{pa 620 2156}%
\special{pa 586 2160}%
\special{pa 578 2162}%
\special{sp}%
%
\special{pn 8}%
\special{pa 1914 1968}%
\special{pa 1894 1994}%
\special{pa 1868 2008}%
\special{pa 1834 2012}%
\special{pa 1828 2012}%
\special{sp}%
%
\special{pn 8}%
\special{pa 342 2110}%
\special{pa 362 2136}%
\special{pa 390 2150}%
\special{pa 422 2156}%
\special{pa 430 2156}%
\special{sp}%
%
\special{pn 8}%
\special{pa 1618 1968}%
\special{pa 1638 1994}%
\special{pa 1666 2008}%
\special{pa 1698 2012}%
\special{pa 1706 2012}%
\special{sp}%
\put(5.1000,-22.1000){\makebox(0,0){$R_0$}}%
\put(17.6000,-21.0000){\makebox(0,0){$R_0$}}%
%
\special{pn 20}%
\special{pa 1970 370}%
\special{pa 1356 370}%
\special{pa 332 2110}%
\special{pa 946 2110}%
\special{pa 1970 370}%
\special{fp}%
%
\special{pn 20}%
\special{pa 2182 190}%
\special{pa 1566 190}%
\special{pa 542 1932}%
\special{pa 1158 1932}%
\special{pa 2182 190}%
\special{dt 0.054}%
%
\special{pn 8}%
\special{pa 328 2106}%
\special{pa 326 2074}%
\special{pa 326 2042}%
\special{pa 328 2010}%
\special{pa 334 1978}%
\special{pa 342 1948}%
\special{pa 354 1918}%
\special{pa 366 1888}%
\special{pa 380 1858}%
\special{pa 394 1828}%
\special{pa 408 1800}%
\special{pa 424 1772}%
\special{pa 438 1742}%
\special{pa 452 1714}%
\special{pa 460 1696}%
\special{sp}%
%
\special{pn 8}%
\special{pa 854 1220}%
\special{pa 824 1232}%
\special{pa 796 1246}%
\special{pa 768 1260}%
\special{pa 742 1278}%
\special{pa 718 1300}%
\special{pa 696 1322}%
\special{pa 676 1346}%
\special{pa 658 1374}%
\special{pa 640 1402}%
\special{pa 624 1430}%
\special{pa 608 1460}%
\special{pa 604 1466}%
\special{sp}%
\put(4.8000,-4.7200){\makebox(0,0)[rb]{$A$}}%
\put(11.1600,-4.7600){\makebox(0,0)[rb]{$B$}}%
\put(19.9000,-21.1000){\makebox(0,0)[rb]{$C$}}%
\put(13.6000,-21.1000){\makebox(0,0)[rb]{$D$}}%
%
\special{pn 20}%
\special{sh 1}%
\special{ar 1322 1082 10 10 0  6.28318530717959E+0000}%
\special{sh 1}%
\special{ar 1322 1082 10 10 0  6.28318530717959E+0000}%
\put(11.4000,-13.3000){\makebox(0,0){$\bf 0$}}%
\put(13.3000,-11.8000){\makebox(0,0){$x_1$}}%
%
\special{pn 8}%
\special{sh 0.600}%
\special{pa 1564 190}%
\special{pa 2172 198}%
\special{pa 1532 1310}%
\special{pa 1484 1222}%
\special{pa 1972 374}%
\special{pa 1468 358}%
\special{pa 1564 190}%
\special{ip}%
%
\special{pn 8}%
\special{sh 0.600}%
\special{pa 1164 1750}%
\special{pa 1212 1822}%
\special{pa 1164 1934}%
\special{pa 1052 1934}%
\special{pa 1052 1934}%
\special{pa 1164 1750}%
\special{ip}%
%
\special{pn 8}%
\special{pa 1638 1358}%
\special{pa 1620 1332}%
\special{pa 1602 1306}%
\special{pa 1582 1280}%
\special{pa 1560 1256}%
\special{pa 1536 1236}%
\special{pa 1510 1220}%
\special{pa 1482 1218}%
\special{sp}%
%
\special{pn 20}%
\special{sh 1}%
\special{ar 3616 1230 10 10 0  6.28318530717959E+0000}%
\special{sh 1}%
\special{ar 3616 1230 10 10 0  6.28318530717959E+0000}%
%
\special{pn 20}%
\special{pa 4362 1968}%
\special{pa 3748 1968}%
\special{pa 2878 498}%
\special{pa 3492 498}%
\special{pa 3492 498}%
\special{pa 4362 1968}%
\special{fp}%
%
\special{pn 8}%
\special{ar 2796 2106 118 118  5.2513655 6.2831853}%
%
\special{pn 8}%
\special{ar 4358 1962 120 120  3.1415927 4.1662224}%
\put(29.8000,-15.4200){\makebox(0,0){$R_{\alpha}$}}%
\put(29.8000,-20.4600){\makebox(0,0){$\alpha$}}%
\put(40.6000,-18.9400){\makebox(0,0){$\pi-\beta$}}%
%
\special{pn 8}%
\special{pa 3132 1070}%
\special{pa 3150 1098}%
\special{pa 3168 1122}%
\special{pa 3188 1148}%
\special{pa 3210 1172}%
\special{pa 3234 1192}%
\special{pa 3260 1210}%
\special{pa 3290 1210}%
\special{sp}%
%
\special{pn 8}%
\special{pa 2900 702}%
\special{pa 2888 674}%
\special{pa 2874 644}%
\special{pa 2864 614}%
\special{pa 2854 584}%
\special{pa 2848 552}%
\special{pa 2846 522}%
\special{pa 2860 496}%
\special{sp}%
\put(29.9600,-9.1000){\makebox(0,0){$R_{\beta}$}}%
%
\special{pn 8}%
\special{pa 3114 2116}%
\special{pa 3094 2142}%
\special{pa 3068 2156}%
\special{pa 3034 2160}%
\special{pa 3026 2162}%
\special{sp}%
%
\special{pn 8}%
\special{pa 4362 1968}%
\special{pa 4342 1994}%
\special{pa 4316 2008}%
\special{pa 4282 2012}%
\special{pa 4276 2012}%
\special{sp}%
%
\special{pn 8}%
\special{pa 2790 2110}%
\special{pa 2810 2136}%
\special{pa 2838 2150}%
\special{pa 2870 2156}%
\special{pa 2878 2156}%
\special{sp}%
%
\special{pn 8}%
\special{pa 4066 1968}%
\special{pa 4086 1994}%
\special{pa 4114 2008}%
\special{pa 4146 2012}%
\special{pa 4154 2012}%
\special{sp}%
\put(29.6000,-22.2000){\makebox(0,0){$R_0$}}%
\put(42.2000,-20.9000){\makebox(0,0){$R_0$}}%
%
\special{pn 20}%
\special{pa 4418 370}%
\special{pa 3804 370}%
\special{pa 2780 2110}%
\special{pa 3394 2110}%
\special{pa 4418 370}%
\special{fp}%
%
\special{pn 20}%
\special{pa 4630 190}%
\special{pa 4014 190}%
\special{pa 2990 1932}%
\special{pa 3606 1932}%
\special{pa 4630 190}%
\special{dt 0.054}%
%
\special{pn 8}%
\special{pa 2776 2106}%
\special{pa 2774 2074}%
\special{pa 2774 2042}%
\special{pa 2776 2010}%
\special{pa 2782 1978}%
\special{pa 2790 1948}%
\special{pa 2802 1918}%
\special{pa 2814 1888}%
\special{pa 2828 1858}%
\special{pa 2842 1828}%
\special{pa 2856 1800}%
\special{pa 2872 1772}%
\special{pa 2886 1742}%
\special{pa 2900 1714}%
\special{pa 2908 1696}%
\special{sp}%
%
\special{pn 8}%
\special{pa 3302 1220}%
\special{pa 3272 1232}%
\special{pa 3244 1246}%
\special{pa 3216 1260}%
\special{pa 3190 1278}%
\special{pa 3166 1300}%
\special{pa 3144 1322}%
\special{pa 3124 1346}%
\special{pa 3106 1374}%
\special{pa 3088 1402}%
\special{pa 3072 1430}%
\special{pa 3056 1460}%
\special{pa 3052 1466}%
\special{sp}%
%
\special{pn 20}%
\special{sh 1}%
\special{ar 3770 1082 10 10 0  6.28318530717959E+0000}%
\special{sh 1}%
\special{ar 3770 1082 10 10 0  6.28318530717959E+0000}%
\put(36.2000,-13.0600){\makebox(0,0){$\bf 0$}}%
\put(37.9000,-11.9000){\makebox(0,0){$x_1$}}%
%
\special{pn 20}%
\special{pa 4426 1836}%
\special{pa 3812 1836}%
\special{pa 2942 366}%
\special{pa 3556 366}%
\special{pa 3556 366}%
\special{pa 4426 1836}%
\special{dt 0.054}%
%
\special{pn 20}%
\special{sh 1}%
\special{ar 3690 1092 10 10 0  6.28318530717959E+0000}%
\special{sh 1}%
\special{ar 3690 1092 10 10 0  6.28318530717959E+0000}%
\put(36.7000,-10.0000){\makebox(0,0){$y_1$}}%
%
\special{pn 8}%
\special{sh 0.600}%
\special{pa 3988 1102}%
\special{pa 4044 1190}%
\special{pa 3980 1318}%
\special{pa 3916 1238}%
\special{pa 3988 1102}%
\special{ip}%
%
\special{pn 13}%
\special{pa 1360 1560}%
\special{pa 1290 1560}%
\special{fp}%
\special{pa 1330 1620}%
\special{pa 1260 1620}%
\special{fp}%
\special{pa 1290 1680}%
\special{pa 1220 1680}%
\special{fp}%
\special{pa 1260 1740}%
\special{pa 1190 1740}%
\special{fp}%
\special{pa 1400 1500}%
\special{pa 1330 1500}%
\special{fp}%
\special{pa 1430 1440}%
\special{pa 1360 1440}%
\special{fp}%
\special{pa 1470 1380}%
\special{pa 1400 1380}%
\special{fp}%
\special{pa 1500 1320}%
\special{pa 1430 1320}%
\special{fp}%
\end{picture}%
\end{center}
\caption{\it{The two shaded regions in the left figure combine on collapsing the lines AD and BC. The shaded parallelogram in the right figure is double counted.}}
\end{figure}

>From Figure 6 we see that if we collapse the lines $AD$ and $BC$ into one and remove the parallelogram contained between these lines then each of the parallelograms $B_{R_0, R_\alpha}^{0,\alpha}$ and $B_{R_0, R_\alpha}^{0,\alpha}(x_1)$ become isomorphic to $B_{R_0, R_\alpha - R^0_\alpha}^{0,\alpha}$.
Moreover the shaded area which represents $\left( (B_{R_0, R_\alpha}^{0,\alpha}(x_1, x_2) \cup B_{R_0, R_\beta}^{0,\beta}(y_1)) \setminus (B_{R_0, R_\alpha}^{0,\alpha} \cup B_{R_0, R_\beta}^{0,\beta})\right)$
is isomorphic to
$\left(B_{R_0, R_\alpha - R^0_\alpha}^{0,\alpha}(x_1,x_2) \setminus  B_{R_0, R_\alpha - R^0_\alpha}^{0,\alpha}\right)$.

Since $(B_{R_0, R_\alpha}^{0,\alpha} \cup B_{R_0, R_\beta}^{0,\beta})
\subseteq (B_{R_0, R_\alpha}^{0,\alpha}(x_1, x_2) \cup B_{R_0, R_\beta}^{0,\beta})$ and 
$B_{R_0, R_\alpha - R^0_\alpha}^{0,\alpha} (x_1, x_2) \supseteq B_{R_0, R_\alpha - R^0_\alpha}^{0,\alpha}$ we have
\begin{equation}
\label{ll1}
C_{\alpha, \beta} \Delta ({\bf x}_2, y_1) = |B_{R_0, R_\alpha - R^0_\alpha}^{0,\alpha} ({\bf x}_2) \setminus B_{R_0, R_\alpha - R^0_\alpha}^{0,\alpha}|.
\end{equation}

Now observe that a similar result may be obtained when $x_1 = {\bf 0}$, $k=1$ and
$\ell =2$, $y_2 = {\bf 0}$ and $y_1$ such that 
\begin{equation}
\label{y1cond}
|y_1^\alpha | \leq R_\alpha - 2R_0^\alpha,
\quad |y_1^\beta | \leq R_\beta - 2R_0^\beta.
\end{equation}
In this case we obtain 
\begin{equation}
\label{ll2}
C_{\alpha, \beta} \Delta (x_1, {\bf y}_2) = |B_{R_0, R_\beta - R^0_\beta}^{0,\beta} ({\bf y}_2) \setminus B_{R_0, R_\beta - R^0_\beta}^{0,\beta}|.
\end{equation}

In case both $k =2$ and $\ell=2$ with $x_1$ and $y_1$ satisfying (\ref{x1cond}) and (\ref{y1cond}) we see from Figure 6 that if we add the areas obtained in (\ref{ll1}) and (\ref{ll2}) there is double counting of the shaded parallelogram with sides of length $|x_1^\beta|$ and $|y_1^\alpha|$ and area
$|x_1^\alpha||y_1^\beta| \sin(\beta - \alpha)$. Thus we have
$C_{\alpha, \beta} \Delta ({\bf x}_2, {\bf y}_2) = |B_{R_0, R_\alpha - R^0_\alpha}^{0,\alpha} ({\bf x}_2) \setminus B_{R_0, R_\alpha - R^0_\alpha}^{0,\alpha}| + |B_{R_0, R_\beta - R^0_\beta}^{0,\beta} ({\bf y}_2) \setminus B_{R_0, R_\beta - R^0_\beta}^{0,\beta}| - |x_1^\beta||y_1^\alpha| \sin(\beta - \alpha)$. 

In general, for any $k$ and $\ell$, we see that if 
\begin{equation}
\label{ll4}
M({\bf x}_k) \leq R_\alpha - 2R^0_\alpha, \mbox{ and }
M({\bf y}_\ell) \leq R_\beta - 2R^0_\beta
\end{equation}
there will be many such shaded areas which will be double counted. These areas need not be all distinct and the total area of this double counted region is at most $M({\bf x}_k^\beta) M({\bf y}_\ell^\alpha) \sin(\beta - \alpha)$. 
Now note that the condition (\ref{llcond1}) guarantees that (\ref{ll4}) holds. Hence Lemma \ref{lemma4.3} (i) follows.

The remaining parts of the lemmas follow from similar arguments and are explained through Figures 7 and 8. $\qquad$ $\qed$

Lemma \ref{lemma4.4} follows similarly and its proof is omitted.

\begin{figure}[h]
\label{Fig7}
\begin{center}\leavevmode
\unitlength 0.1in
\begin{picture}( 42.6400, 20.0300)(  0.5200,-22.1300)
%
\special{pn 20}%
\special{sh 1}%
\special{ar 1142 1234 10 10 0  6.28318530717959E+0000}%
\special{sh 1}%
\special{ar 1142 1234 10 10 0  6.28318530717959E+0000}%
%
\special{pn 20}%
\special{pa 208 2156}%
\special{pa 976 2156}%
\special{pa 2064 320}%
\special{pa 1296 320}%
\special{pa 1296 320}%
\special{pa 208 2156}%
\special{fp}%
%
\special{pn 8}%
\special{ar 1814 1746 148 148  3.1415927 4.1708214}%
%
\special{pn 8}%
\special{ar 214 2150 148 148  5.2539566 6.2831853}%
\put(4.1200,-20.7200){\makebox(0,0){$\alpha$}}%
\put(14.7000,-16.7000){\makebox(0,0){$\pi-\beta$}}%
%
\special{pn 8}%
\special{pa 1776 1484}%
\special{pa 1790 1514}%
\special{pa 1804 1542}%
\special{pa 1816 1572}%
\special{pa 1826 1602}%
\special{pa 1836 1632}%
\special{pa 1842 1664}%
\special{pa 1844 1696}%
\special{pa 1840 1726}%
\special{pa 1828 1742}%
\special{sp}%
%
\special{pn 8}%
\special{pa 1718 1420}%
\special{pa 1702 1394}%
\special{pa 1686 1366}%
\special{pa 1670 1338}%
\special{pa 1650 1312}%
\special{pa 1630 1288}%
\special{pa 1608 1264}%
\special{pa 1584 1244}%
\special{pa 1556 1230}%
\special{pa 1534 1232}%
\special{sp}%
%
\special{pn 20}%
\special{pa 464 728}%
\special{pa 1238 728}%
\special{pa 1828 1746}%
\special{pa 1052 1746}%
\special{pa 1060 1746}%
\special{pa 464 728}%
\special{fp}%
%
\special{pn 8}%
\special{sh 0}%
\special{pa 1680 1362}%
\special{pa 1764 1362}%
\special{pa 1764 1446}%
\special{pa 1680 1446}%
\special{pa 1680 1362}%
\special{ip}%
%
\special{pn 8}%
\special{sh 0}%
\special{pa 1744 1458}%
\special{pa 1840 1458}%
\special{pa 1840 1568}%
\special{pa 1744 1568}%
\special{pa 1744 1458}%
\special{ip}%
\put(17.5600,-14.5800){\makebox(0,0){$R_{\beta}$}}%
%
\special{pn 8}%
\special{pa 1462 1740}%
\special{pa 1482 1766}%
\special{pa 1506 1786}%
\special{pa 1536 1794}%
\special{pa 1570 1798}%
\special{pa 1572 1798}%
\special{sp}%
%
\special{pn 8}%
\special{pa 208 2156}%
\special{pa 228 2182}%
\special{pa 252 2202}%
\special{pa 282 2210}%
\special{pa 316 2214}%
\special{pa 316 2214}%
\special{sp}%
%
\special{pn 8}%
\special{pa 1814 1740}%
\special{pa 1796 1766}%
\special{pa 1770 1786}%
\special{pa 1740 1794}%
\special{pa 1706 1798}%
\special{pa 1706 1798}%
\special{sp}%
%
\special{pn 8}%
\special{pa 580 2156}%
\special{pa 560 2182}%
\special{pa 536 2202}%
\special{pa 506 2210}%
\special{pa 472 2214}%
\special{pa 470 2214}%
\special{sp}%
\put(16.6000,-18.6000){\makebox(0,0){$R_0$}}%
\put(3.9000,-22.7000){\makebox(0,0){$R_0$}}%
%
\special{pn 20}%
\special{pa 432 524}%
\special{pa 1206 524}%
\special{pa 1796 1542}%
\special{pa 1020 1542}%
\special{pa 1028 1542}%
\special{pa 432 524}%
\special{dt 0.054}%
%
\special{pn 20}%
\special{pa 150 2054}%
\special{pa 918 2054}%
\special{pa 2006 216}%
\special{pa 1238 216}%
\special{pa 1238 216}%
\special{pa 150 2054}%
\special{dt 0.054}%
%
\special{pn 8}%
\special{pa 630 1248}%
\special{pa 752 1248}%
\special{fp}%
%
\special{pn 8}%
\special{pa 996 518}%
\special{pa 1084 658}%
\special{fp}%
%
\special{pn 8}%
\special{sh 0.600}%
\special{pa 1002 524}%
\special{pa 1174 524}%
\special{pa 1092 664}%
\special{pa 1092 664}%
\special{pa 1092 664}%
\special{pa 1092 664}%
\special{pa 1002 524}%
\special{fp}%
%
\special{pn 8}%
\special{sh 0.600}%
\special{pa 694 1138}%
\special{pa 630 1240}%
\special{pa 752 1240}%
\special{pa 694 1138}%
\special{fp}%
\put(6.6000,-4.8000){\makebox(0,0)[lb]{$R_\beta^0$}}%
%
\special{pn 8}%
\special{pa 988 512}%
\special{pa 970 484}%
\special{pa 946 466}%
\special{pa 916 458}%
\special{pa 882 454}%
\special{pa 880 454}%
\special{sp}%
%
\special{pn 8}%
\special{pa 432 512}%
\special{pa 452 484}%
\special{pa 476 466}%
\special{pa 506 458}%
\special{pa 540 454}%
\special{pa 540 454}%
\special{sp}%
%
\special{pn 20}%
\special{sh 1}%
\special{ar 3396 1228 10 10 0  6.28318530717959E+0000}%
\special{sh 1}%
\special{ar 3396 1228 10 10 0  6.28318530717959E+0000}%
%
\special{pn 20}%
\special{pa 2460 2150}%
\special{pa 3228 2150}%
\special{pa 4316 312}%
\special{pa 3548 312}%
\special{pa 3548 312}%
\special{pa 2460 2150}%
\special{fp}%
%
\special{pn 8}%
\special{ar 4068 1740 148 148  3.1415927 4.1708214}%
%
\special{pn 8}%
\special{ar 2468 2144 148 148  5.2539566 6.2831853}%
\put(26.6500,-20.6600){\makebox(0,0){$\alpha$}}%
\put(37.3000,-16.6000){\makebox(0,0){$\pi-\beta$}}%
%
\special{pn 8}%
\special{pa 4028 1478}%
\special{pa 4042 1506}%
\special{pa 4056 1536}%
\special{pa 4068 1564}%
\special{pa 4078 1596}%
\special{pa 4088 1626}%
\special{pa 4094 1656}%
\special{pa 4098 1688}%
\special{pa 4094 1720}%
\special{pa 4080 1736}%
\special{sp}%
%
\special{pn 8}%
\special{pa 3972 1414}%
\special{pa 3956 1386}%
\special{pa 3938 1358}%
\special{pa 3922 1332}%
\special{pa 3902 1306}%
\special{pa 3882 1282}%
\special{pa 3860 1258}%
\special{pa 3836 1238}%
\special{pa 3808 1224}%
\special{pa 3786 1224}%
\special{sp}%
%
\special{pn 20}%
\special{pa 2716 722}%
\special{pa 3492 722}%
\special{pa 4080 1740}%
\special{pa 3306 1740}%
\special{pa 3312 1740}%
\special{pa 2716 722}%
\special{fp}%
%
\special{pn 8}%
\special{sh 0}%
\special{pa 3932 1356}%
\special{pa 4016 1356}%
\special{pa 4016 1440}%
\special{pa 3932 1440}%
\special{pa 3932 1356}%
\special{ip}%
%
\special{pn 8}%
\special{sh 0}%
\special{pa 3996 1452}%
\special{pa 4092 1452}%
\special{pa 4092 1560}%
\special{pa 3996 1560}%
\special{pa 3996 1452}%
\special{ip}%
\put(40.0900,-14.5200){\makebox(0,0){$R_{\beta}$}}%
%
\special{pn 8}%
\special{pa 3716 1734}%
\special{pa 3734 1760}%
\special{pa 3758 1780}%
\special{pa 3788 1788}%
\special{pa 3822 1792}%
\special{pa 3824 1792}%
\special{sp}%
%
\special{pn 8}%
\special{pa 2460 2150}%
\special{pa 2480 2176}%
\special{pa 2504 2196}%
\special{pa 2534 2204}%
\special{pa 2568 2208}%
\special{pa 2570 2208}%
\special{sp}%
%
\special{pn 8}%
\special{pa 4068 1734}%
\special{pa 4048 1760}%
\special{pa 4024 1780}%
\special{pa 3994 1788}%
\special{pa 3960 1792}%
\special{pa 3958 1792}%
\special{sp}%
%
\special{pn 8}%
\special{pa 2832 2150}%
\special{pa 2814 2176}%
\special{pa 2790 2196}%
\special{pa 2760 2204}%
\special{pa 2726 2208}%
\special{pa 2724 2208}%
\special{sp}%
\put(39.0000,-18.5000){\makebox(0,0){$R_0$}}%
\put(26.4000,-22.5000){\makebox(0,0){$R_0$}}%
%
\special{pn 20}%
\special{pa 2684 518}%
\special{pa 3460 518}%
\special{pa 4048 1536}%
\special{pa 3274 1536}%
\special{pa 3280 1536}%
\special{pa 2684 518}%
\special{dt 0.054}%
%
\special{pn 20}%
\special{pa 2404 2048}%
\special{pa 3172 2048}%
\special{pa 4260 210}%
\special{pa 3492 210}%
\special{pa 3492 210}%
\special{pa 2404 2048}%
\special{dt 0.054}%
%
\special{pn 8}%
\special{sh 0.600}%
\special{pa 3312 518}%
\special{pa 3434 518}%
\special{pa 3318 722}%
\special{pa 3190 728}%
\special{pa 3312 518}%
\special{fp}%
\end{picture}%
\end{center}
\caption{\it{The shaded triangles in the left figure give the last two terms in (\ref{r-tricond}), while the shaded parallelogram in the right figure is double counted.}}
\end{figure}
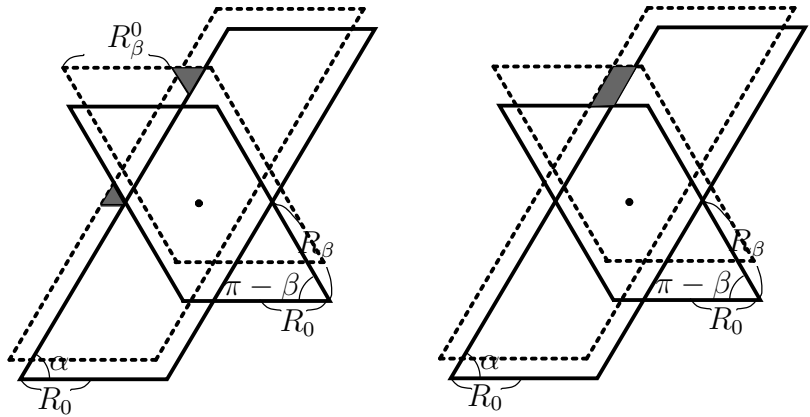
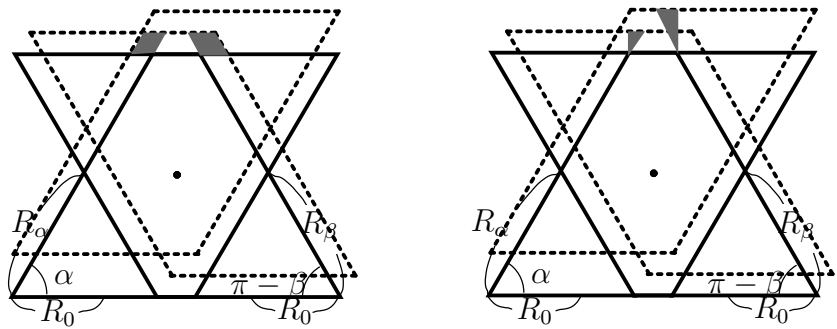
\begin{figure}[h]
\label{Fig8}
\begin{center}\leavevmode
\unitlength 0.1in
\begin{picture}( 47.4000, 15.9200)(  0.9200,-23.1200)
%
\special{pn 20}%
\special{sh 1}%
\special{ar 1400 1592 10 10 0  6.28318530717959E+0000}%
\special{sh 1}%
\special{ar 1400 1592 10 10 0  6.28318530717959E+0000}%
%
\special{pn 8}%
\special{ar 2240 2232 184 184  3.1415927 4.1719695}%
%
\special{pn 8}%
\special{ar 536 2216 184 184  5.2528085 6.2831853}%
\put(8.0800,-21.2800){\makebox(0,0){$\alpha$}}%
\put(18.8000,-21.7600){\makebox(0,0){$\pi-\beta$}}%
%
\special{pn 8}%
\special{pa 2192 1904}%
\special{pa 2206 1934}%
\special{pa 2220 1962}%
\special{pa 2232 1992}%
\special{pa 2244 2022}%
\special{pa 2254 2052}%
\special{pa 2264 2082}%
\special{pa 2272 2112}%
\special{pa 2278 2144}%
\special{pa 2278 2176}%
\special{pa 2272 2206}%
\special{pa 2256 2228}%
\special{sp}%
%
\special{pn 8}%
\special{pa 2120 1824}%
\special{pa 2104 1796}%
\special{pa 2088 1770}%
\special{pa 2070 1742}%
\special{pa 2054 1716}%
\special{pa 2034 1690}%
\special{pa 2014 1666}%
\special{pa 1994 1642}%
\special{pa 1970 1620}%
\special{pa 1946 1600}%
\special{pa 1918 1586}%
\special{pa 1890 1588}%
\special{sp}%
%
\special{pn 8}%
\special{pa 664 1808}%
\special{pa 682 1782}%
\special{pa 700 1756}%
\special{pa 720 1730}%
\special{pa 740 1704}%
\special{pa 760 1680}%
\special{pa 782 1656}%
\special{pa 802 1634}%
\special{pa 828 1614}%
\special{pa 854 1596}%
\special{pa 884 1586}%
\special{pa 910 1590}%
\special{sp}%
%
\special{pn 20}%
\special{pa 552 960}%
\special{pa 1520 960}%
\special{pa 2256 2232}%
\special{pa 1288 2232}%
\special{pa 1296 2232}%
\special{pa 552 960}%
\special{fp}%
%
\special{pn 8}%
\special{sh 0}%
\special{pa 2072 1752}%
\special{pa 2176 1752}%
\special{pa 2176 1856}%
\special{pa 2072 1856}%
\special{pa 2072 1752}%
\special{ip}%
%
\special{pn 8}%
\special{sh 0}%
\special{pa 2152 1872}%
\special{pa 2272 1872}%
\special{pa 2272 2008}%
\special{pa 2152 2008}%
\special{pa 2152 1872}%
\special{ip}%
\put(21.4400,-18.6400){\makebox(0,0){$R_{\beta}$}}%
%
\special{pn 20}%
\special{pa 2240 960}%
\special{pa 1272 960}%
\special{pa 536 2232}%
\special{pa 1504 2232}%
\special{pa 1496 2232}%
\special{pa 2240 960}%
\special{fp}%
%
\special{pn 8}%
\special{pa 600 1902}%
\special{pa 584 1930}%
\special{pa 572 1960}%
\special{pa 560 1990}%
\special{pa 546 2018}%
\special{pa 536 2048}%
\special{pa 526 2078}%
\special{pa 518 2110}%
\special{pa 512 2140}%
\special{pa 510 2172}%
\special{pa 516 2204}%
\special{pa 532 2224}%
\special{sp}%
%
\special{pn 8}%
\special{sh 0}%
\special{pa 536 1864}%
\special{pa 624 1864}%
\special{pa 624 1968}%
\special{pa 536 1968}%
\special{pa 536 1864}%
\special{ip}%
%
\special{pn 13}%
\special{sh 0}%
\special{pa 640 1752}%
\special{pa 688 1752}%
\special{pa 688 1840}%
\special{pa 640 1840}%
\special{pa 640 1752}%
\special{ip}%
\put(6.3200,-18.4800){\makebox(0,0){$R_{\alpha}$}}%
%
\special{pn 8}%
\special{pa 1784 2232}%
\special{pa 1804 2260}%
\special{pa 1826 2282}%
\special{pa 1854 2294}%
\special{pa 1886 2302}%
\special{pa 1912 2304}%
\special{sp}%
%
\special{pn 8}%
\special{pa 528 2232}%
\special{pa 548 2260}%
\special{pa 570 2282}%
\special{pa 598 2294}%
\special{pa 630 2302}%
\special{pa 656 2304}%
\special{sp}%
%
\special{pn 8}%
\special{pa 2240 2240}%
\special{pa 2222 2268}%
\special{pa 2200 2290}%
\special{pa 2172 2302}%
\special{pa 2140 2310}%
\special{pa 2112 2312}%
\special{sp}%
%
\special{pn 8}%
\special{pa 1016 2240}%
\special{pa 998 2268}%
\special{pa 976 2290}%
\special{pa 948 2302}%
\special{pa 916 2310}%
\special{pa 888 2312}%
\special{sp}%
\put(20.1600,-23.0400){\makebox(0,0){$R_0$}}%
\put(7.7600,-23.1200){\makebox(0,0){$R_0$}}%
%
\special{pn 20}%
\special{pa 632 848}%
\special{pa 1600 848}%
\special{pa 2336 2120}%
\special{pa 1368 2120}%
\special{pa 1376 2120}%
\special{pa 632 848}%
\special{dt 0.054}%
%
\special{pn 20}%
\special{pa 2248 736}%
\special{pa 1280 736}%
\special{pa 544 2008}%
\special{pa 1512 2008}%
\special{pa 1504 2008}%
\special{pa 2248 736}%
\special{dt 0.054}%
%
\special{pn 20}%
\special{sh 1}%
\special{ar 3896 1584 10 10 0  6.28318530717959E+0000}%
\special{sh 1}%
\special{ar 3896 1584 10 10 0  6.28318530717959E+0000}%
%
\special{pn 8}%
\special{ar 4736 2224 184 184  3.1415927 4.1719695}%
%
\special{pn 8}%
\special{ar 3032 2208 184 184  5.2528085 6.2831853}%
\put(33.0400,-21.2000){\makebox(0,0){$\alpha$}}%
\put(43.7600,-21.6800){\makebox(0,0){$\pi-\beta$}}%
%
\special{pn 8}%
\special{pa 4688 1896}%
\special{pa 4702 1926}%
\special{pa 4716 1954}%
\special{pa 4728 1984}%
\special{pa 4740 2014}%
\special{pa 4750 2044}%
\special{pa 4760 2074}%
\special{pa 4768 2104}%
\special{pa 4774 2136}%
\special{pa 4774 2168}%
\special{pa 4768 2198}%
\special{pa 4752 2220}%
\special{sp}%
%
\special{pn 8}%
\special{pa 4616 1816}%
\special{pa 4600 1788}%
\special{pa 4584 1762}%
\special{pa 4566 1734}%
\special{pa 4550 1708}%
\special{pa 4530 1682}%
\special{pa 4510 1658}%
\special{pa 4490 1634}%
\special{pa 4466 1612}%
\special{pa 4442 1592}%
\special{pa 4414 1578}%
\special{pa 4386 1580}%
\special{sp}%
%
\special{pn 8}%
\special{pa 3160 1800}%
\special{pa 3178 1774}%
\special{pa 3196 1748}%
\special{pa 3216 1722}%
\special{pa 3236 1696}%
\special{pa 3256 1672}%
\special{pa 3278 1648}%
\special{pa 3298 1626}%
\special{pa 3324 1606}%
\special{pa 3350 1588}%
\special{pa 3380 1578}%
\special{pa 3406 1582}%
\special{sp}%
%
\special{pn 20}%
\special{pa 3048 952}%
\special{pa 4016 952}%
\special{pa 4752 2224}%
\special{pa 3784 2224}%
\special{pa 3792 2224}%
\special{pa 3048 952}%
\special{fp}%
%
\special{pn 8}%
\special{sh 0}%
\special{pa 4568 1744}%
\special{pa 4672 1744}%
\special{pa 4672 1848}%
\special{pa 4568 1848}%
\special{pa 4568 1744}%
\special{ip}%
%
\special{pn 8}%
\special{sh 0}%
\special{pa 4648 1864}%
\special{pa 4768 1864}%
\special{pa 4768 2000}%
\special{pa 4648 2000}%
\special{pa 4648 1864}%
\special{ip}%
\put(46.4000,-18.6400){\makebox(0,0){$R_{\beta}$}}%
%
\special{pn 20}%
\special{pa 4736 952}%
\special{pa 3768 952}%
\special{pa 3032 2224}%
\special{pa 4000 2224}%
\special{pa 3992 2224}%
\special{pa 4736 952}%
\special{fp}%
%
\special{pn 8}%
\special{pa 3096 1894}%
\special{pa 3080 1922}%
\special{pa 3068 1952}%
\special{pa 3056 1982}%
\special{pa 3042 2010}%
\special{pa 3032 2040}%
\special{pa 3022 2070}%
\special{pa 3014 2102}%
\special{pa 3008 2132}%
\special{pa 3006 2164}%
\special{pa 3012 2196}%
\special{pa 3028 2216}%
\special{sp}%
%
\special{pn 8}%
\special{sh 0}%
\special{pa 3032 1856}%
\special{pa 3120 1856}%
\special{pa 3120 1960}%
\special{pa 3032 1960}%
\special{pa 3032 1856}%
\special{ip}%
%
\special{pn 13}%
\special{sh 0}%
\special{pa 3136 1744}%
\special{pa 3184 1744}%
\special{pa 3184 1832}%
\special{pa 3136 1832}%
\special{pa 3136 1744}%
\special{ip}%
\put(30.4000,-18.4000){\makebox(0,0){$R_{\alpha}$}}%
%
\special{pn 8}%
\special{pa 4280 2224}%
\special{pa 4300 2252}%
\special{pa 4322 2274}%
\special{pa 4350 2286}%
\special{pa 4382 2294}%
\special{pa 4408 2296}%
\special{sp}%
%
\special{pn 8}%
\special{pa 3024 2224}%
\special{pa 3044 2252}%
\special{pa 3066 2274}%
\special{pa 3094 2286}%
\special{pa 3126 2294}%
\special{pa 3152 2296}%
\special{sp}%
%
\special{pn 8}%
\special{pa 4736 2232}%
\special{pa 4718 2260}%
\special{pa 4696 2282}%
\special{pa 4668 2294}%
\special{pa 4636 2302}%
\special{pa 4608 2304}%
\special{sp}%
%
\special{pn 8}%
\special{pa 3512 2232}%
\special{pa 3494 2260}%
\special{pa 3472 2282}%
\special{pa 3444 2294}%
\special{pa 3412 2302}%
\special{pa 3384 2304}%
\special{sp}%
\put(45.1200,-22.9600){\makebox(0,0){$R_0$}}%
\put(32.7200,-23.0400){\makebox(0,0){$R_0$}}%
%
\special{pn 20}%
\special{pa 3128 840}%
\special{pa 4096 840}%
\special{pa 4832 2112}%
\special{pa 3864 2112}%
\special{pa 3872 2112}%
\special{pa 3128 840}%
\special{dt 0.054}%
%
\special{pn 20}%
\special{pa 4744 728}%
\special{pa 3776 728}%
\special{pa 3040 2000}%
\special{pa 4008 2000}%
\special{pa 4000 2000}%
\special{pa 4744 728}%
\special{dt 0.054}%
%
\special{pn 8}%
\special{sh 0.600}%
\special{pa 1280 960}%
\special{pa 1344 848}%
\special{pa 1224 848}%
\special{pa 1160 960}%
\special{pa 1288 952}%
\special{pa 1280 960}%
\special{ip}%
%
\special{pn 8}%
\special{sh 0.600}%
\special{pa 1600 848}%
\special{pa 1672 968}%
\special{pa 1520 968}%
\special{pa 1456 848}%
\special{pa 1600 848}%
\special{ip}%
%
\special{pn 8}%
\special{sh 0.600}%
\special{pa 4024 960}%
\special{pa 4024 728}%
\special{pa 3904 720}%
\special{pa 4024 960}%
\special{ip}%
%
\special{pn 8}%
\special{sh 0.600}%
\special{pa 3760 970}%
\special{pa 3760 840}%
\special{pa 3850 840}%
\special{pa 3760 970}%
\special{ip}%
\end{picture}%
\caption{\it{The shaded areas are double counted and is deducted in (\ref{r-trapcond}).}}
\end{center}
\end{figure}




\noindent Rahul Roy\\
Indian Statistical Institute,\\
7 SJS Sansanwal Marg,\\
New Delhi 110016, INDIA.\\
e-mail:- rahul@isid1.isid.ac.in \\

\vspace{.5cm}

\noindent Hideki Tanemura\\
Department of Mathematics and Informatics,\\
Faculty of Science,\\
Chiba University,\\
1-33, Yayoi-cho, Inage-ku,\\
Chiba 263-8522, JAPAN.\\
e-mail:- tanemura@math.s.chiba-u.ac.jp\\

\end{document}